\documentclass[a4paper,11pt]{article}
\usepackage{graphicx}
\usepackage{amssymb}
\usepackage{amsfonts}
\usepackage{amsmath}
\usepackage{undertilde}
\usepackage{color}
\usepackage{url}
\catcode`\@=11 \@addtoreset{equation}{section}

\catcode`\@=12

\newtheorem{Theorem}{Theorem}[section]
\newtheorem{Proposition}{Proposition}[section]
\newtheorem{Lemma}{Lemma}[section]
\newtheorem{Corollary}{Corollary}[section]
\newtheorem{Remark}{Remark}[section]
\newtheorem{Definition}{Definition}[section]

\newcommand{\bTheorem}[1]{
\begin{Theorem} \label{T#1} }
\newcommand{\eT}{\end{Theorem}}

\newcommand{\bProposition}[1]{
\begin{Proposition} \label{P#1}}
\newcommand{\eP}{\end{Proposition}}

\newcommand{\bLemma}[1]{
\begin{Lemma} \label{L#1} }
\newcommand{\eL}{\end{Lemma}}

\newcommand{\bCorollary}[1]{
\begin{Corollary} \label{C#1} }
\newcommand{\eC}{\end{Corollary}}

\newcommand{\bRemark}[1]{
\begin{Remark} \label{R#1} }
\newcommand{\eR}{\end{Remark}}

\newcommand{\bDefinition}[1]{
\begin{Definition} \label{D#1} }
\newcommand{\eD}{\end{Definition}}

\newcommand{\bFormula}[1]{
\begin{equation} \label{#1}}
\newcommand{\eF}{\end{equation}}

\newcommand{\Divh}{{\rm div}_h}
\newcommand{\Gradh}{\nabla_h}

\newcommand{\Ov}[1]{\overline{#1}}
\newcommand{\av}[1]{ \left\{ #1 \right\}}

\newcommand{\DC}{C^\infty_c}

\newcommand{\aleq}{\stackrel{<}{\sim}}

\newcommand{\vr}{\varrho}

\newcommand{\vu}{\vc{u}}
\newcommand{\vc}[1]{{\bf #1}}

\newcommand{\Div}{{\rm div}_x}
\newcommand{\Grad}{\nabla_x}

\newcommand{\tn}[1]{\mbox {\F #1}}
\newcommand{\dx}{{\rm d} {x}}
\newcommand{\dt}{{\rm d} t }

\newcommand{\ju}[1]{[[ #1 ]]}

\newcommand{\intO}[1]{\int_{\Omega} #1 \ \dx}

\newcommand{\vv}{\vc{v}}

\newcommand{\ep}{\varepsilon}

\font\F=msbm10 scaled 1000
\newcommand{\R}{\mbox{\F R}}

\definecolor{Cgrey}{rgb}{0.85,0.85,0.85}
\definecolor{Cblue}{rgb}{0.50,0.85,0.85}
\definecolor{Cred}{rgb}{1,0,0}
\definecolor{fancy}{rgb}{0.10,0.85,0.10}

\newcommand\Cbox[2]{%
    \newbox\contentbox%
    \newbox\bkgdbox%
    \setbox\contentbox\hbox to \hsize{%
        \vtop{
            \kern\columnsep
            \hbox to \hsize{%
                \kern\columnsep%
                \advance\hsize by -2\columnsep%
                \setlength{\textwidth}{\hsize}%
                \vbox{
                    \parskip=\baselineskip
                    \parindent=0bp
                    #2
                }%
                \kern\columnsep%
            }%
            \kern\columnsep%
        }%
    }%
    \setbox\bkgdbox\vbox{
        \color{#1}
        \hrule width  \wd\contentbox %
               height \ht\contentbox %
               depth  \dp\contentbox
        \color{black}
    }%
    \wd\bkgdbox=0bp%
    \vbox{\hbox to \hsize{\box\bkgdbox\box\contentbox}}%
    \vskip\baselineskip%
}

\date{}

\begin{document}

%%%%%%%%%%%%%%%%%%%%%%%%%%%%%%%%

\title{Construction of weak solutions to compressible Navier--Stokes equations with general inflow/outflow boundary conditions via a numerical approximation}

\author{Young-Sam Kwon
\thanks{ {The work of the first author was partially supported by NRF-2019H1D3A2A01101128 and  NRF2020R1F1A1A01049805.}}   \and
Antonin Novotn{\' y}\thanks{The work of the second author was supported by Brain Pool program funded by the Ministry of Science and ICT through the National Research Foundation of Korea, NRF-2019H1D3A2A01101128.}}

\maketitle

\bigskip

\centerline{Department of Mathematics, Dong-A University} 
\centerline{Busan 49315, Republic of Korea, ykwon@dau.ac.kr}

\medskip
\centerline{University of Toulon, IMATH, EA 2134,  BP 20139} 
\centerline{839 57 La Garde, France, novotny@univ-tln.fr}

\maketitle

\bigskip

\begin{abstract}
The construction of weak solutions to compressible Navier-Stokes equations via a numerical method (including a rigorous proof of the convergence) is in a short supply, and so far, available only for one sole numerical scheme suggested in Karper [{\em Numer. Math.}, 125(3) : 441--510, 2013] for the
no slip boundary conditions and the isentropic pressure with adiabatic coefficient $\gamma>3$. Here we consider the same problem for the general non zero inflow-outflow boundary conditions, which is definitely more appropriate setting from the point of view of applications, but which is essentially more involved as far as the existence of weak solutions is concerned. There is a few recent proofs of existence of weak solutions in this setting, but none of them  is performed via a numerical method. The goal of this paper is to fill this gap.

The existence of weak solutions on the continuous level requires several tools of functional and harmonic analysis and differential geometry whose numerical counterparts are not known. Our main strategy therefore consists in rewriting of the numerical scheme in its variational form modulo remainders and to apply and/or to adapt to  the new variational formulation the tools developed in the theoretical analysis. In addition to the result, which is new, the synergy between numerical and theoretical analysis is the main originality of the present paper.  

%We modify the numerical scheme of Karper in order to accommodate general boundary conditions and  prove the convergence of its numerical solutions to weak solutions %for the same range of adiabatic coefficients as Karper did. Since the proofs on continuous level require several mathematical tools which do not have a counterpart %in the  numerical analysis (as compensated compactness and commutator lemmas), the leading idea is to rewrite the numerical scheme in terms of a variational %formulation modulo appearance of  remainders in such a way, that the methodology of the "continuous proofs" introduced in Lions[{\it Mathematical Topics in Fluid %Dynamics, Vol. 2, OUP}, 1998]  could be applied. This passes, among others,
%through the derivation of the effective viscous flux identity and an extensive use of the technique of renormalized solutions to the continuity equation. 

%Especially, the latter tool has to be, however, essentially modified to accommodate the non-zero inflow-outflow boundary conditions and the discrete solution. 
%Also, in this respect, the theory presented in the present paper is new, and of independent and broader interest since applicable to the numerical treatment of %transport equations with inflow-outflow, in general, with transporting coefficients of low regularity.
\end{abstract}


\section{Introduction}
Evolution of the density $\mathfrak{r}=\mathfrak{r}(t,x)$ and velocity $\mathfrak{u}=\mathfrak{u}(t,x)$ through the time interval $I=[0,T)$, $T>0$, $t\in [0,T)$
in a bounded (Lipschitz) domain $\Omega$, $x\in \Omega$ of a viscous compressible fluid can be described by the
Navier-Stokes equations. Unlike most of the theoretical  literature, in  this paper, we consider the Navier-Stokes system with the nonzero inflow-outflow 
boundary conditions, which is, from the point of view of applications, a more appropriate setting, than the "standard" setting with the no-slip or
Navier boudary conditions. The equations read: 
\begin{equation} \label{NS1}
\partial_t {\mathfrak{r}} + \Div ({\mathfrak{r}} {\mathfrak{u}}) = 0,
\end{equation}
$$
\partial_t ({\mathfrak{r}} {\mathfrak{u}}) + \Div ({\mathfrak{r}} {\mathfrak{u}} \otimes {\mathfrak{u}}) + \Grad p({\mathfrak{r}}) = 
\Div \mathbb{S}(\Grad\mathfrak{u}), 
$$
$$
\mathbb{S}(\Grad{\mathfrak{u}})=\mu\Grad{\mathfrak{u}} +(\mu+\lambda){\rm div}{\mathfrak{u}}\mathbb{I},\ \mu>0,\, \lambda+\frac 23\mu>0 \footnote{
One usually writes $\mathbb{S}(\Grad\vu)$ in its frame indifferent form, $\mathbb{S}(\Grad\vu)=\mu(\Grad\mathfrak{u}+(\Grad\mathfrak{u})^T) +\lambda{\rm div}{\mathfrak{u}}\mathbb{I}$. Both writing are equivalent in the strong formulation but not in the weak formulation later and in the numerical scheme.}
$$
with initial and boundary conditions,
\bFormula{ibc}
\mathfrak{r}(0)=\mathfrak{r}_0,\ \mathfrak{ u}(0)={\mathfrak{ u}}_0,\ \mathfrak{ u}|_{\partial\Omega}=\mathfrak{ u}_B,\
\mathfrak{ r}|_{\Gamma_{\rm in}}=\mathfrak{r}_B,
\eF
where (for simplicity)\footnote{The regularity of initial conditions could be relaxed to
$$
0<\mathfrak{r}_0\in L^\gamma(\Omega),\;{\bf{\mathfrak{ u}}_0}\in L^1(\Omega),\;\mathfrak{r}_0\mathfrak{u}_0^2\in L^1(\Omega).
$$
The last condition in (\ref{ruB}) means that $\mathfrak{u}_B\cdot\vc n=f|_{\partial\Omega}$ with $f\in C(R^3)$ in all points
of $\partial\Omega$ where $\vc n$  is defined.
}
\bFormula{ru0}
0<\mathfrak{r}_0\in C(\overline\Omega),\;{\bf{\mathfrak{ u}}_0}\in C(\overline\Omega),
\eF
{
\bFormula{ruB}
0<\mathfrak{r}_B\in C^1(\overline\Omega),\;{\bf{\mathfrak{ u}}_B}\in C^2(R^3),\; \;\mathfrak{u}_B\cdot\vc n\in C(\partial\Omega)
\eF
($\vc n$ is the outer normal to $\partial\Omega$)
}
are given initial and boundary data, and
{
\bFormula{inB}
\Gamma^{\rm in}={\{x\in \partial\Omega|\;\mathfrak{u}_B\cdot\vc n< 0\}}
\eF
is the inflow boundary.  For the further use, we denote
\bFormula{outB}
\Gamma^{\rm out}= {\{x\in \partial\Omega|\;\mathfrak{u}_B\cdot\vc n> 0\}},\;\Gamma^0=\partial\Omega\setminus(\Gamma^{\rm in}\cup\Gamma^{\rm out})
\eF
the outflow and slip boundary, respectively.
}

We suppose
\bFormula{pres1}
p\in C^1[0,\infty),\; p(0)=0,\; p'(\vr)>0\,\mbox{for all $\vr>0$}
\eF
(and we extend $p$ by zero to the negative real line so that $p\in C(R)$ if needed).  Further
\bFormula{pres2}
\underline \pi+\underline p\vr^{\gamma-1} \le  p'(\vr)\le\overline p\vr^{\gamma-1}+\overline \pi,\ \mbox{for all $\vr >0$ with some $\gamma>1$},
\eF
where $0\le\underline\pi<\overline\pi$, $0<\underline p\le\overline p$.

We associate to $p$ its Helmholtz function $H$,
\bFormula{H}
H(\vr)=\vr\int_1^\vr\frac{p(z)}{z^2}{\rm d}z\;\mbox{in particular}\;{ \vr H'(\vr)-H(\vr)=p(\vr).}
\eF

It is to be noticed that an iconic example of the isentropic pressure $p(\vr)=a\vr^\gamma$, $a>0$ complies with hypotheses (\ref{pres1}--\ref{pres2}).

The main goal of the paper is to construct {\em weak solutions} to problem (\ref{NS1}--\ref{pres2}) {\em via a numerical scheme.}

We remark that for the no-slip (or Navier or periodic) boundary conditions, existence of weak solutions to the Navier-Stokes system (\ref{NS1}--\ref{pres2}) is known  in the case $\gamma>3/2$, and it is nowadays a standard result, see \cite{FNP}
(and monographs by Lions \cite{Li}, Feireisl \cite{Fe},  Novotny, Straskraba \cite{NoSt}) for the no-slip boundary conditions. The same result for the
nonzero inflow-outflow boundary conditions is more recent, see \cite{ChJiNo}, \cite{ChoNoY}, \cite{KwNo} preceded by Girinon \cite{Gi}. In all these works, the weak solutions are constructed via  a several level approximations of the original system by more regular systems of PDEs and none of these approximations is a numerical one.

In the case of no-slip boundary conditions, Karper \cite{Ka} (and later on Feireisl et al.
in a different  context \cite{FeKaNo}, \cite{FeKaMi}, see also the comprehensive monograph \cite{FeKaPo}) constructed weak solutions to
the problem (\ref{NS1}--\ref{ru0}) via a specific finite volume/finite element
method -- which has been originally proposed  by Karlsen, Karper \cite{KaKa0},\cite{KaKa1}, \cite{KaKa2} and which we will call Karper's scheme -- at least for large values of the adiabatic coefficients, namely $\gamma>3$. So far, no other scheme was proved to have this convergence property, meaning that  the proof is strongly dependent on the structure of the scheme. 

The task in this paper is to construct the weak solutions to the compressible Navier-Stokes system with general boundary data via an adaptation of the Karper's scheme that we have suggested in \cite{KwNoNUM} in order to accommodate the non homogenous boundary data, for the same range of adiabatic coefficients as Karper did. In view of the newness of the theoretical results in this flow setting, this is still a challenging open question. To reach the same convergence property for the adiabatic coefficients $\gamma\le 3$
is even more challenging problem, which however, seems, so far, to be out of reach of our existing knowledge, even in the case of no-slip boundary data.

We wish to mention also the  works by Gallouet et al. \cite{GaHeLa}, Eymard et al. \cite{EyGaHeLa}, Perrin et al. \cite{PeSa} who construct via  a numerical scheme the stationary weak solutions with no-slip boundary data. In this paper, we deal exclusively with the non steady solutions for which the mathematical features of the  problem are essentially different even in the "simple" no-slip case. Moreover, it is known, that the steady problem with general inflow-outflow data is ill-posed, unless the pressure law is the so called hard sphere pressure law, cf. Feireisl et al. \cite{FeNoIHP}.

The unconditional  results on convergence of numerical
approximations to system (\ref{NS1}--\ref{inB}) --or to a similar systems without dissipation-- are in a short supply in the mathematical literature, and except the very recent result \cite{KwNoNUM}, they all deal with the no-slip or periodic
boundary conditions. In this respect, it is adequate to mention the monographs Feistauer et al. \cite{Feist}, \cite{FeistFeSt}, Kr\"oner \cite{Kr},
Eymard et al. \cite {EyGaHe}, works of
Feistauer et al. \cite{FeistFeLu}, \cite{FeistFeLuWa}, Gallouet et al. \cite{FeGa}, \cite{GaGaLaHe}, Jovanovic \cite{Jo}, Tadmor et al.
\cite{Ta0}, \cite{Ta1}, \cite{TaZh} and the works on error estimates \cite{GaHeMaNo} followed by \cite{GHMN-MAC}, \cite{MaNo}, \cite{FeHoMaNo}.
The convergence to Young measure valued solutions (in the spirit of the work  by Feireisl et al. \cite{FeGwGwWi})--and promoted in numerical analysis by  Fjordholm et al. \cite{FjMiTa0}, \cite{FjMiTa1}, \cite{FjKaMiTa}-- has been proved by Feireisl et al.  \cite{FeLu} followed by similar investigations in \cite{FeLuMi}, \cite{FeLuMiSh}, \cite{HoSh} and pushed further in Feireisl et al. in \cite{FeLuMiz}.  The convergence to similar type of solutions -- with Reynold defect (introduced in Abbatielo et al. \cite{AbFeNo}) has  been  shown in \cite{KwNoNUM}. In contrast with \cite{FeLu}, the latter work provides also a quantitative evaluation of the
convergence rate to the strong solutions. 
{None of the  works mentioned in this paragraph, however, deals with the convergence to weak solutions}.

Our approach is based on the following steps:

\begin{enumerate}
\item We take the finite element/finite volume numerical scheme suggested in \cite[Section 3.1]{KwNoNUM} (which is 
a modified Karper's scheme \cite{Ka} able to accommodate the general boundary conditions) and report from
\cite[Sections 5,6]{KwNoNUM} its fundamental properties: Existence of solutions with positive density, mass and energy balance, and the uniform estimates which can be deduced from these identities. This is the point, where the serious work of this  paper starts.
\item We rewrite the {\em numerical continuity equation} and the {\em numerical momentum equation} in their {\em variational forms} letting appear "remainders" vanishing as space $h$ and time $\Delta t$ discretizations tend to zero.
In accordance with our philosophy of "uniform approach" this still follows closely \cite[section 7]{KwNoNUM}. However, 
in view of further estimates and convergence properties needed for the convergence to weak solutions, the remainders
must be of better "quality" than in \cite{KwNoNUM}. This is the main (technical) point, where the restriction $\gamma>3$ is needed.
\item Once we derive the variational formulation of equations -- which resembles very much to the weak formulation on the continuous level -- our goal is to mimic the proof from the continuous case,
following \cite{ChJiNo}. This is so far the only exploitable strategy. In fact the proof needs on many places
special tools of harmonic analysis which exploit the particular non-linear structure of equations (as. e.g. the compensated compactness or commutator lemmas) and some
pieces of differential geometry (in particular, to treat the inflow field outside the inflow boundary), which are not currently available on the discrete level. 
This goes through several steps (designed on the continuous level for the first time by Lions \cite{Li} as far as compactness is concerned, and in \cite{ChJiNo} as far as the treatment of the inflow boundary is concerned) - and the principal output of this analysis is the strong convergence of the density sequence. These steps are:
\begin{enumerate}
\item Derivation of improved estimates of density sequence (by using the Bogovskii operator), which are good enough to
eliminate the possible concentrations in the density (or more exactly) in the pressure sequence.
\item Derivation of adequate numerical integration by parts formulas for treatment of the viscous terms in the momentum equation.
\item Derivation of a specific commutator lemma (a consequence of the celebrated Div-Curl lemma) for treatment of the material derivative of the linear momentum.
\item Using the three previous points leads to the derivation of the so called effective visous flux identity.
\item Adaptation of the DiPerna-Lions transport theory of renormalized solutions  \cite{DL} to the continuity equation
with non-zero inflow outflow boundary conditions in the spirit of \cite[Lemma 3.1]{ChJiNo}, on one hand, and to the
discrete solutions of the numerical scheme, on the other hand.  This part (presented in Sections \ref{741}--\ref{743}) is therefore of independent 
and broader interest since it is applicable, to numerical treatment of transport equations with transporting coefficients of low regularity, in general.
\item Once the concentrations in the density sequence eliminated, in order to prove the strong convergence, we must
still eliminate the oscillations. This is done by evaluating the defect of the sequence $\vr\log\vr$ exploiting the last
two properties and the effective viscous flux identity in the spirit of the approach of Lions \cite{Li}, which must, however, be modified in a non-trivial way, to take into account the non homogenous boundary conditions.
\end{enumerate}
\end{enumerate}

All convergence results are {\em "unconditional"}: to obtain them {\em we use only a priory estimates derived from the numerical scheme}. No a posteriori bounds are needed, in contrast with the most mathematical literature about the subject.
They are obtained through {\em an extensive  application of the functional analysis and theory of PDEs in the numerical
analysis which goes far beyond the standard approaches current in this domain of research.}

We consider,  the problem on a polygonal domain $\Omega$ covered by meshes that "fit" to the domain  and to its inflow and outflow boundaries,
cf. Sections \ref{mesh} and \ref{Upw}. In particular, we allow a finite number of inflow and outflow portions of the boundary and each of them must be "flat". 
In spite of this fact, it can be surprisingly "singular and realistic", as e.g. composed of an union of a finite number flat grids. We do not treat the problem of discretization error arising from "unfitted meshes" on general domains and on general "unflat" inflow-outflow boundaries. In the case of zero boundary velocity a solution has been proposed e.g. in the monograph \cite{FeKaPo}. In the non zero inflow-outflow situation, this is, in this context, an excellent open problem, where
even reasonable conditions of how the mesh has to "fit" to the physical domain and its inflow-outflow boundaries in order to conserve the property of convergence to weak solutions, remains to be determined.

Some more complex hydrodynamical models of compressible fluids with  similar
structure of convective terms which include transport and continuity equations--as e.g. models of compressible fluids with non-linear stress, see Abbatiello et al. \cite{AbFeNo}, fluids of compressible polymers, see Barrett and S\"uli \cite{BaSu}, models of compressible magneto-hydrodynamics, see \cite{HuWang}, or multi-fluid models with differential closure, see \cite{NoSCM}--can be treated on the basis of the methodology introduced in this paper. These studies would  certainly be of
a non negligible interest.

The paper is organized as follows. After introducing the bounded energy weak solutions in Section \ref{DD1}, we explain, in the next two sections the numerical setting to approximate them. The main result about the convergence of numerical solutions to weak solutions is stated in Theorem \ref{TN2}. The remaining parts of the paper are devoted to its proof. We derive the mass and energy balance and deduce the uniform estimates in Section \ref{SEst}. In Section \ref{NC}
we rewrite the algebraic numerical setting of the continuity and momentum  equations in their variational formulations. This will allow us to treat the
problem in a similar way as it is usually treated on the "continuous level": This treatment includes the improvement of density estimates in Section \ref{SEBOG}
and the investigation of the convergence in Section \ref{SConv} whose goal is to obtain the strong convergence of the density sequence. This follows the Lions methodology \cite{Li} passing  through the derivation of the effective viscous flux identity (proved in Section \ref{73}) and through the applications of the theory of renormalized weak solutions in Section \ref{74}. The latter is inspired by DiPerna-Lions \cite{DL} but has to be essentially modified in order to accommodate the
general boundary conditions (see Lemma \ref{LP2}, which is of independent interest). The proof requires many specific mathematical and numerical tools. They are gathered in the Appendix for reader's convenience.

We finish this introductory part by a remark concerning the notation: The special functional spaces are always defined in the text. For the classical Lebesgue, Sobolev, Bochner spaces and their duals, we use the standard notation, see e.g. Evans \cite{Ev}. Strong convergence in a Banach space is always denoted "$\to$", while "$\rightharpoonup$" means the weak convergence and "$\rightharpoonup_*$" means the star-weak convergence.

\subsection{Weak solutions}

\begin{Definition}{\rm [Bounded energy weak solutions]} \label{DD1}

The quantity $[\mathfrak{r}, \mathfrak{u}]$ is \emph{ bounded energy weak solution} of the problem (\ref{NS1}--\ref{inB}) in $(0,T) \times \Omega$
if the following is satisfied:
\begin{enumerate}
\item  
\[
0\le\mathfrak{r} \in C_{\rm weak}([0,T]; L^\gamma(\Omega)) \cap L^\gamma(0,T; L^\gamma(\partial \Omega; |\mathfrak{u}_B \cdot \vc{n}|
{\rm d}S_x)),\ \ \mbox{with some}\ \gamma > 1,
\]
\[
 \mathfrak{v}=\mathfrak{u} - \mathfrak{u}_B  \in L^2(0,T; W^{1,2}_0(\Omega; R^3)),
\mathfrak{m}:=\mathfrak{r}\mathfrak{u} \in C_{\rm weak}([0,T]; L^{\frac{2 \gamma}{\gamma + 1}}(\Omega; R^3));
\]
\bFormula{D1}
\mathfrak{r}{\mathfrak{u}}^2\in L^\infty(0,T;L^1(\Omega)),\;p(\mathfrak{r})\in L^1((0,T)\times\Omega));\footnote{We say that $r\in C_{\rm weak}([0,T];X)$, $X$ a Banach space, if
$r:[0,T]\mapsto X$ is defined on $[0,T]$, $r\in L^\infty(0,T;X)$ and
$<{\cal F},r>_{X^*,X}\in C[0,T]$ with any ${\cal F}$ in the dual space $X^*$ to $X$.}
\eF
\item The continuity equation
\begin{equation} \label{D2}
\intO{ {\mathfrak{r}} \varphi (\tau,\cdot) } +
\int_0^\tau \int_{\Gamma_{\rm out}} \varphi {\mathfrak{r}} {\mathfrak{u}}_B \cdot \vc{n} \ {\rm d}  S_x
+
\int_0^\tau \int_{\Gamma_{\rm in}} \varphi {\mathfrak{r}}_B {\mathfrak{u}}_B \cdot \vc{n} \ {\rm d} S_x
\end{equation}
$$
= \intO{ {\mathfrak{r}}_0 \varphi (0,\cdot) }+
\int_0^\tau \intO{ \Big[ {\mathfrak{r}} \partial_t \varphi + {\mathfrak{r}} {\mathfrak{u}} \cdot \Grad \varphi \Big] } \dt,\ {\mathfrak{r}}(0,\cdot) = {\mathfrak{r}}_{0},
$$
holds for any $0 \leq \tau \leq T$, and any test function $\varphi \in C^1([0,T] \times \Ov{\Omega})$;

\item The momentum equation
\bFormula{D3}
\intO{ {\mathfrak{r}} {\mathfrak{u}} \cdot \phi(\tau,\cdot) }  -\intO{ {\mathfrak{r}}_0 {\mathfrak{u}}_0 \cdot \phi(\tau,\cdot) }
\eF
$$
=\int_0^\tau \intO{ \Big[ {\mathfrak{r}} {\mathfrak{u}} \cdot \partial_t \phi + {\mathfrak{r}} {\mathfrak{u}} \otimes {\mathfrak{u}} : \Grad \phi
+ p({\mathfrak{r}}) \Div \phi - \mathbb{S}(\Grad{\mathfrak{u}}) : \Grad \phi \Big] } {\rm d}t
$$
holds for any $0 \leq \tau \leq T$ and any $\phi \in C_c^1([0,T] \times {\Omega}; R^3)$.
\item
The energy inequality
\begin{equation} \label{D4}
\intO{\left[ \frac{1}{2} {\mathfrak{r}} |{\mathfrak{v}}|^2 + H({\mathfrak{r}}) \right](\tau) }
+\int_0^\tau \intO{ \mathbb{S}(\Grad {\mathfrak{u}}):\Grad\mathfrak{u} } \dt
\end{equation}
$$
+\int_0^\tau  \int_{\Gamma_{\rm out}} H({\mathfrak{r}})  {\mathfrak{u}}_B \cdot \vc{n} \ {\rm d} S_x \dt 
+
\int_0^\tau\int_{\Gamma_{\rm in}} H({\mathfrak{r}}_B)  {\mathfrak{u}}_B \cdot \vc{n} \ {\rm d}S_x \dt	
$$
$$
\leq 
\intO{\left[ \frac{1}{2} {\mathfrak{r}}_0 |{\mathfrak{v}}_0|^2 + H({\mathfrak{r}}_0) \right] } 
-\int_0^\tau \intO{ \left[ {\mathfrak{r}} {\mathfrak{u}} \otimes {\mathfrak{u}} + p({\mathfrak{r}}) \mathbb{I} \right]  :  \Grad {\mathfrak{u}}_B } \dt
$$
$$
+ \int_0^\tau   \intO{ {{\mathfrak{r}}} {\mathfrak{u}} \cdot \Grad {\mathfrak{u}}_B   \cdot {\mathfrak{u}}_B  }
\dt 
+ \int_0^\tau  \intO{ \mathbb{S}(\Grad\mathfrak{u}) : \Grad {\mathfrak{u}}_B } \dt 
$$
holds 
for a.a. $\tau\in (0,T)$ with 
${\mathfrak{v}}_0=\mathfrak{u}_0 - {\mathfrak{u}}_B$.
\end{enumerate}
\end{Definition}

\section{Numerical setting}
\label{Nn}

We consider the same numerical setting as in \cite{KwNoNUM}.

\subsection{Mesh}
\label{mesh}

We suppose that the physical space is a \emph{polyhedral bounded domain} $\Omega \subset R^3$ that admits a \emph{tetrahedral} mesh ${\cal T}={\cal T}_h$;
the individual elements in the mesh will be denoted by $K=K_h \in {\cal T}$ (closed sets) and their gravity centers by $x_K$. Faces in the mesh are denoted as $\sigma=\sigma_h$ (close sets in $\R^2$) and their gravity centers by $x_\sigma$, whereas ${\cal E}={\cal E}_h$ is the set of all faces. \footnote{In the sequel
we shall omit in the notation the dependence on the "size" $h$ whenever there is no danger of confusion.} We  denote by ${\cal E}(K)$ the set of all faces
of $K\in {\cal T}$. Moreover, the set of faces in ${\cal E}$ belonging to the boundary $\partial \Omega$ is denoted ${\cal E}_{{\rm ext}}$, while ${\cal E}_{{\rm int}} = {\cal E} \setminus {\cal E}_{{\rm ext}}$. 

We denote by $h_K$ the diameter of $K$ and by $\mathfrak{h}_K$ the radius of the largest ball included in $K$. We call $h= \sup_{K\in{\cal T}}h_K$ the size of the mesh  and denote
$\mathfrak{h}= \inf_{K\in{\cal T}}{\mathfrak{h}}_K$. We also denote by ${\cal E}(K)$ the set of all faces
of $K\in {\cal T}$.

For two numerical quantities $a$, $b$, we shall write
\[
a \aleq b \ \mbox{if} \ a \leq c b, \ c > 0 \ \mbox{a constant}, \ a \approx b \ \mbox{if} \ a \aleq b \ \mbox{and} \ b \aleq a.
\]
Here, ``constant'' typically means a generic quantity independent of the size $h$ of the mesh and the time step $\Delta t$ used in the numerical scheme as well as other parameters as the case may be.

In addition, we require the mesh to be admissible in the sense of Eymard et al. \cite[Definition 2.1]{EyGaHe}:
\begin{enumerate}
\item For $K, L \in {\cal T}$, $K \ne L$, the intersection $K \cap L$ if non-empty is either a vertex, or an edge, or a face $\sigma \in {\cal E}$. In the latter case, we write
$\sigma=K|L$.
\item There holds
$$
h\approx \mathfrak{h}
$$
\end{enumerate}

We denote by $\vc{n}_{\sigma, K}$ the unit normal to the face $\sigma\in {\cal E}(K)$
outwards to $K$. On the other hand we associate to each element $\sigma \in {\cal E}$ 
a fixed normal vector $\vc{n}={\vc n}_\sigma$. If $\sigma\in{\cal E}_{\rm ext}$ then
$\vc n_\sigma$ is always the outer normal to $\partial\Omega$.

\subsection{Piecewise constant finite elements}\label{Q}

We introduce the space
\bFormula{Q0}
Q (\Omega)=Q_h(\Omega) = \left\{ g \in L^1(\Omega) \ \Big| \ 
g|_K = a_K \in \tn R \right\}
\eF
of piecewise constant functions (and we tacitly extend $g$ by zero outside $\Omega$, if convenient) along with the associated projection
\[
\Pi^Q=\Pi_h^Q: L^1(\Omega) \to Q(\Omega), \ \widehat v|_K:=\Pi^Q [v]|_K 
= v_K:= \frac{1}{|K|} \int_K{v}{\rm d}x.
\]

For a function $g\in Q(\Omega)$ and any $\sigma\in {\cal E}_{\rm int}$, we denote
\bFormula{gpm}
g^+_\sigma:=g^+_{\vc n_\sigma} = \lim_{\delta \to 0+ } g(x_\sigma + \delta \vc{n} ),\
g^-_\sigma:=g_{\vc n_\sigma} = \lim_{\delta \to 0+ } g(x_\sigma - \delta \vc{n} ).
\eF
Further, we define  the jumps and mean values over $\sigma$ 
(relative to $\vc n_\sigma$),
\begin{equation}\label{jump}
\ju{g}_\sigma=\ju{g}_{\sigma,\vc{n}}:= g^+_\sigma - g^-_\sigma,
\av{g}_\sigma:= \frac{1}{2} \left( g^+ + g^- \right).
\end{equation}

\subsection{Crouzeix-Raviart finite elements}\label{V}

A differential operator $D$ acting on the $x-$variable will be discretized as
\[
D_h v|_K = D (v|_K) \ \mbox{for any}\ v \ \mbox{differentiable on each element}\ K\in {\cal T}.
\]

The \emph{Crouzeix-Raviart finite element spaces} (see Brezzi and Fortin \cite{BreFor}, among others) are defined as
\bFormula{n1}
V (\Omega)=V_h(\Omega): = \Big\{ {v} \in L^1 (\Omega) \ \Big| \ {v}|_K = \mbox{affine function},
\eF
$$
 \ \int_\sigma v|_K{\rm dS}_x= \int_\sigma v|_L{\rm dS}_x \ \mbox{for any}\ \sigma =K|L\in {\cal E}_{\rm int} \Big\},
$$
together with
\bFormula{n1+}
V_{0} (\Omega)=V_{h,0} = \left\{ {v} \in V(\Omega) \ \Big| \ \int_{\sigma} {v} \ {\rm dS}_x = 0\ \mbox{for any}\ \sigma \in {\cal E}_{\rm ext} \right\}.
\eF
We extend any function $v\in V_0(\Omega)$ by $0$ outside $\Omega$, if convenient.

Next,
we introduce the associated projection
\begin{equation}\label{PiV}
\Pi_V=\Pi_h^V : W^{1,1}(\Omega) \to V (\Omega), \ {\tilde v}= \Pi^V[v]:=\sum_{\sigma\in{\cal E}}
v_\sigma\phi_\sigma,
\end{equation}
where
\[
v_\sigma= 
\frac 1{|\sigma|}\int_\sigma{ v }{\rm d}S_x 
\]
and $\{\phi_\sigma\}_{\sigma\in {\cal E}}\subset V(\Omega)$ is a basis in $V(\Omega)$
given by
$$
\frac 1{|\sigma'|}\int_{\sigma'}\phi_\sigma=\delta_{\sigma,\sigma'},\;(\sigma,\sigma')\in {\cal E}^2.
$$

\subsection{Convective terms, upwinds}\label{Upw}

Suppose, that we are given
 $\vu_B\in V(\Omega,R^3)$. We define
\begin{equation}\label{in-out}
{\cal E}^{\rm in}=\{\sigma\in {\cal E}_{\rm ext}\,|\,\vu_{B,\sigma}\cdot\vc n<0\},\quad
{\cal E}^{\rm out}=\{\sigma\in {\cal E}_{\rm ext}\,|\,\vu_{B,\sigma}\cdot\vc n> 0\},\quad
{\cal E}^0=\{\sigma\in {\cal E}_{\rm ext}\,|\,\vu_{B,\sigma}\cdot\vc n= 0\}.
\end{equation}
We say that the {\em mesh fits to the inflow-outflow boundaries} if
\bFormula{fit}
\overline{\Gamma^{\rm in}}=\cup_{\sigma\in {\cal E}^{\rm in}}\sigma,\quad
\overline{\Gamma^{\rm out}}=\cup_{\sigma\in {\cal E}^{\rm out}}\sigma.
\eF
We see that
\bFormula{in-out+}
{\rm int}_2\Gamma^0={\rm int}_2\Big(\cup_{\sigma\in {\cal E}^0}\sigma\Big)\;\mbox{and}\;{\cal E}_{\rm ext}={\cal E}^{\rm in}\cup {\cal E}^{\rm out}\cup
{\cal E}^{0} 
\eF
Here and hereafter, for $A\subset \partial\Omega$, ${\rm int}_2A$ denotes the interior of $A$ with respect to the trace topology
of $R^3$ on $\partial\Omega$.

We  define for  any $\sigma\in{\cal E}_{\rm int}$ for any
$g\in Q(\Omega)$  its upwind value on any $\sigma\in{\cal E}_{\rm int}$,
\begin{equation}\label{gup}
g_\sigma^{\rm up}=\left\{
\begin{array}{c}
 g_K\;\mbox{if $\sigma= K|L\in {\cal E}_{\rm int}$ and 
 $\vc u_\sigma\cdot\vc n_{\sigma,K}\ge 0$},\\
 g_L\;\mbox{if $\sigma= K|L\in {\cal E}_{\rm int}$ and 
 $\vc u_\sigma\cdot\vc n_{\sigma,K}< 0$}.
 %\\
 %\vr_K\;\mbox{if $\sigma\in {\cal E}^{\rm out}\cap{\cal E}(K)$}
 %\\
 %\vr_{B,\sigma}\;\mbox{if $\sigma\in {\cal E}_{\rm in}\cap{\cal E}(K)$}
\end{array}
\right\}.
\end{equation}

Finally, we associate to any face $\sigma\in{\cal E}_{\rm int}$
%$\cup{\cal E}^{\rm out}$  
the \emph{upwind} operator 
${\rm Up}_{\sigma}[{ g}, \vu]:={\rm Up}_{\sigma,\vc n}[{ g}, \vu]$ defined as
\bFormula{Up}
{\rm Up}_{\sigma,\vc n}[{ g} ,\vu] = { g}^- [\vu_\sigma \cdot \vc{n}]^+ + { g}^+ [\vu_\sigma \cdot \vc{n}]^-,\ \mbox{where}\ [c]^+ = \max \{ c, 0 \}, \ [c]^-= \min \{ c, 0 \},
\eF
and to any face $\sigma\in {\cal E}(K)$ the specific flux $F_{\sigma,K}$ (outwards the element $K$) defined as
\bFormula{F}
F_{\sigma,K}[g,\vu]=g_\sigma^{\rm up}\vu_\sigma\cdot{\vc n}_{\sigma,K}.
\eF

\subsection{Time discretization}
For simplicity, we shall consider the constant time step $\Delta t>0$ where $T=N\Delta t$, $N\in \tn{N}$ and
we set 
\bFormula{Ik}
I_k=(\tau_{k-1}, \tau_k],\; \tau_k=k\Delta t,\; k\in \mathbb{Z}.
\eF
 Suppose that we have functions $v^k:\Omega\to R$,
$k=0,\ldots, N$. For convenience, we set $v^k(x)=v^0(x)$ if $k\le 0$ and $v^k(x)=v^{N}(x)$ if $k>N$, and
we introduce numbers
$$
D_tv^k(x)=\frac{v^k(x)-v^{k-1}(x)}{\Delta t}, \; k\in \mathbb{Z}.
$$
Finally, we define
\bFormula{pwt}
v(t,x)=\sum_{k\in\mathbb{Z}}1_{I_k}(t)v^k(x),\; \quad D_tv(t,x)=\sum_{k\in\mathbb{Z}}1_{I_k}(t)D_tv^k(x),
\eF
{
\bFormula{pwt+}
\utilde{v}(t,x)=\sum_{k\in\mathbb{Z}} 1_{I_k}(t)\Big(v^{k-1}(x)+(t-(k-1)\Delta t)D_tv^k(x)\Big),\;\mbox{so that
$\partial_t\utilde{v}(t,x)=D_tv(x)$. }
\eF 
}
In the sequel, we denote by $L_{\Delta t}(0,T; Q_h(\Omega)):=L(0,T;Q(\Omega))$ resp.
$L_{\Delta t}(0,T; V_h(\Omega)):=L(0,T;V(\Omega))$ the spaces of piecewise constant functions from $[0,T]$
to $Q(\Omega)$ and $V(\Omega)$, respectively (constant on each $I_k$ and extended by the value in $I_0$
to the negative real axes and by the value in $I_{N}$ to $[T,\infty)$).

\section{Numerical scheme, main result}
\label{NN}

We shall construct weak solutions of problem (\ref{NS1}-\ref{inB}) by using the numerical approximation
suggested in \cite[Section 3.1]{KwNoNUM} which generalizes to the non-homogenous boundary conditions the scheme
originally suggested by Karlsen, Karper \cite{KaKa0}, \cite{KaKa1}, \cite{KaKa2}, \cite{Ka}.

\subsection{Numerical scheme}
We  are given the approximations of the initial and boundary conditions
\bFormula{N1}
\vr^0=\vr^0_{h} = \Pi^Q_h [\mathfrak{r}_0], \ \vu^0=\vu^0_{h} = \Pi^V[{\bf \mathfrak{u}}_0], \ \vr_B=\vr_{B,h}=\Pi^Q[\mathfrak{r}_B],\
\vu_B=\vu_{B,h}=\Pi^V[{\bf \mathfrak{u}}_B].
\eF
We are searching for
$$
\vr_{h,\Delta t}(t,x)=\sum_{k=1}^{N}\sum_{K\in{\cal T}}1_{I_k}(t)1_K(x)\vr^k_{K,h,\Delta t},
\;
\vu_{h,\Delta t}(t,x)=\sum_{k=1}^{N}\sum_{K\in{\cal T}}1_{I_k}(t)1_K(x)\vu^k_{K,h,\Delta t}
$$
where
\bFormula{N1+}
\vr^k\in Q(\Omega),
\vr^k_{h,\Delta t}> 0,\;\vu^k_{h,\Delta t}\in V(\Omega;R^3),\;\vc v^k=\vu^k-\vu_B\in
V_0(\Omega;R^3),\; k=1,\ldots,N
\eF
such that the following algebraic equations
(for the unknowns $\vr^k_K$, $\vu^k_\sigma$, $k=1,\ldots, N$, $K\in {\cal T}$,
$\sigma\in {\cal E}$) are satisfied:\footnote{In what follows,
we omit the indexes ``$h$'' and/or ``$\Delta t$'' and write simply $\vr^k$ instead of $\vr^k_{h,\Delta t}$,
$\vr$ instead of $\vr_{h,\Delta t}$,
etc.,  in order to avoid the cumbersome notation, whenever there is no danger of confusion.}
\begin{enumerate}
\item{\textsc{Approximation of the continuity equation}}
\bFormula{N2}
\intO{ D_t \vr^k \phi } + \sum_{K\in {\cal T}}\sum_{\sigma \in {\cal E}(K)\cap{\cal E}^{\rm int}}
\int_\sigma F_{\sigma,K}(\vr^k,\vu^k)\phi{\rm d}S_x
\eF
$$
+
\sum_{\sigma\in{\cal E}^{\rm out}}\int_\sigma \vr^{k}\vu_{B,\sigma}\cdot\vc n_\sigma
\phi{\rm d}S_x
+\sum_{\sigma\in {\cal E}^{\rm in}}\int_{\sigma}
\vr_{B}\vu_{B,\sigma}\cdot\vc n_{\sigma}\phi{\rm d}S_x
$$
$$
{ +  h^\omega\sum_{\sigma\in{\cal E}_{\rm int}}\int_\sigma[[\vr^k]]_{\sigma}[[\phi]]_\sigma{\rm d}S_x}=0
$$
%$$
% {\tc +\kappa h^{\omega} \sum_{K\in {\cal T}}\sum_{\sigma \in {\cal E}(K)\cap{\cal E}^{\rm int}}
%\int_\sigma F_{\sigma,K}(\vr^k,\vc n_\sigma^{\rm up})\phi{\rm d}S_x}
%$$
for all $\phi \in Q (\Omega)$, where  $\omega>0$.
\item{\textsc{Approximation of the momentum equation}}
\bFormula{N3}
\intO{ D_t (\vr^k \widehat{\vv}^k) \cdot \widehat\phi }
+\sum_{K\in {\cal T}}\sum_{\sigma \in {\cal E}(K)\cap{\cal E}_{\rm int}}\int_\sigma F_{\sigma,K}(\vr^k\widehat{\vv}^k,\vu^k)\cdot\widehat\phi{\rm d}S_x
\eF
$$
{+ \sum_{K\in {\cal T}}\sum_{\sigma \in {\cal E}(K)}\int_\sigma\vr^k\widehat\vu^k\cdot{\vc n}_{\sigma,K}\vu_B\cdot\widehat\phi
{\rm d} S_x} 
%+\intO{\vr^k\widehat\vu^k\cdot\nabla_h(\vu_B\cdot\widehat\phi)}
$$
$$
+
\sum_{\sigma\in{\cal E}^{\rm out}}\int_\sigma \vr^{k}\vu_{B,\sigma}\cdot\vc n_\sigma\widehat\vv^k\cdot\widehat\phi{\rm d}S_x
+
\sum_{\sigma\in {\cal E}^{\rm in}}\int_\sigma \vr_B\vu_{B,\sigma}\cdot\vc n_\sigma\widehat\vv^k\cdot\widehat\phi{\rm d}S_x
$$
$$
+ \intO{ \Big(\mathbb{S}(\nabla_h\vu^k):\Gradh \phi 
-
p_h(\vr^k) \Divh \phi\Big) }
 +  h^\omega \sum_{\sigma\in{\cal E}_{\rm int}}\int_\sigma[[\vr^k]]_{\sigma}\{\widehat\vv^k\}_{\sigma}[[\widehat\phi]]_\sigma{\rm d}S_x
 =0.
$$
for any $\phi \in V_{0}(\Omega;R^3)$. In the above,
\bFormula{ph}
p_h(\vr)=p(\vr) +h^\beta \vr^2, 
\eF
where $\beta>0$ is a small number which will be determined later.
\end{enumerate}

\bRemark{tr}
\begin{enumerate}
\item It is to be noticed that the background linear momentum $\vr\vu_B$ in the momentum equation is not "upwinded". If it were "upwinded"
we would loose the derivation of uniform estimates from the energy balance. This observation seems to have an universal character valid for 
the discretizations of the inflow/outflow problems via the finite volume methods, in general.
\item In agreement with (\ref{NS1}), here and in the sequel, $\intO{\mathbb{S}(\nabla_h\vu^k):\nabla_h \phi}$
means exactly $\intO{(\mu\nabla_h\vu^k:\nabla_h\phi+(\mu+\lambda){\rm div}_h\vu^k{\rm div}_h\phi)}.$ This form is important for the estimates. 
Indeed, it is well known that the Korn
inequality does not hold in the  Crouzeix-Raviart finite element space.
\item It is proved in \cite[Theorem 3.1]{KwNoNUM} that under assumptions { (\ref{fit}), (\ref{pres1}--\ref{pres2})}, (\ref{ru0}--\ref{ruB}),
the algebraic system (\ref{N1}--\ref{N3}) admits at least one solution $(\vr_{h,\Delta t},\vu_{h,\Delta t})$. 
Any of its solutions has a strictly positive density. The uniqueness of the numerical solutions to system (\ref{N1}--\ref{N3}) is, however,
an open problem.
\item If in the condition (\ref{pres2}) $\underline\pi>0$ 
then the perturbation (\ref{ph}) is not needed for the proof of convergence (and we can take in the numerical scheme simply $p$ 
at place of $p_h$).
\end{enumerate}
\eR

\subsection{Main results}

The main result deals with the case $h \approx\Delta t$. It guarantees a {\it convergence of a  subsequence of numerical solutions to a weak solution.}  
%Before announcing the statement, we
%shall specify the properties of the inflow boundary: We shall say that the inflow boundary
%is {\em mesh compatible} if
%\bFormula{Nfit}
%\Gamma^{\rm in}=\cup_{k=1}^{\overline k} \Gamma^{\rm in}_k,
%\eF
%where:
%\begin{enumerate}
%\item For all $k=1,\ldots,\overline k$, $\Gamma^{\rm in}_k$ is a domain (open and connected set) in a two dimensional
%hyperplane and $\Gamma_k^{\rm in}\cap\Gamma^{\rm in}_j=\emptyset$.
%\item If $k,j=1,\ldots,\overline k$, $k\neq j$ and $\Gamma^{\rm}_k$, $\Gamma^{\rm}_j$
%are domains in the same hyperplane, then $\overline{\Gamma^{\rm}_k}\cap\overline{\Gamma^{\rm}_j}=\emptyset$.
%\item There is $\overline s>0$ such that $B^+(\Gamma^{\rm in}_k;\overline s)\cap
%B^+(\Gamma^{\rm in}_j;\overline s)=\emptyset$, where
%$$
%B^+(\Gamma^{\rm in};\overline s):=\{z=x_B+s\vc n(x_B),\;s\in (0,\overline s),\; x_B\in \Gamma_k^{\rm in}\}.
%$$
%\end{enumerate}

\bTheorem{N2}
Let $h=\Delta t$. Suppose that the pressure satisfies assumptions (\ref{pres1}--\ref{pres2}) with $\gamma>3$ and that the initial and boundary conditions verify
(\ref{ru0}--\ref{ruB}). Suppose further that the mesh fits to the inflow-outflow boundaries, cf. 
(\ref{fit}). 
%and that the inflow boundary is mesh compatible, see (\ref{Nfit}).

 Consider a sequence of numerical solutions $[\vr_h,\vu_h]\in L_{\Delta t}(0,T; Q_h(\Omega))\times L_{\Delta t}(0,T,V_h(\Omega))$ of the problem (\ref{N1}--\ref{N3}) with
$\beta\in (0,\min\{\frac 12,\frac{2\gamma-6}\gamma\})$ and $\omega\in (0,1-\beta)$. Then we have:
There exists a subsequence $[\vr_h,\vu_h=\vv_h+\vu_B]$ (not relabeled) such that
$$
\vr_h\to\mathfrak{r}\;\mbox{in $L^q(0,T;L^q(\Omega))$},\; 1\le q<\gamma+1,
$$
\bFormula{L2L6}
\vv_h\rightharpoonup\mathfrak{v}\;\mbox{in $L^2(0,T;L^6(\Omega))$},
\eF
$$
\nabla_h\vv_h\rightharpoonup\nabla\mathfrak{v}\;\mbox{in $L^2(0,T;L^2(\Omega))$},
$$
where the couple $[\mathfrak{r},\mathfrak{u}=\mathfrak{v}+\mathfrak{u}_B]$
is a weak solution of the problem (\ref{NS1}--\ref{inB}) in the sense of Definition \ref{DD1}.
\eT
\bRemark{tr1}
\begin{enumerate}
\item In agreement with Item 4 in Remark \ref{Rtr}, if in the condition (\ref{pres2}) $\underline\pi>0$ 
then we can take $p_h=p$ (no parameter $\beta$ is needed). Then Theorem \ref{TN2} holds with $\omega\in (0,1)$.
%\item  We denote
%$$
%\Gamma^{\rm bd}=(\overline{\Gamma^{\rm out}}\cap\Gamma^0)\cup
%(\overline{\Gamma^{\rm in}}\cap\Gamma^0)\cup (\overline{\Gamma^{\rm out}}\cap\overline{\Gamma^{\rm in}}).
%$$
%Since in our case $\Gamma^{\rm bd}$ is an union of finite number of (possibly intersecting) segments and isolated points, 
%we have for all $s>0$,
%\bFormula{small}
%U_s:=\Big\{x\in R^3\,\Big |\,{\rm dist}(x,\Gamma^{\rm bound})<s\Big\},\;
%|U_s|\aleq s^2
%\eF
%where, clearly, $U_s$ is an open set. This property plays an important role in the proof of strong convergence of the density sequence within the
%process of the construction of weak solutions.
\item  The authors of \cite{AbFeNo} have introduced the notion of weak solutions with Reynolds defect and, proved, under assumption
\bFormula{pres4}
H(0)=0,\; H-\underline a p,\quad \overline a p-H\;\mbox{are convex functions for some $0<\underline a<\overline a$}
\eF
the weak strong uniqueness principle for such solutions. 
(It is to be noticed that  the iconic isentropic pressure $p(\vr)=a\vr^\gamma$, $a>0$, $\gamma>1$ complies with (\ref{pres4}) as well.) 
\begin{enumerate}
\item
Any weak solution of problem  is a weak solution with the Reynolds defect $0$. Therefore, 
if the problem (\ref{NS1}--\ref{inB}) admits a strong solution $[r,\vc U=\vc V+\mathfrak{u}_B]$ in the class 
\bFormula{tests}
\vc V \in C^1([0,T] \times \Ov{\Omega}; R^3),\,\Grad^2\vc V\in  C([0,T] \times \Ov{\Omega}; R^3),\
\vc V|_{\partial \Omega} = 0,
\eF
$$
r \in C^1([0, T] \times \Ov{\Omega}), \ \underline r:=\inf_{(0,T) \times \Omega} r > 0.
$$
then the  limit in (\ref{L2L6})
$[\mathfrak{r},\mathfrak{u}]$ is equal to $[r,\vc U]$, see \cite[Theorem 6.3]{AbFeNo}.  In this case the whole sequence $[\vr_h,\vu_h]$ converges to $[r,\vc U]$. 
\item Moreover, there exists $\alpha=\alpha(\gamma,\omega,\beta)>0$  and a positive number $C$ dependent on 
\bFormula{C}
\underline r,\;\overline r:=\sup_{I\times\Omega}r,\; \|\Grad r,\partial_t r, \vc V,\Grad\vc V, \Grad^2\vc V,\partial_t \vc V, \mathfrak{u}_B,
\Grad\mathfrak{u}_B, \Grad^2\mathfrak{u}_B\|_{C([0,T]\times\overline\Omega)}
\eF
such that
$$
\Big[{\cal E}\Big(\vr_h,\widehat\vv_h| r,\vc V\Big)\Big]_0^\tau +\int_0^\tau\Big(\|\vu_h-\tilde{\vc U}\|_{L^2(\Omega)}^2+
\|\nabla_h(\vu_h-\tilde{\vc U)}\|_{L^2(\Omega)}^2\Big) \aleq c \,h^\alpha, 
$$
where
\bFormula{calE}
{\cal E}(\vr,\vu|r,\vc U)=\intO{\Big(\frac 12\vr|\vu-\vc U|^2{\rm d}x+E(\vr|r)\Big)},\;
E(\vr|r)=H(\vr)-H'(r)(\vr-r)-H(r),\footnote{Indeed, the functional ${\cal E}$ is always positive and vanishing if
$\vu=\vc U$, $\vr=r$,
offering thus a natural evaluation of the "distance" between vector fields $\vu$, $\vc U$ and positive scalar fields $\vr$ and $r$.} 
\eF
see \cite[Lemma 10.2]{KwNoNUM} and compare with Gallou{e}t at al. \cite{GaHeMaNo}.
\end{enumerate}
\item Local in time existence of strong solutions to problem (\ref{NS1}--\ref{inB}) notably with non zero inflow outflow
boundary conditions is discussed in Valli, Zajaczkowski \cite{VaZa}.
Existence of weak solutions to the problem (\ref{NS1}--\ref{ruB}) in its full generality, has been obtained in
\cite{ChJiNo}, \cite{ChoNoY} for $\gamma >3/2$. Theorem \ref{TN2} provides an alternative proof of existence of weak solutions via a numerical scheme for adiabatic coefficients $\gamma>3$.
Weak solutions with Reynolds defect have been constructed 
in \cite{AbFeNo} and \cite{KwNoNUM} for the whole range $1<\gamma<\infty$ of the adiabatic coefficients. 
\end{enumerate}
\eR

\section{Energy balance and uniform estimates}\label{SEst}
%%%%%%%%%%%%%%%%%%%%%
We denote by $H_h$ the Helmholtz function corresponding to the pressure $p_h$. It reads
\bFormula{Hh}
H_h(\vr)=H(\vr) +h^\beta\vr^2
\eF
We report the mass and energy balance proved in \cite[Lemma 6.1]{KwNoNUM}.

\bLemma{ebalance}{\rm [Mass conservation and energy balance]}
Suppose that the pressure satisfies assumptions (\ref{pres1}). Then any solution of $(\vr,\vu)$ of the
algebraic system (\ref{N1}--\ref{N3}) satisfies for all $m=1,\ldots,N$ the following:
\begin{enumerate}
 \item Mass conservation
 \bFormula{mbalance}
 \vr^m>0,\;
 \intO{\vr^m}+\Delta t\sum_{k=1}^m\sum_{\sigma\in{\cal E}^{\rm out}}\int_\sigma\vr^k\vu_{B,\sigma}\cdot\vc n_\sigma{\rm d}S_x=
 \intO{\vr^0}-\Delta t\sum_{k=1}^m\sum_{\sigma\in{\cal E}^{\rm in}}\int_\sigma\vr_B\vu_{B,\sigma}\cdot\vc n_\sigma{\rm d}S_x.
 \eF
\item Energy balance. There exist $\overline\vr^{k-1,k}\in Q(\Omega)$, $\overline\vr_K^{k-1,k}\in [\min\{\vr_K^{k-1},\vr_K^k\},
\max\{\vr_K^{k-1},\vr_K^k\}]$, $K\in{\cal T}$ and $\overline\vr^{k,\sigma}\in [\min\{\vr_\sigma^{-},\vr_\sigma^+\}, 
\max\{\vr_\sigma^{-},\vr_\sigma^+\}]$, $\sigma\in{\cal E}_{\rm int}$, $k=1,\ldots,N$, such that
\bFormula{ebalance}
\intO{\frac 12\vr^m |\widehat\vv^m|^2} + \intO{H_h(\vr^m) }
+\Delta t  \sum_{k=1}^m\intO{\mathbb{S}(\nabla_h\vv^k):\nabla_h\vv^k}
\eF
 $$
 +
 \frac 12\sum_{k=1}^m\intO{\Big( \vr^{k-1} |\widehat\vv^k-\widehat\vv^{k-1}|^2+H_h''(\overline\vr^{k-1,k})
  |\vr^k-\vr^{k-1}|^2\Big)
 }
 $$
 $$
 +\frac {\Delta t}2\sum_{k=1}^m\sum_{\sigma\in {\cal E}_{\rm int}}\int_\sigma |{\rm Up}_\sigma(\vr^k,\vu^k)|\,\ju{\widehat\vv^k}_\sigma^2{\rm d}S_x 
+ h^\omega { {\Delta t} \sum_{k=1}^m}\sum_{\sigma\in {\cal E}_{\rm int}}\int_\sigma[[\vr^k]]_\sigma[[H_h'(\vr^k)]]_\sigma{\rm d} S_x
$$
 $$
+\frac {\Delta t} 2\sum_{k=1}^m\sum_{\sigma \in {\cal E}_{\rm int}}
\int_\sigma H_h''(\overline\vr^{k,\sigma})\ju{\vr^k}_\sigma^2|\vu^k_\sigma\cdot\vc n_{\sigma}|{\rm d}S_x
 +
 {\Delta t}\sum_{k=1}^m\sum_{\sigma\in {\cal E}^{\rm out}}\int_\sigma \vr^k
 \vu^k_\sigma\cdot\vc n_\sigma|\widehat\vv^k|^2{\rm d}S_x
$$
$$
 + {\Delta t}\sum_{k=1}^m\sum_{\sigma \in {\cal E}^{\rm out}}\int_\sigma H_h(\vr^k)\vu_{B,\sigma}\cdot\vc n_\sigma
{\rm d}S_x  
+ \Delta t\sum_{k=1}^m\sum_{\sigma \in {\cal E}^{\rm in}} E_{H_h}(\vr_B|\vr^k)
|\vu_{B,\sigma}\cdot\vc n_\sigma|
{\rm d}S_x
$$
$$
 =\intO{\Big[\frac 12\vr^0 |\widehat\vv^0|^2 + H_h(\vr^0)\Big] } 
 +
 \Delta t\sum_{k=1}^m\sum_{\sigma \in {\cal E}^{\rm in}}\int_\sigma H_h(\vr_B)|\vu_{B,\sigma}\cdot\vc n_\sigma|
{\rm d}S_x
$$
$$
{ +\Delta t  \sum_{k=1}^m\intO{\mathbb{S}(\nabla_h\vu^k):\nabla_h\vu_B}}
 -\Delta t\sum_{k=1}^m
  \intO{ p_h(\vr^k){\rm div}_h\vu_{B}}
 $$
 $$
{
 -\Delta t\sum_{k=1}^m\intO{\vr^k\widehat\vu^k\cdot\nabla_h\vu_B\cdot\widehat\vv^k}
 - {\Delta t}\sum_{k=1}^m\sum_{\sigma\in {\cal E}^{\rm in}}\int_\sigma \vr_B\vu_{B,\sigma}\cdot\vc n_\sigma |\widehat\vv^{k}|^2{\rm d}S_x}.
$$
%$$
%{\tc +\kappa h^\omega\,\Delta t\sum_{k=1}^m \sum_{\sigma\in{\cal E}_{\rm int}} \int_{\sigma}[[\vr^k]]_\sigma\{\widehat\vu_B\}_\sigma[[\widehat\vv^k]]_\sigma{\rm d}S_x}
%$$
 \end{enumerate}
 \eL

As in \cite[Lemma 6.2]{KwNoNUM}, we can readily deduce from Lemma \ref{Lebalance} the uniform estimates. They are gathered in
the next lemma.
 
\bLemma{estimates}
We denote $I=[0,T)$,
$$
\overline E_0=\sup_{h\in (0,1)}E_{0,h},\;\; E_{0,h}= \intO{ \left[ \frac{1}{2} \vr^0_h | \widehat{\vu}^0_h |^2 +  H_h(\vr_0)\right] }.
$$
We suppose that the pressure satisfies assumption (\ref{pres1}--\ref{pres2}).

Then there exists a number $\mathfrak{D}$,
$
\mathfrak{D}:= \mathfrak{D}(\|\mathfrak{u}_B\|_{C^1(\overline\Omega)},\|\mathfrak{r}_B\|_{C(\overline\Omega)}, \overline E_{0}, T,\Omega)\aleq 1,
$
such that any solution $(\vr,\vu)=(\vr_{h,\Delta t}, \vu_{h,\Delta t})$ of the discrete problem (\ref{N1}--\ref{N3}) admits the following estimates:
\bFormula{e0}
\vr>0, \|\vr\|_{L^\infty(I;L^\gamma(\Omega))} \aleq 1 ,\; \|p(\vr), H(\vr)\|_{L^\infty(I;L^1(\Omega))}\aleq 1,
\eF
\bFormula{e1}
\|\vr|\mathfrak{u}_{B}\cdot\vc n|^{1/\gamma}\|_{L^\gamma(I;L^\gamma(\partial\Omega))}\aleq 1, \|H(\vr)|\mathfrak{u}_{B}\cdot\vc n|\|_{L^1(I;L^1(\partial\Omega))}\aleq 1,
\eF
\bFormula{e0+}
\|h^{\beta/2}\vr\|_{L^\infty(I;L^2(\Omega))} \aleq 1 ,
\eF
\bFormula{e1+}
\|h^{\beta/2}\vr|\mathfrak{u}_{B}\cdot\vc n|^{1/2}\|_{L^2(I;L^\gamma(\partial\Omega))}\aleq 1,
\eF
\bFormula{e2}
{\rm sup}_{\tau \in (0,T)} \| \sqrt{ \vr } \widehat{\vu} (\tau, \cdot) \|_{L^\infty(I;L^2(\Omega)} \aleq 1,
\eF
\bFormula{e3}
\|\nabla_h \vv\|_{L^2(I\times\Omega)} \aleq 1,\quad \|\vv\|_{L^2(I; L^6(\Omega))} \aleq 1,
\eF
\bFormula{e4}
\sum_{k \geq 0} \intO{   \left[ h^\beta
\left| \vr^k - \vr^{k-1}  \right|^2 + \vr^{k-1}\left| \widehat{\vv}^k - \widehat{\vv}^{k-1}  \right|^2 \right] } \aleq 1,
\eF
\bFormula{e5}
\sum_{\sigma\in {\cal E}_{\rm int}}\int_0^T\int_\sigma |{\rm Up}_\sigma(\vr,\vu)
| \, [[\vv]]^2_\sigma{\rm d}S_x{\rm d}t \aleq 1,
\eF
\bFormula{e6}
 h^\beta\sum_{\sigma\in{\cal E}_{\rm int}} \int_0^T \int_{\sigma} 
 \ju{ \vr  }_\sigma^2 | \vu_\sigma \cdot \vc{n}| {\rm dS}_x  \dt \aleq 1,\quad
h^{\omega+\beta}\sum_{\sigma\in{\cal E}_{\rm int}} \int_0^T \int_{\sigma} 
 \ju{ \vr  }_\sigma^2  {\rm dS}_x  \dt \aleq 1.\footnote{The second estimate in (\ref{e3}) follows from the fisrt one by virtue of the discrete Sobolev inequality (\ref{SobolevV0}). In the notation (\ref{?}) the latter formula says 
$$
h^{\frac{\omega+\beta+1} 2}\|\vr_h\|_{L^2(0,T;Q^{1,2}(\Omega))}\aleq 1.
$$
in terms of the Sobolev "broken" norm.}
\eF
%Finally,
%\bFormula{e7}
%h^\omega\sum_{\sigma\in{\cal E}_{\rm int}} \int_0^T \int_{\sigma} 
% \ju{ \vr  }_\sigma^\gamma {\rm dS}_x  \dt \aleq 1.
%\eF
\eL

\section{Consistency}
\label{NC}

Having collected all the available uniform bounds, our next task is to verify that our numerical method is \emph{consistent} with the variational formulation of the original problem. In this task, we follow the strategy of \cite{KwNoNUM} having in mind the following observation traced back to Karper \cite{Ka}: It is well known from the existence proofs of weak solutions, that further estimates and identities -- namely more integrability of density and the so called effective viscous flux identity -- obtained
by testing of the momentum equation by a special test functions involving solution of the continuity equation (which are not as regular as test functions in \cite{KwNoNUM}) are needed. Due to this fact, we have to admit in the variational formulations of the
continuity and momentum numerical methods less regular test functions which allow
these tests.

\subsection{Consistency of the continuity equation}\label{CE}

The goal of this section is to prove the following lemma.
\bLemma{c-consistency}{\rm [Consistency for the continuity equation]}
Let $h=\Delta t$ and let the pressure satisfy the hypotheses (\ref{pres1}--\ref{pres2}) with $\gamma>3$. 
Then we have:
\begin{enumerate}
\item
For any $\beta\in (0,\min\{\frac 12,\frac {2\gamma-6}\gamma\})$ and any $\omega>0$ there exists $\alpha_C>0$ 
 such that  any numerical solution 
of problem (\ref{N1}--\ref{N3})   satisfies the continuity equation in the variational form
\bFormula{ce-v+}
\intO{\utilde{\vr}\phi(\tau)}-\intO{{\vr}^0\phi(0)}-\int_{0}^\tau\intO{\Big(\utilde{\vr}\partial_t\phi+\vr{ \vu}\cdot\Grad\phi\Big)}{\rm d}t
\eF
$$
+  \int_{0}^{\tau}\int_{\Gamma^{\rm out}}
\vr\mathfrak{u}_{B}\cdot \vc n\phi{\rm d}S_x{\rm d}t
+  \int_{0}^{\tau}\int_{\Gamma^{\rm in}}
\vr_B\mathfrak{u}_{B}\cdot \vc n\phi{\rm d}S_x{\rm d}t
=\int_{0}^{\tau}< R^C_{h}(t),\phi>{\rm d}t 
$$
with any $\tau\in(0,T]$ and all $\phi\in C^1([0,T]\times\overline\Omega)$,
{ where $\utilde \vr$ is defined in (\ref{pwt+}),}
 where 
the remainder $R^C_h (t)$ satisfies 
\bFormula{NC1}
\left| \left< R^C_h (t), \phi \right> \right| \aleq h^{\alpha_C} r^C_h(t) \| \Grad \phi \|_{L^\infty(0,T;L^{\frac{6 \gamma}{\gamma + 6}}(\Omega))}, \
\| r^C_h \|_{L^2(0,T)} \aleq 1.
\eF
\item Alternatively, (\ref{ce-v+}) can be rewritten in the form
\bFormula{ce-v}
D_t\vr+ F= R_h^C\;\mbox{in $L^2(0,T; [W^{1,\frac{6\gamma}{\gamma+6}}(\Omega)]^*)$},
\eF
where
{
$$
F=F_h,\; <F_h,\phi>=-\int_{0}^T\intO{\vr\vu\cdot\Grad\phi}{\rm d}t
+  \int_{0}^{T}\int_{\Gamma^{\rm out}}
\vr\mathfrak{u}_{B}\cdot \vc n\phi{\rm d}S_x{\rm d}t
+  \int_{0}^{T}\int_{\Gamma^{\rm in}}
\vr_B\mathfrak{u}_{B}\cdot \vc n\phi{\rm d}S_x{\rm d}t,
$$
}
and
 $F_h$ is bounded in $L^2(0,T;[W^{1,\frac{6\gamma}{\gamma+6}}(\Omega)]^*)$,
$R_h^C\to 0$ in $L^2(0,T;[W^{1,\frac{6\gamma}{\gamma+6}}(\Omega)]^*)$.
\end{enumerate}
\eL
\bRemark{rc1}
Visiting the proof of Lemma \ref{Lc-consistency} we find that 
$$
\alpha_C\ge\min\Big\{\omega,\frac 14-\frac\beta 2, \frac {\gamma-3}\gamma-\frac\beta 2,
\frac{5\gamma-12}{6\gamma}\Big\}.
$$
\eR
\noindent
{\bf Proof of Lemma \ref{Lc-consistency}}\\
For $\phi \in C^1([0,T]\times\Ov{\Omega})$, we take $\widehat\phi(t)$ as a test function in the discrete continuity equation (\ref{N2}).  It is proved
in Section 7.1 in \cite[formula (7.2)]{KwNoNUM} (essentially using formula (\ref{Up3}) and (\ref{pwt+})) that equation (\ref{N2}) can be rewritten in the variational form
(\ref{ce-v+}),
where
$$
<{R}^C_h,\phi>=
 \sum_{K\in {\cal T}}\sum_{\sigma\in {\cal E}(K)\cap{\cal E}_{\rm int} }\int_\sigma ( \phi - \widehat\phi ) \ju{\vr}_{{\vc n}_{\sigma,K}}  \ [\vu_\sigma \cdot \vc{n}_{\sigma,K}]^- {\rm d}S_x  
$$
$$
+\sum_{K\in {\cal T}} 
\sum_{\sigma\in {\cal E}(K)}\int_{\sigma} (\widehat\phi-\phi) \vr(\vu - \vu_\sigma ) \cdot \vc{n}_{\sigma,K} \ {\rm dS}_x { +\intO{\vr(\phi-\widehat\phi){\rm div}_h\vu}}
$$
$$
{ - h^\omega\sum_{\sigma\in {\cal E}_{\rm int}}\int_\sigma [[\vr]]_\sigma [[\widehat\phi]]_\sigma {\rm d} S_x}
+  \sum_{\sigma\in {\cal E}^{\rm in}}\int_\sigma
(\vr-\vr_B)\vu_{B,\sigma}\cdot \vc n_{\sigma}(\widehat\phi-\phi){\rm d}S_x.
$$
%$$
%{\tc +  \sum_{\sigma\in {\cal E}^{\rm out}}\int_\sigma
%\vr(\mathfrak{u}_B-\mathfrak{u}_{B,\sigma})\cdot \vc n_{\sigma}\phi{\rm d}S_x
%+  \sum_{\sigma\in {\cal E}^{\rm in}}\int_\sigma
%\vr_B(\mathfrak{u}_B-\mathfrak{u}_{B,\sigma})\cdot \vc n_{\sigma}\phi{\rm d}S_x.
%}
%{\tb + \intO{\vr(\vu-\widehat\vu)\cdot\Grad\phi}}.
%$$

The goal now is to estimate conveniently all terms in the remainder $<R^C_h,\phi>$. To do this we shall use the tools evoked in Section \ref{Preln} and employ the bounds (\ref{e0}--\ref{e6}). 
\begin{enumerate}
\item 
 {\it The first term in $<{R}^C_{h},\phi>$} is bounded from above by
\bFormula{cest0}
\Big[ h^\beta\sum_{\sigma\in{\cal E}_{\rm int}}  \int_{\sigma} 
\ju{ \vr  }_\sigma^2 | \vu_\sigma \cdot \vc{n}_\sigma| {\rm dS}_x \Big]^{1/2}
\eF
$$
\times\Big[h^{-\beta}\sum_{K\in {\cal T}} 
\sum_{\sigma\in {\cal E}(K)\cap{\cal E}_{\rm int}}\int_\sigma  ( \widehat\phi - \phi )^2  |\vu_\sigma \cdot \vc{n}_{\sigma}| {\rm d}S_x\Big]^{1/2}, 
$$
where the first term in the product is controlled by (\ref{e6}). As far as for the second term in the  product, we have, if $\gamma>3$,
$$
h^{-\beta}\Big|\sum_{K\in {\cal T}} 
\sum_{\sigma\in {\cal E}(K)\cap{\cal E}_{\rm int}}\int_\sigma  ( \widehat\phi - \phi )^2  |\vu_\sigma \cdot \vc{n}_{\sigma}| {\rm d}S_x\Big|
$$
$$
\aleq h^{-\beta}\sum_{K\in {\cal T}} 
\sum_{\sigma\in {\cal E}(K)\cap{\cal E}_{\rm int}}\|\phi-\widehat\phi\|^2_{L^{\frac{6\gamma}{\gamma+6}}(\sigma)}\|\vu\|_{L^{\frac{3\gamma}{2\gamma-6}}(\sigma)}
$$
$$
 \aleq h^{1-\beta}\sum_{K\in {\cal T}} \|\vu\|_{L^{\frac{3\gamma}{2\gamma-6}}(K)} \|\nabla\phi\|^2_{L^{\frac{6\gamma}{\gamma+6}}(K)}
\aleq h^{1-\beta}\|\vu\|_{L^{\frac{3\gamma}{2\gamma-6}}(\Omega)} \|\nabla\phi\|^2_{L^{\frac{6\gamma}{\gamma+6}}(\Omega)}
$$
\begin{equation}\label{cest1}
\aleq h^{1-\beta+\min\{0,\frac{2\gamma-6)}{\gamma}-\frac 12\}}(\Delta t)^{-1/2} \|(\Delta t)^{1/2}\vu\|_{L^6(\Omega)} \|\nabla\phi\|^2_{L^{\frac{6\gamma}{\gamma+6}}(\Omega)}
\end{equation}
where $\|(\Delta t)^{1/2}\vu\|_{L^6(\Omega)}$ is bounded in $L^\infty((0,T))$ by virtue of (\ref{timeG}).

Indeed, to get the second line we have employed the H\"older and Jensen inequalities on $\sigma$, third line employs the trace estimates (\ref{trace}--
\ref{traceFE}) and one of Poincar\'e type inequalities listed in (\ref{PoincareV}) together with the H\"older's inequality. Finally, one concludes by the
"negative" interpolation estimates (\ref{interG}).

\item By the same token, the absolute value of {\it the second term in $<{ R}^C_{h}(t),\phi>$}  
$$
\sum_{K\in {\cal T}} 
\sum_{\sigma\in {\cal E}(K)}\int_{\sigma} (\widehat\phi-\phi) \vr(\vu - \vu_\sigma ) \cdot \vc{n}_{\sigma,K} \ {\rm dS}_x
$$
is, if $\gamma>3$, bounded from above by
\bFormula{cest2}
\aleq h\|\vr\|_{L^{\frac{6\gamma}{2\gamma-6}}(\Omega)}\|\nabla_h\vu\|_{L^2(\Omega)}\|\Grad\phi\|_{L^{\frac{6\gamma}{\gamma+6}}(\Omega)}
\eF
$$
\aleq h^{1+\min\{0,\frac{\gamma-3}\gamma-\frac 3\gamma\}}\|\vr\|_{L^{\gamma}(\Omega)}\|\nabla_h\vu\|_{L^2(\Omega)}
\|\Grad\phi\|_{L^{\frac{6\gamma}{\gamma+6}}(\Omega)}.
$$
{ {\it The third term 
%and the last term 
in $<{ R}^C_{h}(t),\phi>$} admits the same estimate.}
\item {\it The artificial density diffusion term in $<R^C_h(t),\phi>$} can be processed in the similar way,
\bFormula{cestad}
   h^\omega\Big|\sum_{\sigma\in {\cal E}_{\rm int}}\int_\sigma [[\vr]]_\sigma [[\widehat\phi]]_\sigma {\rm d} S_x\Big|
\aleq  h^{\omega} \sum_{K\in {\cal T}}\sum_{\sigma\in {\cal E}(K)\cap{\cal E}_{\rm int}}\|\vr\|_{L^2(\sigma)}\|[[\widehat\phi]]_\sigma\|_{L^2(\sigma)}
\aleq h^{\omega}\|\vr\|_{L^2(\Omega)}
\|\Grad\phi\|_{L^2(\Omega)}.
\eF
%for any value $\gamma>1$, where the first expression in the product is controlled by (\ref{e6}) in $L^2((0,T))$.
%
\item {\it The boundary term in  $<R^C_h(t),\phi>$} is controlled in the following way (without loss of generality, we perform the calculation for the first term with $\vr$ only -instead of $(\vr-\vr_B)$): If $\gamma\ge 12/5$, 
$$
 \Big|\sum_{\sigma\in {\cal E}^{\rm in}}\int_\sigma
\vr\vu_{B,\sigma}\cdot \vc n_{\sigma,K}(\widehat\phi-\phi){\rm d}S_x\Big|\aleq \sum_{K\in {\cal T},\, {\cal E}(K)\cap {\cal E}^{\rm in}\neq\emptyset}
\sum_{\sigma\in {\cal E}^{\rm in}}\|\vr\|_{L^\gamma(\sigma)}\|\widehat\phi-\phi\|_{L^{\gamma'}(\sigma)}
$$
$$
\aleq  \sum_{K\in {\cal T},\, {\cal E}(K)\cap {\cal E}^{\rm in}\neq\emptyset}\|\vr\|_{L^\gamma(K)}\|\nabla\phi\|_{L^{\gamma'}(K)}
\aleq \|\vr\|_{L^\gamma(\Omega)}\|\nabla\phi\|_{L^{\gamma'}({\cal U})}
$$
\bFormula{cest3}
\aleq 
\|\vr\|_{L^\gamma(\Omega)}\|\nabla\phi\|_{L^{\frac{6\gamma}{\gamma+6}}({\cal U})}|{\cal U}|^{\frac{5\gamma-12}{6\gamma}}
\aleq \|\nabla\phi\|_{L^{\frac{6\gamma}{\gamma+6}}({\Omega})} h^{\frac{5\gamma-12}{6\gamma}}
\eF
where ${\cal U}=\cup_{K\in {\cal T},\,{\cal E}(K)\cap {\cal E}^{\rm in}\neq\emptyset} K$.
Indeed, in the passage from the first to the second line, we have used trace estimates (\ref{trace}), (\ref{traceFE})
together with the first line in (\ref{PoincareV}), and for the rest the H\"older's inequalities as well as the fact that
the Lebesgue measure $|{\cal U}|$ of ${\cal U}$ is $\approx h$.
%
%\item {\it The last term in $<R^c_{h,\Delta t}(t),\phi>$} admits for any $\gamma>3$ the %bound
%$$
%|\intO{\vr(\vu-\widehat\vu)\cdot\Grad\phi}|\aleq h %\|\vr\|_{L^3(\Omega)}\|\nabla_h\vu\|_{L^2(\Omega)}\|\Grad\phi\|_{L^{6/5}(\Omega)}
%$$
%\bFormula{cest4}
%\aleq %h\|\vr\|_{L^\gamma(\Omega)}\|\nabla_h\vu\|_{L^2(\Omega)}\|\Grad\phi\|_{L^{\frac{6\gamma}{%\gamma+6}}(\Omega)}.
%\eF
%\end{enumerate}
 \item Resuming the results of calculations in (\ref{cest0}--\ref{cest3})  yields the first statement of Lemma \ref{Lc-consistency}. Formulas (\ref{ce-v}) and (\ref{NC1}) are thus proved.
\item Due to (\ref{pwt+}) and integration by parts, equation (\ref{ce-v}) can
be rewritten, in particular, in the form
$$
\int_0^T\intO{\Big(D_t\vr\phi-\vr{ \vu}\cdot\Grad\phi\Big)}{\rm d}t +  
\int_0^T\int_{\Gamma^{\rm out}}
\vr\mathfrak{u}_{B}\cdot \vc n \phi{\rm d}S_x{\rm d}t
$$
$$
+  \int_0^T\int_{\Gamma^{\rm in}}
\vr_B\mathfrak{u}_{B}\cdot \vc n\phi{\rm d}S_x{\rm d}t=\int_0^T<{ R}^C_h(t),\phi>{\rm d}t 
$$
with any $\phi\in C^1([0,T]\times\Omega)$.

By virtue of (\ref{e0}), (\ref{e3})
$$
|\intO{\vr{\vu}\cdot\Grad\phi}|\aleq \|\nabla\phi\|_{L^2(I; L^{\frac{6\gamma}{5\gamma-6}}(\Omega))}\aleq \|\nabla\phi\|_{L^2(I; L^{\frac{6\gamma}{\gamma+6}}(\Omega))},
$$
where we have used the continuous imbedding $L^{\frac{6\gamma}{\gamma+6}}(\Omega)\hookrightarrow
L^{\frac{6\gamma}{5\gamma-6}}(\Omega)$.
Seeing that, due to the trace theorem and the continuous inmbedding $W^{1,\frac{6\gamma}{\gamma+6}}(\Omega)\hookrightarrow
W^{1,\gamma'}(\Omega)$,
$$
\Big|\int_{\Gamma^{\rm out}}\vr\mathfrak{u}_B\cdot{\vc n}\phi{\rm d}S_x\Big|
\aleq \|\vr |\mathfrak{u}_B\cdot{\vc n}|^{1/\gamma}\|_{L^\gamma(\partial\Omega)}
\|\phi\|_{L^{\gamma'}(\partial\Omega)}\aleq \|\phi\|_{W^{1,\gamma'}(\Omega)}
\aleq \|\phi\|_{W^{1,\frac{6\gamma}{\gamma+6}}(\Omega)},
$$
we can rewrite the latter equation as identity (\ref{ce-v+}). This yields
the second statement of Lemma \ref{Lc-consistency}. The lemma is proved.
\end{enumerate}

\subsection{Momentum equation}

The goal of this section is to prove the following lemma.

\bLemma{m-consistency}{\rm [Consistency for the momentum equation]}
Let $h=\Delta t$ and let the pressure satisfy assumptions (\ref{pres1}--\ref{pres2}) with $\gamma>3$. Suppose that $\omega>0$.
Then we have: 
\begin{enumerate}
\item { For any $\beta\in (0,\min\{\frac 12,\frac{2\gamma-6}\gamma\})$ there exists $q_0>1$ such that for any  $q\in [1,q_0]$ there exists  $\alpha_M=\alpha_M(\beta,q,\omega)>0$,} 
 such that  any numerical solution 
of problem (\ref{N1}--\ref{N3})  
satisfies the momentum equation in the variational form, 
\bFormula{me-v+}
\intO{\utilde{\vr\widehat\vv}\cdot\phi(\tau,x)}-\intO{{\vr^0\widehat\vv^0}\cdot\phi(0,x)}
- \int_{0}^\tau\intO{\Big[\utilde{\vr\widehat\vv}\cdot\partial_t\phi+ \Big(\vr { {\vu }}\otimes \widehat{\vu}+p_h(\vr)\mathbb{I}\Big):\Grad\phi\Big] }{\rm d}t
\eF
$$
+\int_{0}^\tau\intO{\vr\widehat\vu\cdot\nabla_h(\vu_B\cdot\phi)} 
+
\int_{0}^\tau\intO{\mathbb{S}(\nabla_h\vu):\Grad\phi }{\rm d}t=
\int_{0}^\tau<{ R}^M_{h,\Delta t}(t),\phi>{\rm d}t
$$
with any $\tau\in (0, T]$ and all $\phi\in C_c^1([0,T]\times\Omega,R^3)$,
where 
the remainder $R^M_{h,h} (t)$ satisfies 
\bFormula{NM1}
\left| \left< R^M_{h,h} , \phi \right> \right| \aleq h^{\alpha_M} r^M_h(t) \| \Grad \phi \|_{L^{\gamma}(\Omega)}, \
{ \| r^M_h \|_{L^{q}(0,T)}} \aleq 1.
\eF 
\item Alternatively, (\ref{me-v+}) can be rewritten as follows:
\bFormula{me-v}
\intO{ D_t (\vr \widehat{\vv}) \cdot \phi } - \intO{ \vr { \vu } \otimes \widehat{\vu} : \Grad \phi } 
 +\intO{\vr\widehat\vu\cdot\nabla_h(\vu_B\cdot\phi)} 
\eF
$$
 +
\intO{\mathbb{S}(\nabla_h\vu):\Grad\phi } -\intO{ p_h(\vr) \Div \phi } 
 =<{ R}^M_{h,\Delta t},\phi>\;\mbox{in $(0,T]$}
$$
with any $\phi\in C^1([0,T]\times \Omega)$.
\end{enumerate} 
\eL
\bRemark{rm1}
Visiting the proof of Lemma \ref{Lm-consistency} we find that 
$$
\alpha_M\ge \min\Big\{\omega,\frac1 2-\frac\beta 4,\frac{2\gamma-6}\gamma, 
\frac{5\gamma-12}{6\gamma}\Big\}
$$
with the choice $q=1$.
\eR
{\bf Proof of Lemma \ref{Lm-consistency}} \\ 
We take
\[
\tilde\phi(t), \ \phi \in C^1_c([0,T]\times{\Omega}; R^3),
\]
as a test function in the discrete momentum equation (\ref{N3}). Seeing that, in accordance with (\ref{vv1}), (\ref{proj}),
\[
\intO{\mathbb{S}(\nabla_h\vu^k):\nabla_h\tilde\phi }=\intO{\mathbb{S}(\nabla_h\vv^k):\Grad\phi }
,\quad
\intO{ p(\vr^k) { \Divh} \tilde\phi } = \intO{ p(\vr^k) \Div \phi},
\]
we may rewrite  (\ref{N3})--by using  the formula (\ref{Up3})
and rearranging conveniently several terms-- and thus obtain the
first identity (\ref{me-v}) and then (\ref{me-v+}) integrating (\ref{me-v}) by parts,
where
$$
 <{R}^M_{h,\Delta t},\phi>
 = \intO{ D_t( \vr \widehat{\vv}) \cdot (\phi - \widehat{\tilde\phi}) } 
+ \sum_{K\in{\cal T}} \sum_{\sigma\in {\cal E}(K)\cap{\cal E}_{\rm int}}\int_\sigma (\phi - \widehat{\tilde\phi}) \cdot 
\ju{ \vr\widehat\vv }_{\sigma,\vc n_{\sigma,K}} \ [\vu_\sigma \cdot \vc{n}_{\sigma,K}]^- {\rm d}S_x 
$$
$$
+ \sum_{K \in {\cal T}} \sum_{\sigma\in {\cal E}(K)}\int_\sigma \vr(\widehat{\tilde\phi}-\phi)\cdot\widehat\vv (\vu -
\vu_\sigma )\cdot \vc{n}_{\sigma,K}{\rm d}S_x + \intO{ \vr(\phi - \widehat{\tilde\phi}) \cdot\widehat\vv  \Divh \vu }
$$
$$
%{\tb + \intO{\vr(\vu-\widehat\vu)\cdot\Grad\phi\cdot\widehat\vv}}
+ \intO{\vr\widehat\vu\cdot\nabla_h\vu_B\cdot(\phi-\widehat{\tilde\phi})}
{
 +\intO{\vr\widehat\vu\cdot\nabla\phi\cdot(\vu_B-\widehat\vu_B)}}
$$
$$
{-  h^\omega \sum_{\sigma\in{\cal E}_{\rm int}}\int_\sigma[[\vr^k]]_{\sigma}\{\widehat\vu^k\}_{\sigma}[[\widehat{\tilde\phi}]]_\sigma{\rm d}S_x}
+  \sum_{\sigma\in {\cal E}^{\rm in}}\int_\sigma
(\vr_B-\vr)\vu_{B,\sigma}\cdot \vc n_{\sigma}\widehat\vu\cdot (\phi - \widehat{\tilde\phi})=\sum_{i=1}^8I_i 
$$
We refer the reader to \cite[Formulas (7.1) and (7.2)]{KwNoNUM} for more details.

Our goal is to estimate conveniently all terms in $<R^M_{h,h},\phi>$.

\begin{enumerate}
\item {\it  Most terms in $<R^M_{h,\Delta t}(t),\phi>$ contain the expression $\phi-\widehat{\tilde\phi}$}. We notice that by virtue of (\ref{PoincareV-}), (\ref{errorV}), (\ref{globalV}),
$$
\|\phi-\widehat{\tilde\phi}\|_{L^q(\Omega)}\aleq h\|\Grad\phi\|_{L^q(\Omega)},\; 1\le q\le\infty.
$$
\item {\it Estimate of the term with the time derivatives ($I_1$)}:
For the error in the time derivative, we obtain
\bFormula{mest0}
\left| \intO{ D_t (\vr \widehat{\vv}) \cdot (\phi - \widehat{\tilde\phi} ) }\right| \aleq h^{\frac 12-\frac \beta 2}  A_h,\; \mbox{$A_h$ bounded in $L^2(0,T)$}.
\eF
Indeed, we split the term  $\intO{ D_t (\vr^k \widehat{\vv}^k) \cdot (\phi - \widehat{\tilde\phi} ) }$ into two parts,
\[
\intO{ \sqrt{ \vr^{k-1} } \sqrt{ \vr^{k-1}} \frac{ \vv^k - \vv^{k-1} }{\Delta t} \cdot ( \phi - \widehat{\tilde\phi}) }
+ \intO{ \frac{ \vr^k - \vr^{k-1} }{\Delta t} \vv^{k} \cdot ( \phi - \widehat{\tilde\phi}) },
\]
where, by virtue of H\"older's inequality, for any $\gamma>1$, the first term is bounded by
\[
\aleq
h (\Delta t)^{-1/2}\| \vr^{k-1} \|_{L^\gamma(\Omega)}^{1/2} \left(\Delta t \intO{ \vr^{k-1} \left( \frac{ \vv^{k - 1} - \vv^{k-1} }{\Delta t} \right)^2 } \right)^{1/2}
\left\| \Grad\phi  \right\|_{L^{\frac{2 \gamma}{\gamma - 1}}(\Omega)}
\]
while the second one is
\[
\aleq h^{1-\frac \beta 2}(\Delta t)^{-1/2}  \left(  \Delta t\, h^\beta \intO{ \left(\frac{ \vr^k - \vr^{k-1} }{ \Delta t } \right)^2 } \right)^{1/2} \| \vv^k \|_{L^6(\Omega;R^3)}
 \| \Grad \phi \|_{L^3(\Omega)},
\]
where the integrals of both expressions are controlled by means of (\ref{e4}).

\item {\it Estimate of the term with the jump of $\vr\widehat\vv$ ($I_2$)}:
This essentially amounts to estimate two terms,
$$
I_{2,1}= \sum_{K\in{\cal T}} \sum_{\sigma\in {\cal E}(K)\cap{\cal E}_{\rm int}}\int_\sigma \vr^+_{\sigma,\vc n_{\sigma,K}}(\phi - \widehat{\tilde\phi}) \cdot 
\ju{ \widehat\vv }_{\sigma,\vc n_{\sigma,K}} \ [\vu_\sigma \cdot \vc{n}_{\sigma,K}]^- {\rm d}S_x 
$$
and
$$
I_{2,2}= \sum_{K\in{\cal T}} \sum_{\sigma\in {\cal E}(K)\cap{\cal E}_{\rm int}}\int_\sigma (\phi - \widehat{\tilde\phi}) \cdot 
\widehat\vv^-_{\sigma,\vc n_{\sigma,K}}\ju{ \widehat\vr }_{\sigma,\vc n_{\sigma,K}} \ [\vu_\sigma \cdot \vc{n}_{\sigma,K}]^- {\rm d}S_x,
$$ 
where
$$
|I_{2,1}|\aleq \Big[\sum_{K\in{\cal T}} \sum_{\sigma\in {\cal E}(K)\cap{\cal E}_{\rm int}}\int_\sigma \vr^+_{\sigma,\vc n_{\sigma,K}}  
\ju{ \widehat\vv}^2_{\sigma,\vc n_{\sigma,K}} \ |[\vu_\sigma \cdot \vc{n}_{\sigma,K}]^-| {\rm d}S_x\Big]^{1/2}\times 
$$
\bFormula{I21}
\Big[\sum_{K\in{\cal T}} \sum_{\sigma\in {\cal E}(K)\cap{\cal E}_{\rm int}}\int_\sigma \vr^+_{\sigma,\vc n_{\sigma,K}}(\phi - \widehat{\tilde\phi})^2  \ |[\vu_\sigma \cdot \vc{n}_{\sigma,K}]^-| {\rm d}S_x \Big]^{1/2}
\eF
with the first term in the product controlled by (\ref{e5}) in $L^2((0,T))$, while
$$
|I_{2,2}|\aleq \Big[h^\beta\sum_{K\in{\cal T}} \sum_{\sigma\in {\cal E}(K)\cap{\cal E}_{\rm int}}\int_\sigma  
\ju{ \widehat\vr }_{\sigma,\vc n_{\sigma,K}}^2 \ |[\vu_\sigma \cdot \vc{n}_{\sigma,K}]^-| {\rm d}S_x\Big]^{1/2}\times
$$
\bFormula{I22}
\Big[h^{-\beta}\sum_{K\in{\cal T}} \sum_{\sigma\in {\cal E}(K)\cap{\cal E}_{\rm int}}\int_\sigma  (\phi - \widehat{\tilde\phi})^2 \cdot 
|\widehat\vv^-_{\sigma,\vc n_{\sigma,K}}|^2 \ |[\vu_\sigma \cdot \vc{n}_{\sigma,K}]^-| {\rm d}S_x\Big]^{1/2}
\eF
with the first term in the product controlled by (\ref{e6}), i.e. belonging to $L^2((0,T))$. It will be therefore enough to estimate
the expressions under the second square roots of $I_{2,1}$, $I_{2,2}$, respectively, and to consider only their "leading parts" (with $\vu$ replaced by $\vv$). 

We have, if $\gamma>3$, compare with (\ref{cest1}),
$$
\Big|\sum_{K\in{\cal T}} \sum_{\sigma\in {\cal E}(K)\cap{\cal E}_{\rm int}}\int_\sigma \vr^+_{\sigma,\vc n_{\sigma,K}}(\phi - \widehat{\tilde\phi})^2  
\ |[\vv_\sigma \cdot \vc{n}_{\sigma,K}]^-| {\rm d}S_x \Big|
$$
\bFormula{mest1}
\aleq \sum_{K\in{\cal T}} \sum_{\sigma\in {\cal E}(K)\cap{\cal E}_{\rm int}}\|\phi - \widehat{\tilde\phi}\|^2_{L^\gamma(\sigma)}
\|\vr\|_{L^\gamma(\sigma)}\|\vv\|_{L^{\frac{\gamma}{\gamma-3}}(\sigma)}
\eF
$$
\aleq h\sum_{K\in{\cal T}} \|\vr\|_{L^\gamma(K)}\|\vv\|_{L^{\frac{\gamma}{\gamma-3}}(K)}\|\Grad\phi\|_{L^\gamma(K)}^2\aleq
{ h^{1+\min\{0,\frac{3(\gamma-3)}\gamma-\frac 12\}}\|\vv\|_{L^6(\Omega)} \|\Grad\phi\|_{L^\gamma(\Omega)}^2}.
$$

Likewise
$$
h^{-\beta}\sum_{K\in{\cal T}} \sum_{\sigma\in {\cal E}(K)\cap{\cal E}_{\rm int}}\int_\sigma  (\phi - \widehat{\tilde\phi})^2 
|\widehat\vv^-_{\sigma,\vc n_{\sigma,K}}|^2 \ |[\vv_\sigma \cdot \vc{n}_{\sigma,K}]^-| {\rm d}S_x
$$
{
\bFormula{mest2}
\aleq h^{1-\beta}\|\vv\|^3_{L^{\frac{3\gamma}{\gamma-2}}(\Omega)}\|\Grad\phi|^2_{L^\gamma(\Omega)}
\eF
$$
\aleq  h^{1-\beta +\min\{0,\frac{3\gamma-6}\gamma-\frac 32\}} (\Delta t)^{-\frac 12\frac{1+\ep}{1+\ep/3}}
\|(\Delta t)^{\frac 16\frac{1+\ep}{1+\ep/3}}\vv\|_{L^{6}(\Omega)}^3\|\Grad\phi\|^2_{L^\gamma(\Omega)},\;\ep>0,
$$
where $\|(\Delta t)^{\frac 16\frac{1+\ep}{1+\ep/3}}\vv\|_{L^{6}(\Omega)}^3$ is bounded in $L^{1+\ep/3}(0,T)$ according to (\ref{timeG}).
Clearly, if $\gamma>3$ and $0<\beta< \min\{\frac 12,\frac{2\gamma-6}\gamma\}$ then $\ep>0$ can be chosen so small that $1-\beta +$ $\min\{0,\frac{3\gamma-6}\gamma-\frac 12\}$ 
$-\frac 12\frac{1+\ep}{1+\ep/3}>0$.}

\item {\it The upper bound of the third term ($I_3$)} is determined by the upper  bound of
$$
|\sum_{K \in {\cal T}} \sum_{\sigma\in {\cal E}(K)}\int_\sigma \vr(\phi - \widehat{\tilde\phi})\cdot\widehat\vv (\vv -
\vv_\sigma )\cdot \vc{n}_{\sigma,K}{\rm d}S_x|.
$$
If $\gamma>3$ it is bounded by
$$
\aleq
h \|\vr\|_{L^{\frac{3\gamma}{\gamma-3}}(\Omega)}\|\vv\|_{L^6(\Omega)}\|\nabla_h\vv\|_{L^2(\Omega)} \|\Grad\phi\|_{L^\gamma(\Omega)}
$$
\bFormula{mest3}
\aleq h^{1+\min\{0,\frac {\gamma-6}\gamma\}}(\Delta t)^{-\frac 12\frac\ep{1+\ep/2}} \|\vr\|_{L^\gamma(\Omega)}\|(\Delta t)^{\frac 14\frac\ep{1+\ep/2}}\nabla_h\vv\|^2_{L^2(\Omega)} \|\Grad\phi\|_{L^\gamma(\Omega)},\;\ep>0,
\eF
{ where $\|(\Delta t)^{\frac 14\frac\ep{1+\ep/2}}\nabla_h\vv\|^2_{L^2(\Omega)}$ is bounded in $L^\frac{2+\ep}2(0,T)$.}
\item {\it The bounds of term $I_4$} are determined by the bounds of
$$
\Big|\intO{ \vr(\phi - \widehat{\tilde\phi}) \cdot\widehat\vv  \Divh \vv }\Big|.
$$
They are exactly the same as in the previous case. The same is true for the terms $I_5$--$I_6$, since they have the same structure.

 \item The artificial viscosity  term $|I_7|= h^\omega |\sum_{\sigma\in{\cal E}_{\rm int}}$ $\int_\sigma[[\vr^k]]_{\sigma}\{\widehat\vv^k\}_{\sigma}[[\widehat\phi]]_\sigma{\rm d}S_x|$ is bounded by
\bFormula{mestav}
 h^{\omega}\|\vr\|_{L^3(\Omega)}\|\vu\|_{L^3(\Omega)}
\|\nabla\phi\|_{L^3(\Omega)}\aleq h^{\omega}\|\vr\|_{L^\gamma(\Omega)}\|\vu\|_{L^6(\Omega)}
\|\nabla\phi\|_{L^\gamma(\Omega)}.
\eF
\item The last term to be evaluated is the boundary term whose decay is determined by 
$$
I_8= 
\sum_{\sigma\in {\cal E}^{\rm in}}\int_\sigma
(\vr_B-\vr)\vu_{B,\sigma}\cdot \vc n_{\sigma}\widehat\vv\cdot (\widehat{\tilde\phi}-\phi) 
{\rm d}S_x
$$
We have, by the same reasoning (\ref{cest3}) for the similar term (again we write $|\vr_B-\vr|\le \vr_B+\vr$
and consider only the more difficult term with $\vr$ instead of $(\vr_B-\vr)$), without loss of generality):
If $\gamma\ge 12/5$, the absolute value admits the bound
\bFormula{mest4}
\aleq h^{\frac{5\gamma-12}{6\gamma}} \|\vr\|_{L^\gamma(\Omega)} \|\vv\|_{L^{6}(\Omega)} 
\|\Grad\phi\|_{L^\gamma(\Omega)}.
\eF
\end{enumerate}
Putting together (\ref{mest0}--\ref{mest4}) finishes the proof of Lemma \ref{Lm-consistency}.

\section{Improved estimates of density}\label{SEBOG}

As already mentioned, from the existence proofs of weak solutions in the continuous case  we know
that the uniform estimates derived in Lemma \ref{Lestimates} are not enough to pass to the limit
in the term containing $p(\vr)$ in the momentum equation. The first step in order to perform this task
is to improve the integrability of the pressure from $L^\infty(I;L^1(\Omega))$ to $L^p(I\times\Omega)$
with some $p>1$. We have to do the same also on the discrete level.

Mimicking the continuous case, we shall use the explicit solution known as Bogovskii's operator $\mathcal{B}$, see Bogovskii \cite{Bog}. Its properties are recalled in Lemmas \ref{LBog}, \ref{LBog+}.

The pressure estimates are obtained by taking
\[
\phi = \mathcal{B} \left[ \vr_h - \frac{1}{|\Omega|} \intO{ \vr_h } \right]
\]
in the discrete momentum equation (\ref{me-v})
($\phi\in L^\infty(I; W_0^{1,\gamma}(\Omega))$ by virtue of (\ref{Z2}) and (\ref{e0}), and it is hence admissible test function, cf. (\ref{NM1})):
\bFormula{Z4}
\int_0^T \intO{ p(\vr)\vr }\dt + h^\beta\int_0^T \intO{ \vr^{3} }\dt= \frac{1}{|\Omega|} \intO{ \vr } \int_0^T \intO{ p_h(\vr) }{ \dt}
\eF
$$
+
\intO{ \vr \widehat{\vv} \cdot \mathcal{B} \left[ \vr - \frac{1}{|\Omega|} \intO{ \vr } \right] (T, \cdot) } -
\intO{ \vr^0 \widehat{\vv}^0 \cdot \mathcal{B} \left[ \vr^0 - \frac{1}{|\Omega|} \intO{ \vr^0} \right] }
$$
\[
-  \int_0^T \intO{  \vr(t-\Delta t) \widehat{\vv}( t-\Delta t ) \cdot \mathcal{B} [D_t \vr(t) ] } \ \dt
- \int_0^T \left< R^M_{h,h} , \mathcal{B} \left[ \vr - \frac{1}{|\Omega|} \intO{ \vr } \right] \right> \ \dt
\]
\[
 - \int_0^T \intO{ (\vr {\vu} \otimes \widehat\vu) : \Grad  \mathcal{B} \left[ \vr - \frac{1}{|\Omega|} \intO{ \vr } \right]  } \ \dt +\int_0^T\intO{\vr\widehat\vu\cdot\nabla_h\Big(\vu_B\cdot\mathcal{B} \left[ \vr - \frac{1}{|\Omega|} \intO{ \vr } \right]\Big) }\dt
\]
\[
+ \int_0^T \intO{  \mathbb{S}(\Gradh \vu) : \Grad \ \mathcal{B} \left[ \vr - \frac{1}{|\Omega|} \intO{ \vr} \right]} \dt,
\]
where we have used the discrete "integration by parts"
\bFormula{dip}
\int_0^T \intO{ D_t (\vr \widehat{\vv})  \cdot \mathcal{B} \left[ \vr - \frac{1}{|\Omega|} \intO{ \vr } \right] } \ \dt= 
-  \int_0^T \intO{  \vr(t-\Delta t) \widehat{\vv}( t-\Delta t) \cdot \mathcal{B} [D_t \vr (t)] } \ \dt
\eF
$$
+\intO{ \vr_h \widehat{\vv} \cdot \mathcal{B} \left[ \vr - \frac{1}{|\Omega|} \intO{ \vr} \right] (T, \cdot) } -
\intO{ \vr^0 \widehat{\vv}^0 \cdot \mathcal{B} \left[ \vr^0 - \frac{1}{|\Omega|} \intO{ \vr^0 } \right] }.
$$

We observe that the expression on the right-hand side of (\ref{Z4}) is bounded uniformly for $h \to 0$. 
Indeed: {
Due to (\ref{JensenV}) and (\ref{e3}),
\bFormula{hatv}
\widehat\vv\;\mbox{is bounded in}\; L^2(I;L^6(\Omega)).
\eF 
}
By virtue of (\ref{e0}), (\ref{hatv}) (and, in combination eventually with (\ref{ruB})),
\bFormula{ru}
\vr\widehat\vv,\; \vr\widehat\vu \;\mbox{are bounded in}\; L^\infty(I;L^{\frac{2\gamma}{\gamma+1}}(\Omega))\cap L^2(I;L^{\frac{6\gamma}{\gamma+6}}(\Omega)).
\eF
Further, 
\bFormula{ruu}
\vr\vu\otimes\widehat\vu, \,\vr\widehat\vu\otimes\widehat\vu,\; \mbox{is bounded in}\; L^2(I;L^{\frac{6\gamma}{4\gamma+3}}(\Omega))\cap
L^1 (I;L^{\frac{3\gamma}{\gamma+3}}(\Omega)).
\eF
Moreover, thanks to  { (\ref{e0}), (\ref{NM1}), (\ref{Z2})},
\[
\int_0^T \left< R^M_{h,h} , \mathcal{B} \left[ \vr_h - \frac{1}{|\Omega|} \intO{ \vr_h } \right] \right> \ \dt \to 0
\ \mbox{as}\ h \to 0.
\]
and by (\ref{ce-v}) and (\ref{Z3})
$$
{\cal B}[D_t\vr]\;\mbox{is bounded in}\; L^2(I;L^{\frac{6\gamma}{5\gamma-6}}(\Omega)).
$$ 
Employing these facts and the H\"older inequality in each term at the right hand side of (\ref{Z4}), we show that it is bounded. As a conclusion, we have,
\bFormula{rg+1}
\|\vr\|_{L^{\gamma+1}(I\times\Omega)}\aleq 1,\quad \|h^{\beta/3}\vr\|_{L^3(I\times\Omega)}\aleq 1.
\eF

\section{Convergence }\label{SConv}

We denote by $[\vr_h,\vu_h=\vv_h+\tilde{\mathfrak{u}}_B]$, $h>0$ a sequence of numerical solution to  the scheme (\ref{N1}--\ref{N3})
where we extend $\vv_h$ by $0$ outside $\Omega$. We want to show that
 there is a  subsequence with weak limit $[{\mathfrak{r}},\mathfrak{u}]$, such that the couple $[{\mathfrak{r}},\mathfrak{u}]$ is a 
weak solution of the continuous problem (\ref{NS1}--\ref{inB}) in the sense of Definition
\ref{DD1}.
\subsection{Weak limits. Continuity equation}
Recalling regularity (\ref{ru0}--\ref{ruB}) of the initial and boundary data, we deduce from (\ref{errorQ}--\ref{PoincareV-}), (\ref{globalV}) (\ref{errorV}),
\bFormula{co0-}
\vr_{B,h}\to\mathfrak{r}_B\;\mbox{in $L^q(\Omega)$ and $L^q(\partial\Omega)$, $1\le q\le\infty$},
\eF 
$$
\widehat\vu_{B,h}\to \mathfrak{u}_B,\;\vu_{B,h}\to\mathfrak{u}_B,\;\nabla_h\vu_{B,h}\to\Grad\mathfrak{u}_B\;\mbox{in $L^q(\Omega)$, $1\le q\le \infty$}, 
$$
$$
\vr^0_h\to \mathfrak{r}_0\;\mbox{in $L^q(\Omega)$}, \;\vu^0_h\to\mathfrak{u}_0\;\;\mbox{in $L^q(\Omega)$},\;
1\le q\le\infty.
$$

Recalling (\ref{e3}) and (\ref{hatv}), we infer\footnote{All convergences in this section hold for a chosen subsequences of the original sequence; for the sake of simplicity, we do not relabel.}
\bFormula{co1+}
\widehat\vv_h\rightharpoonup\mathfrak {v},\;
\vv_h\rightharpoonup\mathfrak{v}
\;\mbox{in $L^2(0,T;L^6(\Omega))$},
\eF
where the limit of both sequences is the same by virtue of (\ref{globalV}). { Further},
%We have by the second inequality in (\ref{jumps}),
$$
\intO{\vv_h\Grad\phi}=-\intO{\nabla_h\vv_h\phi}+I^{1}_h +I^{2}_h,\;\phi\in C^1((0,T)\times\overline\Omega)),
$$
where
{
$$
I^{1}_h:=\sum_{K\in {\cal T}}\sum_{\sigma\in {\cal E}(K)}\int_\sigma\vv_h\vc n_{\sigma,K}\phi{\rm d}S_x=
\sum_{K\in {\cal T}}\sum_{\sigma\in {\cal E}(K)\cap{\cal E}_{\rm int}}\int_\sigma(\vv_h-\vv_{h,\sigma})\vc n_{\sigma,K}(\phi-\phi_\sigma){\rm d}S_x,\;
I^{2}_h=
 \int_{\partial\Omega}\vv_h\vc n\phi{\rm d}S_x
$$
admit the bounds
\bFormula{dod*}
 |I^{1}_h|\aleq h\|\nabla_h\vv_h\|_{L^2(0,T;L^2(\Omega))}\|\phi\|_{L^2(0,T;W^{1,2}(\Omega))}, \;
|I^2_h|\aleq h^{1/3}\|\nabla_h\vv_h\|_{L^2(0,T;L^2(\Omega))}\|\phi\|_{L^2(0,T;W^{1,2}(\Omega))}
\eF
by virtue of the H\"older inequality, trace estimates (\ref{trace}--\ref{traceFE}), 
%first inequality in (\ref{jumps}) 
the first inequality in (\ref{PoincareV-}), and in the second estimate also
the standard Sobolev imbedding $W^{1,2}(\Omega)\hookrightarrow L^6(\Omega)$,  the fact that $\vv_h\in V_0(\Omega;R^3)$ and that $|\cup_{K\cap {\cal E}_{\rm ext}\neq \emptyset} K|\aleq h.$}
We deduce from this and from (\ref{e3}), 
\bFormula{co2}
\nabla_h\vv_h\rightharpoonup\Grad\mathfrak{v}\;\mbox{in $L^2(0,T;L^2(\Omega))$},\;\mbox{and $\mathfrak{v}\in
L^2(0,T;W^{1,2}_0(\Omega))$.}
\eF

It is the consequence of (\ref{pwt+}) and (\ref{e0}), (\ref{e3}), (\ref{rg+1}) resp. (\ref{ru}), (\ref{hatv}), and (\ref{globalV}) that, in particular,
\bFormula{etildevr}
\utilde{\vr_h}\;
\mbox{ is bounded in}\; L^\infty(0,T;L^\gamma(\Omega))\cap L^{\gamma+1}((0,T)\times\Omega),
\eF
$$
\utilde{\vr_h}\widehat\vv_h,\;
\utilde{\vr_h\widehat\vv_h}\;
\mbox{ is bounded in}\; L^\infty(0,T;L^{\frac{2\gamma}{\gamma+1}}(\Omega))\cap L^2(0,T;L^{\frac{6\gamma}{\gamma+6}}(\Omega)),
$$
$$
\utilde{\vr_h}\vv_h \;
\mbox{ is bounded in}\; L^2(0,T;L^{\frac{6\gamma}{\gamma+6}}(\Omega)),
$$
$$
\utilde{\vr_h\widehat\vv_h}\vv_h,\; \utilde{\vr_h\widehat\vv_h}\widehat\vv_h\;\mbox{is bounded in $L^2(0,T; L^{\frac{6\gamma}{4\gamma+3}}(\Omega))$}.
$$
We make a little detour and calculate by using  (\ref{pwt+}) and (\ref{e4}),
\bFormula{interpol}
\|\utilde{\vr_h}-\vr_h\|_{L^1((0,T)\times\Omega)} \aleq \int_0^T\|{\vr_h}(t-\Delta t)-\vr_h(t)\|_{L^1(\Omega)}
{\rm d}t \aleq (\Delta t)^{1/2}h^{-\frac \beta 2},
\eF
$$
\|\utilde{\vr_h\widehat\vu_h}-\vr_h\widehat\vu_h\|_{L^1((0,T)\times\Omega)}
\aleq \int_0^T\|{\vr_h}\widehat\vu_h(t-\Delta t)-\vr_h\widehat\vu_h(t)\|_{L^1(\Omega)}{\rm d}t
\aleq (\Delta t)^{1/2}h^{-\frac \beta 2}.
$$
We also recall the interpolation formula
\bFormula{interpol+}
\|z\|^2_{L^2(0,T;L^{2}(\Omega))}\aleq \|z\|^{\alpha}_{L^\infty(0,T;L^{1}(\Omega))}\,\|z\|^{\alpha}_{L^1(0,T;L^{1}(\Omega))}
\|z\|^{2(1-\alpha)}_{L^2(0,T;L^{\frac{6\gamma}{\gamma+6}}(\Omega))},\;\frac 12=\alpha+\frac{(1-\alpha)(\gamma+6)}{6\gamma}.
\eF
This will become useful later.

Let us come back with (\ref{e0}), (\ref{e1}) and (\ref{ru}) to the continuity equation (\ref{ce-v+}). We find that for all $\phi\in C^1_c(\Omega)$,
$$
\intO{\utilde{\vr_h}(t)\phi}=A^\phi_h(t)+B^\phi_h(t),
$$
where 
$$
\tau\mapsto
A^\phi_h(\tau)=\intO{{\vr}^0\phi(0)}+\int_{0}^\tau\intO{\vr\widehat\vu\cdot\Grad\phi}{\rm d}t
$$
is equi-bounded and equi-continuous in $C([0,T])$, while
$$\Big(\tau\mapsto
B^\phi_h(\tau)=\int_{0}^{\tau}< R^C_{h}(t),\phi>{\rm d}t \Big)\to 0\;\mbox{in $C([0,T])$}
$$
by virtue of (\ref{NC1}). Consequently,
by density of $C^1_c(\Omega)$ in $L^{\gamma'}(\Omega)$ and by the Arzela-Ascoli type argument, we get
\bFormula{co4}
\utilde{\vr_h}\to\mathfrak{r}\;\mbox{in $C_{\rm weak}([0,T];L^\gamma(\Omega))$ and in $L^a(0,T; W^{-1,2}(\Omega))$, $1\le a<\infty$},
\eF
where we have used the compact imbedding $L^\gamma(\Omega)\hookrightarrow W^{-1,2}(\Omega)$, estimate (\ref{e0}) and the Lebesgue dominated convergence theorem to deduce from the first converegence relation the second one.
Recalling the bounds (\ref{ru}), (\ref{etildevr}) and writing 
$$
\vv_h=[\vv_h]_{h}+ (\vv_h-[\vv_h]_{h}),
$$
with
$[\cdot]_h$ denoting the standard mollification with a regularizing kernel, and realizing that, according to (\ref{regV})
\bFormula{conv1}
 [\vv_h]_{h}\rightharpoonup\mathfrak{v}\;\mbox{in $L^2(0,T;W^{1,2}(\Omega))$},\quad \|\vv_h-[\vv_h]_{h}\|^2_{L^2(\Omega)}
\aleq h^2\|\Gradh\vv_h\|^2_{L^2(\Omega)}
\eF  
we infer, seeing the last line in (\ref{etildevr})
\bFormula{crv}
\utilde{\vr_h}\vv_h\rightharpoonup_*\mathfrak{r}\mathfrak{v}\;\mbox{in
$L^2(0,T;L^{\frac{6\gamma}{\gamma+6}}(\Omega))$},
\eF
Employing   (\ref{ru}), (\ref{co0-}),  (\ref{etildevr}), (\ref{interpol}), (\ref{interpol+}) with $z=\utilde{\vr_h}-\vr_h$ and
(\ref{globalV}), we deduce from (\ref{crv}),
\bFormula{cru}
{\vr_h}\widehat\vu_h\rightharpoonup_*\mathfrak{r}\mathfrak{u}\;\mbox{in
$L^\infty(0,T;L^{\frac{2\gamma}{\gamma+1}}(\Omega))$},\;\mathfrak{u}=\mathfrak{v}+\mathfrak{u}_B.
\eF
 
Last but not least, due to estimate (\ref{e3}),
\bFormula{cru+}
\vr_h\rightharpoonup\mathfrak{r}\;\mbox{in $L^\gamma(0,T;L^\gamma(\partial\Omega,|\vu_B\cdot\vc n|{\rm d}S_x))$}.
\eF
Finally,
due to Lemma \ref{Lc-consistency}
$$
\int_0^\tau<R^C_{h},\phi>{\rm d}t\to 0\;\mbox{as $h\to 0$}.
$$

At this stage we can pass to the limit in the consistency formulation (\ref{ce-v+}) of the continuity equation
 in order to obtain the weak formulation (\ref{D2}). 

\subsection{Convergence in the momentum equation starts}

According to (\ref{pres2}), (\ref{rg+1}),
\bFormula{co5}
p(\vr_h)\rightharpoonup \overline{p(\mathfrak{ r})}\;\mbox{in $L^{\frac{\gamma+1}\gamma}(I\times\Omega))$},
\eF
where here and in the sequel $\overline {g(\mathfrak{r},\mathfrak{v},\nabla\mathfrak{v})}$ denotes
a weak limit of the sequence $g(\vr_h,\vv_h,\nabla_h\vv_h)$ (in $L^1(I\times\Omega)$).
Coming back with (\ref{e2}--\ref{e3}), (\ref{NM1}), (\ref{ruu}), (\ref{co5})   to the momentum equation (\ref{me-v+}), we obtain by the same Arzela-Ascoli type arguments -- compare with (\ref{co4})--
\bFormula{co6}
\utilde{\vc q}_h\to{\mathfrak{q}}\;\mbox{in $C_{\rm weak}([0,T];L^{\frac{2\gamma}{\gamma+1}}(\Omega))$ and in $L^a(0,T;W^{-1,2}(\Omega))$, $1\le a<\infty$},\;
\vc q_h=\vr_h\widehat\vv_h
\eF
where, due to (\ref{cru}) and (\ref{interpol}),
$$
{\mathfrak{q}}=\mathfrak{r}\mathfrak{v}\;\mbox{a.e. in $(0,T)\times\Omega$}.
$$
Evidently, due to (\ref{co0-}) and (\ref{co6}),
\bFormula{co6+}
{\vc m}_h:=\vr_h\widehat\vu_h\to{\mathfrak{m}}\;\mbox{in $C_{\rm weak}([0,T];L^{\frac{2\gamma}{\gamma+1}}(\Omega))$},
\eF
where
$$
{\mathfrak{m}}=\mathfrak{r}\mathfrak{u}\;\mbox{a.e. in $(0,T)\times\Omega$}\,\mbox{and}\,\mathfrak{u}=\mathfrak{v}+\mathfrak{u}_B.
$$
Now, reasonning as in (\ref{conv1}), (\ref{crv}) we get,
\bFormula{cruu-}
\utilde{\vr_h\widehat\vu_h}\otimes\vu_h\rightharpoonup
\mathfrak{r}\mathfrak{u}\otimes\mathfrak{u}\;\mbox{in $L^2(0,T;L^{\frac{6\gamma}{4\gamma+3}}(\Omega))$},
\eF
where we have used also the last line in (\ref{etildevr}). Finally, by virtue of (\ref{ru}), (\ref{ruu}) the last line in (\ref{etildevr}), (\ref{interpol}), (\ref{interpol+}) with $z= \utilde{\vr_h\widehat\vu_h}-{\vr_h\widehat\vu_h}$ and (\ref{globalV}), we deduce from (\ref{cruu-})
\bFormula{cruu}
{\vr_h\widehat\vu_h}\otimes\widehat\vu_h\rightharpoonup
\mathfrak{r}\mathfrak{u}\otimes\mathfrak{u}\;\mbox{in $L^2(0,T;L^{\frac{6\gamma}{4\gamma+3}}(\Omega))$}.
\eF

Now, we are ready to pass to the limit in the momentum equation (\ref{me-v+}). We obtain
\bFormula{me-d}
\intO{{\mathfrak{q}}\cdot\phi(\tau,x)}-\intO{\mathfrak{r}_0\mathfrak{v}_0\cdot\phi(0,x)}
\eF
$$
= \int_0^\tau\intO{\Big(\mathfrak{r}\mathfrak{v}\cdot\partial_t\phi+\mathfrak{r} { \mathfrak{u} } \otimes \mathfrak{u}:\Grad\phi+\overline{p(\mathfrak{r})}\Div\phi -\mathfrak{r}\mathfrak{u}\cdot\nabla(\mathfrak{u}_B\cdot\phi)- \mathbb{S}(\Grad\mathfrak{u}):\Grad\phi\Big)}{\rm d}t,
$$
where we have used also (\ref{NM1}).
According to (\ref{D2}),
$$
\int_0^\tau\intO{\mathfrak{r}\mathfrak{u}\cdot\Grad(\mathfrak{u}_B\cdot\phi)}{\rm d}t=-\int_0^\tau\intO{\mathfrak{r}\mathfrak{u}_B\cdot\partial_t\phi}{\rm d}t
+\Big[\intO{\mathfrak{r}\mathfrak{u}_B\cdot\phi}\Big]_0^\tau
$$
we obtain from (\ref{me-d}) the formulation
\bFormula{D3*}
\intO{{\mathfrak{m}}\cdot\phi(\tau,x)}-\intO{\mathfrak{r}_0\mathfrak{u}_0\cdot\phi(0,x)}=
\eF
$$
= \int_0^\tau\intO{\Big(\mathfrak{r}\mathfrak{u}\cdot\partial_t\phi+\mathfrak{r} { \mathfrak{u} } \otimes \mathfrak{u}:\Grad\phi+\overline{p(\mathfrak{r})}\Div\phi - \mathbb{S}(\Grad\mathfrak{u}):\Grad\phi\Big)}{\rm d}t.
$$ 

\subsection{Effective viscous flux}\label{73}
\subsubsection{An integration by parts formula}
We denote ${\rm curl}\phi:=\nabla\phi-\nabla^T\phi$. It is easy to verify that,
for all $\vc w\in V(\Omega;R^3)$, $\phi\in W^{1,2}(\Omega;R^3)$,
\bFormula{intpp}
2\intO{\nabla_h\vc w:\Grad\phi}=\intO{\Big({\rm curl}_h\vc w:{\rm curl}_x\phi
+2{\rm div}_h\vc w\Div\phi\Big)}+J_h[\vc w,\phi],
\eF
where
$$
J_h[\vc w,\phi]=2\intO{\Big(\nabla_h\vc w:\Grad^T\phi-{\rm div}_h\vc w\Div\phi\Big)}.
$$
We have on one side,
\bFormula{est0}
|J_h[\vc w,\phi]|\aleq\|\nabla_h\vc w\|_{L^2(\Omega)}\|\Grad\phi\|_{L^2(\Omega)}
\eF
simply by the Cauchy-Schwartz inequality, and on the other side
\bFormula{est1}
|J_h[\vc w,\phi]|\aleq h\|\nabla_h\vc w\|_{L^2(\Omega)}\|\Grad\phi\|_{W^{1,2}(\Omega)}\; 
\mbox{provided $\vc w\in V_0(\Omega;R^3)$, $\phi\in W^{2,2}(\Omega)$}.
\eF
%and
%\bFormula{est2}
%|J_h[\tilde\vc w,\phi]|\aleq h\|\Grad\phi\|_{L^2(\Omega)}\|\Grad^2{\vc w}\|_{L^2(\Omega)},
%\mbox{provided $\phi\in W^{1,2}_0(\Omega)$, $\vc w\in W^{2,2}(\Omega)$}.
%\eF

Indeed, to get (\ref{est1}), we calculate by using integration by parts,
$$
\Big|J_h[\vc w,\phi]\Big|=\Big|\sum_{K\in {\cal T}}\sum_{\sigma\in {\cal E}(K)}\int_\sigma
\Big(\vc w\cdot\Grad\phi\cdot\vc n_{\sigma,K}+\vc w\cdot\vc n_{\sigma,K}\Div\phi\Big)
{\rm d}S_x\Big|
$$
$$
= \Big|\sum_{K\in {\cal T}}\sum_{\sigma\in {\cal E}(K)}\int_\sigma
\Big((\vc w-\vc w_\sigma)\cdot(\Grad\phi-(\Grad\phi)_\sigma)\cdot\vc n_{\sigma}+(\vc w-\vc w_\sigma)\cdot\vc n_{\sigma}(\Div\phi-(\Div\phi)_\sigma)\Big)
{\rm d}S_x\Big|
$$
$$
\aleq h^2 \sum_{K\in {\cal T}}\sum_{\sigma\in {\cal E}(K)}
\|\vc w-\vc w_\sigma\|_{L^2(\sigma)}\|\Grad\phi-(\Grad\phi)_\sigma\|_{L^2(\sigma)}
\aleq 
h \|\nabla_h\vc w\|_{L^2(\Omega)}\|\Grad\phi\|_{W^{1,2}(\Omega)},
$$
where we have used the fact that $\vc w_\sigma$ and $\Grad\phi$ are continuous over internal faces, $\vc w_\sigma=0$ if $\sigma\in {\cal E}_{\rm ext}$ and (\ref{PoincareV-}), (\ref{trace}--\ref{traceFE}).

%For (\ref{est2}) we employ (\ref{proj}) and then integration by parts, to write, by the same token,
%$$
%\Big|J_h[\tilde\vc w,\phi]\Big|= \Big|J_h[\vc w,\tilde\phi]\Big|=
%\Bi|\sum_{K\in {\cal T}}\sum_{\sigma\in {\cal E}_{\rm ext}\cap {\cal E}(K)}\int_\sigma\Big(\tilde\phi\cdot
%(\Grad\ww-(\Grad\ww)_\sigma)\cdot\vc n-\vc n_{\sigma}\cdot\tilde\phi
%(\Div\vc w-(\Div\vc w)_\sigma)\Big){\rm d}S_x\Big|
%$$
%$$
%\aleq \sum_{K\in {\cal T}}\|\Grad^2\vc w\|_{L^2(K)}\|\tilde\phi\|_{L^2(K)}\aleq
% \|\Grad^2\vc w\|_{L^2(\Omega)}(\|\tilde\phi-\phi\|_{L^2(\Omega)}+ \|\phi\|_{L^2(K)})
%\aleq h \|\Grad^2\vc w\|_{L^2(\Omega)} \|\Grad\phi\|_{L^2(\Omega)}.
%$$

At the point of conclusion, we deduce from (\ref{est0}--\ref{est1}),
\bFormula{intppfinal}
\intO{\mathbb{S}(\nabla_h\vu):\Grad\phi}=
\frac\mu 2\intO{{\rm curl}_h\vu:{\rm curl}_x\phi}
+ (2\mu+\lambda)\intO{{\rm div}_h\vu \Div\phi} +\frac\mu 2 J_h[\vu,\phi]
\eF
for all $\phi\in W^{1,2}_0(\Omega)$, where $\vu=\vv+\vu_B$ and
\bFormula{A1}
|J_h[\vu,\phi]|\aleq \|\nabla_h\vu\|_{L^2(\Omega)}\|\Grad\phi\|_{L^2(\Omega)},
\eF
\bFormula{A2}
|J_h[\vu,\phi]|\aleq h \|\nabla_h\vu\|_{L^2(\Omega)}\|\Grad\phi\|_{W^{1,2}(\Omega)}\;\mbox{if moreover $\phi\in W^{2,2}(\Omega)$}.
\eF

\subsubsection{Testing of momentum equation}
\noindent
{\it 1. Testing on the discrete level}
\\ \\

 We use in (\ref{me-v}) the functions
\bFormula{dod23}
\phi = \varphi \Grad \Delta^{-1} [\vr_h],\;\varphi(t,x) = \psi (t) \eta(x) , \psi \in \DC(0,T),\ \eta \in \DC(\Omega)
\eF
where the operator $\Grad \Delta^{-1}$ is defined in (\ref{calA}) (and $\vr_h$ si tacitly extended by $0$ outside $\Omega$) as test functions in the momentum equation
(\ref{me-v}). Due to Lemmas
\ref{LcalA}--\ref{LcalA+}, \ref{Lm-consistency} and estimate (\ref{e0}), $\phi
\in L^\infty(0,T; W^{1,\gamma}_0(\Omega))$ is an admissible test function.
We get
\bFormula{Z29}
\int_0^T \intO{  \Big[\varphi p_h(\vr_h) \vr_h - \mathbb{S}(\nabla_h\vu_h):\Grad\phi\Big] } \ \dt 
\eF
$$
=
-\int_0^T \intO{  p_h(\vr_h)\Grad\varphi\cdot\Grad\Delta^{-1}[ \vr_h] } \ \dt
- \int_0^T \left<{ R}^M_{h,h}, \varphi \Grad \Delta^{-1} [\vr_h]  \right> \ \dt
$$
\[
- \int_0^T \intO{ (\vr_h {\vu}_h \otimes \widehat\vv_h) : \Grad \left( \varphi \Grad \Delta^{-1} [\vr_h] \right) } \ \dt
+ \int_0^T \intO{\varphi D_t (\vr_h \widehat{\vv}_h) \cdot \Grad \Delta^{-1} [\vr_h]  } \ \dt
\]
$$
{+\int_0^T\intO{\vr_h(\widehat\vu_h-\vu_h)\cdot\Grad\phi\cdot\vu_B}{\rm d}t}
+\int_0^T\intO{\vr_h\widehat\vu_h\cdot\Grad\phi\cdot(\vu_B-\widehat\vu_B)}{\rm d}t+
\int_0^T\intO{\vr_h\widehat\vu_h\cdot\nabla_h\vu_B\cdot\phi}{\rm d}t.
$$

Furthermore, employing formula (\ref{intppfinal}) at the left-hand side, and the integration by parts formula 
\[
\int_0^T \intO{\varphi D_t (\vr_h \widehat{\vv}_h) \cdot \Grad \Delta^{-1} [\vr_h]  } \ \dt = -
\int_0^T \intO{ \frac{ \varphi (t) - \varphi (t-\Delta t) }{\Delta t} (\vr_h \widehat{\vv}_h)(t-\Delta t) 
\cdot
\Grad \Delta^{-1} [\vr_h(t)]  } \ \dt
\]
\bFormula{ippt}
- \int_0^T \intO{ (\varphi \vr_h  \widehat{\vv}_h) (t- \Delta t) \cdot \Grad \Delta^{-1} [D_t \vr_h(t)] } \ \dt
\eF
\[
= -
\int_0^T \intO{ \frac{ \varphi (t) - \varphi (t-\Delta t) }{\Delta t} (\vr_h \widehat{\vv}_h)(t-\Delta t) \cdot
\Grad \Delta^{-1} [\vr_h(t)] } \ \dt
\]
\[
+ \int_0^T \intO{ (\varphi \vr_h \widehat{\vv}_h)(t- \Delta t) \cdot \Grad \Delta^{-1} \Div ((\vr_h \vu_h)(t)) } \ \dt
\]
\[
+ \int_0^T \intO{ (\varphi \vr_h  \widehat{\vv}_h)(t- \Delta t) \cdot \Grad \Delta^{-1} [ { R}^C_h (t)] } \ \dt,
\]
at the right hand side, we obtain
\bFormula{Z30}
\int_0^T \intO{  \varphi\Big[p_h(\vr_h) \vr_h - (2\mu+\lambda)\vr_h{\rm div}_h\vu_h \Big] } \ \dt 
=
\mathfrak{K}[\vr_h,\vu_h] + \mathfrak{O}[\vr_h,\vu_h]+ \mathfrak{C}[\vu_h,\vr_h,\vr_h\widehat\vv_h],
\eF
where
$$
\mathfrak{ K}[\vr_h,\vu_h]
=
-\int_0^T \intO{\Big[  p_h(\vr_h)\mathbb{I}
-\mathbb{S}(\nabla_h\vu_h) \Big]: \Big(\Grad\varphi\otimes\Grad\Delta^{-1}[ \vr_h]\Big)
} \ \dt
$$
$$
 -
\int_0^T \intO{ \frac{ \varphi (t) - \varphi (t-h) }{h} (\vr_h \widehat{\vv}_h)(t-h) \cdot
\Grad \Delta^{-1} [\vr_h(t)]  } \ \dt
$$
$$
- \int_0^T \intO{ (\vr_h {\vu}_h \otimes \widehat\vv_h) : \left(\Grad \varphi \otimes\Grad \Delta^{-1} [\vr_h(t)] \right) } \ \dt +\int_0^T\intO{\vr_h\widehat\vu_h\cdot\nabla_h\vu_B\cdot\phi}{\rm d}t
$$
and
$$
\mathfrak{ O}[\vr_h,\vu_h]=
\frac\mu 2 \int_0^T J_h[\vu_h,\phi]  {\rm d}t
- \int_0^T \left< { R}^M_{h,h}, \varphi \Grad \Delta^{-1} [\vr_h] \right> \ \dt
$$
$$
+  \int_0^T \intO{ (\varphi \vr_h \widehat{\vv}_h)(t- h) \cdot \Grad \Delta^{-1} [ R^C_h(t) ] } \ \dt+\int_0^T\intO{\vr_h\widehat\vu_h\cdot\Grad\phi\cdot(\vu_B-\widehat\vu_B)}{\rm d}t
$$
{
%$$
%+\int_0^T\intO{\varphi\Big((\vr_h\widehat\vv_h)(t-h)-(\vr_h\widehat\vv_h)(t)\Big)
%\cdot\Grad{\Delta}^{-1}[\vr_h]}{\rm d}t
%$$
$$
+\int_0^T\intO{\vr_h(\widehat\vu_h-\vu_h)\cdot\Grad\phi\cdot\vu_B}{\rm d}t
+ \int_0^T \intO{ \Big((\varphi \vr_h \widehat{\vv}_h)(t- h) -
(\varphi\vr_h\widehat{\vv}_h)(t)\Big)\cdot \Grad \Delta^{-1} \Div (\vr_h \vu_h) } \ \dt
$$
}
and
$$
 \mathfrak{C}[\vu_h,\vr_h,\vr_h\widehat\vv_h]=
\int_0^T \intO{ \varphi\vu_h\cdot\Big[\vr_h\Grad\Delta^{-1}\Grad\cdot[\vr_h\widehat\vv_h]-(\vr_h\widehat\vv_h\cdot\Grad\Delta^{-1}\Grad)[\vr_h]\Big]} \ \dt.
$$
In the very last indentity of the integration by parts formula (\ref{ippt}), we have used the formulation (\ref{ce-v}) of the continuity equation in order to express
$D_t\vr_h$. This can be justified by using Lemmas
 \ref{LcalA}--\ref{LcalA+}.
\\ \\
{\it 2. Testing of  the limiting momentum equation}
\\ \\
 Now, we use the function
$$
\phi = \varphi \Grad \Delta^{-1} [\mathfrak{r}],\;\varphi(t,x) = \psi (t) \eta(x) , \psi \in \DC(0,T),\ \eta \in \DC(\Omega)
$$
as test function in the momentum equation
(\ref{D3*}).  Reasoning in the similar manner as we did when deriving (\ref{Z29}--\ref{Z30}),
we arrive at
\bFormula{Z31}
\int_0^T \intO{  \Big[\varphi \overline{p(\mathfrak{r} )} \mathfrak{r} 
- (2\mu+\lambda)\mathfrak{r}{\rm div}_x\mathfrak{u} \Big] } \ \dt 
=
\underline{\mathfrak{K}}[\mathfrak{r},\mathfrak{u}] + \underline{\mathfrak{C}}[
\mathfrak{u},\mathfrak{r},\mathfrak{r}\mathfrak{v}],
\eF
where
$$
\underline{\mathfrak{ K}}[\mathfrak{r},\mathfrak{u}]
=
-\int_0^T \intO{\Big[  \overline{p(\mathfrak{r})}\mathbb{I}
-\mathbb{S}(\Grad\mathfrak{u}) \Big]: \Grad\varphi\otimes\Grad\Delta^{-1}[ \mathfrak{r}]
} \ \dt
 -
\int_0^T \intO{ \partial_t\varphi \mathfrak{r} \mathfrak{v} \cdot
\Grad \Delta^{-1} [\mathfrak{r}]  } \ \dt
$$
$$
- \int_0^T \intO{ (\mathfrak{r} \mathfrak{u}\otimes \mathfrak{v}) : \left(\Grad \varphi \otimes\Grad \Delta^{-1} [\mathfrak{r}] \right) } \ \dt+ \int_0^T \intO{\mathfrak{r}\mathfrak{u}\cdot
\Grad\mathfrak{u}_B\cdot\phi}{\rm d}t
$$
and
$$
 \underline{\mathfrak{C}}[\mathfrak{u},\mathfrak{r}, \mathfrak{r}\mathfrak{v}]=
\int_0^T \intO{ \varphi\mathfrak{u}\cdot\Big[\mathfrak{r}\Grad\Delta^{-1}\Grad\cdot[\mathfrak{r}\mathfrak{v}]-(\mathfrak{r}\mathfrak{v}\cdot\Grad\Delta^{-1}\Grad)[\mathfrak{r}])\Big]} \ \dt
$$
{\it 3. Effective viscous flux identity}
\\ \\
{\it Step 1: Convergence of $\mathfrak{ C}[\vu_h, \vr_h,\vr_h\widehat\vv_h]$.} 
We set
$$
D_h:=D_h[{\vr_h},{\vr_h\widehat\vv_h}]:=\Big[{\vr_h}\Grad\Delta^{-1}\Grad\cdot[{\vr_h\widehat\vv_h}]-({\vr_h\vv_h}\cdot\Grad\Delta^{-1}\Grad)[{\vr_h}]\Big],\quad \utilde{D_h}=D_h[\utilde{\vr_h},\utilde{\vr_h\widehat\vv_h}]
$$ 
$$
{D}= D[\mathfrak{r},\mathfrak{q}]:= \Big[\mathfrak{r}\Grad\Delta^{-1}\Grad\cdot[\mathfrak{q}]-(\mathfrak{q}\cdot\Grad\Delta^{-1}\Grad)[\mathfrak{r}])\Big],\; \mathfrak{q}=
\mathfrak{r}\mathfrak{v}.
$$
The goal is to show that
\bFormula{limC}
\int_0^T\intO{\varphi{\vv_h}\cdot D_h}{\rm d}t\to \int_0^T\intO{\varphi{\mathfrak{v}}\cdot D}{\rm d}t.
\eF
By virtue of (\ref{etildevr}) and in view of (\ref{interpol}--\ref{interpol+}), the limits
of $\int_0^T\intO{\varphi{\vv_h}\cdot D_h}{\rm d}t$ and $\int_0^T\intO{\varphi{\vv_h}\cdot\utilde{D_h}}{\rm d}t$
are the same. It is therefore enough to show
$$
\int_0^T\intO{\varphi{\vv_h}\cdot\utilde{D_h}}{\rm d}t\to \int_0^T\intO{\varphi{\mathfrak{v}}\cdot D}{\rm d}t.
$$
We shall proceed in several steps:
\begin{enumerate}
\item We shall write
$$
\int_0^T\intO{\varphi\Big(\vv_h\cdot\utilde{D_h}-{\mathfrak{v}}\cdot D\Big)}{\rm d}t=
\int_0^T\intO{\varphi\Big(\vv_h-[\vv_h]_h\Big)\cdot\utilde{D_h}}{\rm d}t 
$$
$$
+
\int_0^T\intO{\varphi\Big([\vv_h]_h-\Big[[\vv_h]_h\Big]_\ep\Big)\cdot\utilde{D_h}}{\rm d}t
+\int_0^T\intO{\varphi\Big[[\vv_h]_h\Big]_\ep\cdot\Big(\utilde{D_h}-D\Big)}{\rm d}t
$$
$$
+\int_0^T\intO{\varphi\Big(\Big[[\vv_h]_h\Big]_\ep-[\mathfrak{v}]_\ep\Big)\cdot D}{\rm d}t
+ \int_0^T\intO{\varphi\Big([\mathfrak{v}]_\ep-\mathfrak{v}\Big)\cdot D}{\rm d}t=\sum_{i=1}^5I_h^{i},
$$
where $\ep>0$,
and treat each integral separately. In the above and hereafter, $[f]_\ep=j_\ep*f$, $\ep>0$ "sufficiently small",
where $j_\ep$ is the standard mollifying kernel over the space variables. 

\item According to (\ref{etildevr}),
\bFormula{B1}
\utilde{D_h},\,D_h,\,D\;\mbox{is bounded in}\;L^\infty(I;L^q(\Omega))\cap L^2(0,T;L^p(\Omega)),
\eF
where $1<q<\frac{2\gamma}{\gamma+3}<2$, $6/5<p<\frac{6\gamma}{\gamma+12}<6$. 
According to (\ref{regV}) and interpolation (expressing $L^{p'}(\Omega)$ via $L^{2}(\Omega)$ and
$L^{6}(\Omega)$ if $p'\ge 2$)
\bFormula{B2}
\|[\vv_h]_h -\vv_h\|_{L^2(0,T;L^{p'}(\Omega))}\aleq h^\alpha,\;\alpha=1\,\mbox{if $p\ge 2$},\;
\alpha=\frac 52-\frac 3p\,\mbox{if $p<2$}.
\eF
Due to (\ref{B1}--\ref{B2})
\bFormula{B3}
|I^1_h|\to 0.
\eF
\item By H\"older inequality, (\ref{regV}), interpolation and Rellich-Kondrachev theorem,
\bFormula{B4}
|I^2_h|\aleq\int_0^T\Big\|[\vv_h]_h-\Big[[\vv_h]_h\Big]_\ep\Big\|_{L^{p'}(\Omega)}
\|\utilde{D_h}\|_{L^p(\Omega)}{\rm d}t
\eF
$$
\aleq
\ep^\alpha\int_0^T\Big\|\nabla_h\vv_h\Big\|^\alpha_{L^{2}(\Omega)}\|\vv_h\|^{1-\alpha}_{L^6(\Omega)}
\|\utilde{D_h}\|_{L^p(\Omega)}{\rm d}t\aleq \ep^\alpha,
$$
where $\alpha$ is the same as in (\ref{B2}).
\item Employing (\ref{co4}), (\ref{co6+}) in combination with Lemma \ref{LcalA+}, and Lemma
\ref{rieszcom}, we obtain
$$
\forall t\in [0,T],\;
\utilde{D_h}(t)\rightharpoonup D(t)\;\mbox{in $L^q(\Omega)$, where $q$ is given in (\ref{B1})}.
$$
Since
the imbedding $L^q(\Omega)\hookrightarrow [W^{1,3}(\Omega)]^*$ is compact, we deduce from the
above that the convergence is strong in $[W^{1,3}(\Omega)]^*$ for all $t\in[0,T]$,
and consequently, in particular, also strong in $L^2(0,T;[W^{1,3}(\Omega)]^*)$, i.e.
\bFormula{B5-}
\utilde{D_h}\to D\;\mbox{in $L^2(0,T; [W^{1,3}(\Omega)]^*)$},
\eF
where the convergence "in time" is handled by the Lebesgue dominated convergence theorem.

On the other hand, since $\Grad^k(j_\ep*[\vv_h]_h)=\Grad^kj_\ep*[\vv_h]_h$, $k\in\tn N$,
we have, in  particular,
\bFormula{B5+}
\Big\|\Big[[\vv_h]_h\Big]_\ep\Big\|_{W^{k,2}}(\Omega)\aleq c(\ep)\|[\vv_h]_h\|_{L^2(\Omega)}
\eF
(with $c(\ep)$ eventually singular as $\ep\to 0+$); whence, by Sobolev imbedding and (\ref{regV}),
$$
\Big\|\Big[[\vv_h]_h\Big]_\ep\Big\|_{L^2(0,T;W^{1,3}(\Omega))}\aleq c(\ep).
$$
In  view of (\ref{B5-}--\ref{B5+}),
\bFormula{B5}
|I^3_h|\to 0
\eF

\item According to (\ref{regV})
$$
[\vv_h]_h\rightharpoonup \mathfrak{w}\;\mbox{in $L^2(0,T; W^{1,2}(\Omega))$ whence in
$L^2(0,T;L^6(\Omega))$},\;\mbox{where $\mathfrak{w}=
\mathfrak{v}$}.
$$
Consequently,
\bFormula{B6}
|I^4_h|=\Big|\int_0^T\intO{\Big([\vv_h]_h-\mathfrak{v}\Big)[ D]_\ep}{\rm d}t\Big|\to 0\;\mbox{as $h\to 0$}.
\eF
\item Finally, by standard properties of mollifiers, $[\mathfrak{v}]_\ep\to\mathfrak{v}$
in $L^2(0,T;L^6(\Omega))$; whence
\bFormula{B7}
|I^5_h|\to 0.
\eF
\end{enumerate}

Resuming the above calculation, we infer,
\bFormula{C*}
{\mathfrak{C}}[\vu_h,\vr_h, \vr_h\widehat\vv_h]\to\underline{\mathfrak{C}}[\mathfrak{u},\mathfrak{r},\mathfrak{r} \mathfrak{v}].
\eF
%We have so far proved that
%\bFormula{C**}
%{\mathfrak{C}}[\widehat\vu_h,\utilde{\vr_h}, \utilde{\vr_h\widehat\vu_h}]\to\underline{\mathfrak{C}}[\mathfrak{u},\mathfrak{r}, \mathfrak{r}\mathfrak{u}]
%\eF

%Now, we set
%$$
%\mathfrak{A}_h^1=\vr_h\Grad\Delta^{-1}\Grad\cdot[\utilde{\vr_h\widehat\vu_h}-{\vr_h\widehat\vu_h}],\quad
%\mathfrak{A}_h^2=(\vr_h-\utilde{\vr_h})\Grad\Delta^{-1}\Grad\cdot[{\vr_h\widehat\vu_h}],
%$$
%$$
%\mathfrak{A}_h^3=(\utilde{\vr_h\widehat\vu_h}-{\vr_h\widehat\vu_h})\cdot\Grad\Delta^{-1}\Grad[\vr_h],\quad
%\mathfrak{A}_h^4=\vr_h\widehat\vu_h\cdot\Grad\Delta^{-1}\Grad\cdot[\utilde{\vr_h}-\vr_h].
%$$
%We have
%$$
%\int_0^T\intO{|[\widehat\vu_h]_\ep\cdot\mathfrak{A}_h^1|}{\rm d}t
%\aleq \|[\widehat\vu_h]_\ep\|_{L^2(0,T;L^{q'}(\Omega)}\|\vr_h\|_{L^\infty(0,T;L^\gamma(\Omega))}
%\|\utilde{\vr_h\widehat\vu_h}-{\vr_h\widehat\vu_h}\|_{L^2(0,T;L^{3/2}(\Omega)}
%$$
%$$
%\aleq \overline c(\ep)(\Delta t)^{\frac 16}h^{-\frac\beta 2}\;\mbox{with $c(\ep)$ tending to $\infty$ as %$\ep\to 0$},
%$$
%where we have use (\ref{interpol}), (\ref{*}) and (\ref{timeG}). The same estimates are valid for
%$\int_0^T\intO{|[\widehat\vu_h]_\ep\cdot\mathfrak{A}_h^i|}{\rm d}t$, $i=2,3,4$.
%
%
%
{\it Step 2: Convergence of $\mathfrak{ K}[\vr_h,\vu_h]$. } 
In the expression $\mathfrak{ K}$, we have grouped the terms which can be handled
relatively easily by the standard compactness arguments. Indeed, 
due to (\ref{co4}) and the second formula in Lemma \ref{LcalA}
we infer
$$
\Grad\Delta^{-1}[\vr_h]\to \Grad\Delta^{-1}[\mathfrak{r}]\;\mbox{in $L^p(0,T; C(\overline\Omega))$, $1\le p<\infty$}. 
$$
Consequently, recalling (\ref{co1+}), (\ref{co2}), (\ref{co4}), (\ref{co5}), (\ref{co6+}), (\ref{cruu}) we conclude that
\bFormula{K*}
\mathfrak{K}(\vr_h,\vu_h)\to\underline{\mathfrak{K}}(\mathfrak{r},\mathfrak{u}).
\eF
{\it Step 3: Convergence of $\mathfrak{ O}[\vr_h,\vu_h]$.}
\\ \\
We want to prove that
\bFormula{O*}
\mathfrak{ O}[\vr_h,\vu_h]\to 0.
\eF

Indeed, for the first term, we write
$$
J_h[\vu_h,\phi]=-J_h\Big[\vu_h,\varphi\nabla\Delta^{-1}\Big[[\vr_h]_h\Big]\Big] +
J_h\Big[\vu_h,\varphi\nabla\Delta^{-1}\Big[[\vr_h]_h-\vr_h\Big]\Big].
$$
In the above, the absolute value of the first term is bounded from above by
$$
\aleq h\,\|\eta\|_{C[0,T]}\|\psi\|_{C^2(\overline\Omega)}\|\Gradh\vu_h\|_{L^2(\Omega)}
\Big(\|[\vr_h]_h\|_{L^2(\Omega)} + \|\Grad[\vr_h]_h\|_{L^2(\Omega)}\Big)
\aleq h \|\Gradh\vu_h\|_{L^2(\Omega)}\times
$$
$$
\Big(\|\vr_h\|_{L^2(\Omega)} + \|\vr_h\|_{Q^{1,2}(\Omega)}\Big)\aleq h^{\frac {1-\omega-\beta} 2} A_h,
\; A_h= \|\Gradh\vu_h\|_{L^2(\Omega)} (\|\vr_h\|_{L^2(\Omega)}+h^{\frac{1+\omega+\beta} 2}\|\vr_h\|_{Q^{1,2}(\Omega)})
$$
due to (\ref{A2}) and Lemmas \ref{LcalA}--\ref{LcalA+} and (\ref{regQ}),
where $A_h$ is bounded in $L^1(0,T)$ by virtue of (\ref{e0}), (\ref{e3}), (\ref{e6}), while,
from the same reasons,
the absolute value of the second term is bounded from above by
$$
 \aleq h\,\|\eta\|_{C^[0,T]}\|\psi\|_{C^1(\overline\Omega)}\|\Gradh\vu_h\|_{L^2(\Omega)}
\|[\vr_h]_h-\vr_h\|_{L^2(\Omega)}\aleq h^{\frac {1-\omega-\beta} 2} A_h,
$$
with 
$A_h=\|\Gradh\vu_h\|_{L^2(\Omega)} h^{\frac{1+\omega+\beta} 2}\|\vr_h\|_{Q^{1,2}(\Omega)}$
bounded in $L^1(0,T)$, where we have used (\ref{A1}) instead of (\ref{A2}). Thus, the first term in $\mathfrak{O}[\vr_h,\vu_h)$ converges to $0$ provided $\omega\in (0,1-\beta)$.

In view of Lemmas \ref{LcalA}--\ref{LcalA+}, the second and third  terms converge to zero due  (\ref{NM1}) and (\ref{NC1}), respectively. Likewise, the fourth term, which is the easiest one. Finally, the convergence to zero of the last term needs to employ
(\ref{interpol}--\ref{interpol+}) and again Lemmas \ref{LcalA}--\ref{LcalA+}. 
\\ \\
{\it Step 4: Effective visous flux identity.} We deduce from (\ref{rg+1}) that
$$
p(\vr_h)\vr_h\rightharpoonup_*\overline{p(\mathfrak{r})\mathfrak{r}}\;\mbox{in 
${\cal M}([0,T]\times\overline\Omega)$}
$$
and
$$
\vr_h\Div\vu_h\rightharpoonup\overline{\mathfrak{r}\Div\mathfrak{u}}\;\mbox{in $L^2(0,T;L^{\frac{2\gamma}{\gamma+2}}(\Omega))$}.
$$
We now calculate,
$$
\lim_{h\to 0} (\ref{Z30}) - (\ref{Z31}).
$$
Taking into account (\ref{C*}), (\ref{K*}) and (\ref{O*}), we obtain that
$$
\overline{p(\mathfrak{r})\mathfrak{r}}\in L^1((0,T)\times\Omega),
$$
and identity
$$
\int_0^T\intO{\psi(t)\eta(x)\Big(\overline{p(\mathfrak{r})\mathfrak{r}} -(2\mu+\lambda)\overline{\mathfrak{r}\Div\mathfrak{u}}\Big)}{\rm d}t=
\int_0^T\intO{\psi(t)\eta(x)\Big({p(\mathfrak{r})\mathfrak{r}} -(2\mu+\lambda){\mathfrak{r}\Div\mathfrak{u}}\Big)}{\rm d}t
$$
for all $\eta\in C^1_c(\Omega)$, $\psi\in C^1_c((0,T))$, which means that
\bFormula{evf}
\overline{p(\mathfrak{r})\mathfrak{r}} -(2\mu+\lambda)\overline{\mathfrak{r}\Div\mathfrak{u}}
={p(\mathfrak{r})\mathfrak{r}} -(2\mu+\lambda){\mathfrak{r}\Div\mathfrak{u}}\;
\mbox{a.e.in $(0,T)\times\Omega$.}
\eF
Since Serre \cite{Se} and Lions \cite{Li}, this equation is called the effective viscous flux identity,
and it expresses the fact that the sum $p(\mathfrak{r})-(2\mu+\lambda)\Div\mathfrak{u}$ has
"better" regularity/summability than its components. 
%\footnote{In the stationary solutions, this quantity is even compact, see Lions \cite{Li} or \cite{No}, \cite{NoPa}.}
 %and it is one of keys in the proof of the strong convergence of the density sequence.

\subsection{Strong convergence of density}\label{74}

\subsubsection{Renormalized continuity equation on the continuous level} \label{741}

The necessary pieces of the DiPerna-Lions transport theory \cite{DL} needed for the existence theory of compressible Navier-Stokes equations
have been generalized to the continuity equation with the  non homogenous boundary data in \cite[Lemma 3.1]{ChJiNo}. The assumptions in \cite{ChJiNo} are more general
as we need here as far as the regularity and shape of the inflow-outflow surfaces is concerned (our inflow and outflow surfaces are flat, while in
\cite{ChJiNo} any union of disjoint  $C^2$-prametrized surfaces is allowed) but they are less general as far as the regularity of their boundaries is concerned
(their boundaries are required to be $C^2$ in \cite{ChJiNo} which condition is not satisfied in our case). We shall
therefore reformulate \cite[Lemma 3.1]{ChJiNo} to our situation and provide an alternative proof.

\begin{Lemma} \label{LP2}
Suppose that $\Omega \subset R^3$ is a polygonal domain with the mesh described in Section \ref{mesh} which fits to the inflow-outflow boundaries, cf. (\ref{fit}),
and let $(0<\mathfrak{r}_B,
\mathfrak{u}_B)\in C(\overline\Omega)\times C^1_c(R^3;R^3)$, $\mathfrak{u}_B\cdot\vc n\in C(\partial\Omega)$, cf. (\ref{ruB}).

Suppose further that couple $(\mathfrak{r},\mathfrak{u})$, $0\le\mathfrak{r}\in C_{\rm weak}([0,T]; L^\gamma(\Omega))\cap L^2((0,T)\times\Omega))$,  
$\mathfrak{u}-\mathfrak{u}_B\in L^2(0,T;$ $W_0^{1,2}$ $(\Omega;R^3))$, $\gamma>1$, satisfies continuity equation in the weak sense (\ref{D2}) with
any test function { $\varphi\in  C^1_c([0,T]$ $\times(\Omega\cup \Gamma^{\rm in}$ 
$\cup {\rm int}_2\Gamma^0))$.} Then we have: 
\begin{enumerate}
\item
The quantity $B(\mathfrak{r})\in C_{\rm weak}([0,T]; L^\gamma(\Omega))\cap C([0,T]; L^p(\Omega))$, $1\le p<\gamma$, and $(\mathfrak{r}, \mathfrak{u})$ is also a renormalized solution of the continuity equation (\ref{D2}), meaning that it verifies equation
\begin{equation}\label{P3}
\intO{ (B(\mathfrak{r}) \varphi) (\tau)} -\intO{ B(\mathfrak{r}_0) \varphi(0)}=
\end{equation}
$$
\int_0^\tau\intO{ \Big({B(\mathfrak{r})\partial_t\varphi}+B(\mathfrak{r}) \mathfrak{u} \cdot \Grad \varphi -\varphi\left( B'(\mathfrak{r}) \mathfrak{r} - B(\mathfrak{r}) \right) \Div \mathfrak{u} \Big)}  {\rm d} t  
- \int_0^\tau\int_{\Gamma^{\rm in}} B(\mathfrak{r}_B) 
\mathfrak{u}_B \cdot \vc{n} \varphi \ {\rm d}S_x {\rm d} t 
$$
with any $\tau\in [0,T]$ 
for any {$\varphi \in  C_c^1([0,T]\times $ 
$({\Omega}\cup\Gamma_{\rm in}$ $\cup {\rm int}_2\Gamma^0))$,} and any $B\in C^1([0,\infty))$, $B'\in L^\infty(0,\infty)$. 
\item There holds
\bFormula{P3*}
\intO{ B(\mathfrak{r}) (\tau)} -\intO{ B(\mathfrak{r}_0) \varphi(0)} +\int_0^\tau\intO{ \left( B'(\mathfrak{r}) \mathfrak{r} - B(\mathfrak{r}) \right) \Div \mathfrak{u} }  {\rm d} t
\eF  
$$
+ \int_0^\tau\int_{\Gamma^{\rm in}} B(\mathfrak{r}_B) 
\mathfrak{u}_B \cdot \vc{n} \varphi \ {\rm d}S_x {\rm d} t \le 0
$$
with any $\tau\in [0,T]$ 
 and any non negative $B\in C^1([0,\infty))$, $B'\in L^\infty(0,\infty)$. 
\end{enumerate}
\end{Lemma}
\begin{Remark}\label{rinflow}
\begin{enumerate}
\item The conditions on the renormalizing function $B$ in Lemma \ref{LP2}
can be relaxed: one can take, e.g., $B\in C([0,\infty))\cap C^1((0,\infty))$, $sB'-B\in C[0,\infty)$ and 
$|B(s)|+ |sB'(s)-B(s)|\aleq c(1+s)^{p}$, $0\le p\le\gamma/2$. 
\end{enumerate}
\end{Remark}
\noindent
{\bf Proof of Lemma \ref{LP2}}\\ \\
We proceed in several steps: First step is a consequence of classical theory and says that
$\mathfrak{r}$ is continuous in time with values in $L^1(\Omega)$ which is important to pass to
the limits in time integrated forms. Steps 2--7 are technical: They serve to produce a convenient extensions of the density and velocity
beyond the inflow and slip boundaries in such a way that the new fields verify continuity equation also in the outer vicinity of these boundaries.
In Step 8, we apply to the extended continuity equation the DiPerna-Lions  regularization technique \cite{DL} in order to obtain the renormalized 
continuity equation (\ref{P3}). The test function $\varphi=1$  is not an admissible test function in the identity (\ref{P3}), however with this test function
(\ref{P3}) turns, under certain circumstances, to an inequality (\ref{P3*}). This is proved in the last Step 9.
\\ \\
{\it Step 1: Immediate consequence of DiPerna-Lions theory}
Since (\ref{D2}) is satisfied, in particular, in the sense of distributions
on $(0,T)\times\Omega$, seeing the regularity of $(\mathfrak{r},\mathfrak{u})$
we deduce from the classical theory, in particular, that $\mathfrak{r}\in C([0,T]; L^p(\Omega))\cap C_{\rm weak}([0,T];L^\gamma(\Omega))$, $1\le p<\gamma$; 
whence the same is true for $B(\mathfrak{r})$, see \cite[Theorems 3, 5]{AN-MP}.
\\ \\
{\it Step 2: The particular geometry of inflow and outflow boundaries:}
We deduce from (\ref{fit}) that
\bFormula{Gin0}
\Gamma^{\rm in}=\cup_{{\mathfrak{k}}_{\rm in}=0}^{\overline{\mathfrak{k}}_{\rm in}}\Gamma^{\rm in}_{{\mathfrak{ k}_{\rm in}}},\; 
{\rm int}_2(\Gamma^{0})=\cup_{{\mathfrak{k}}_{0}=0}^{\overline{\mathfrak{k}}_{0}}\Gamma^0_{{\mathfrak{ k}_{0}}},\; \overline{\mathfrak{k}}_{\rm in},
\overline{\mathfrak{k}}_{0}\in\mathbb{N},
\eF
where $\Gamma^{\rm in}_{{\mathfrak{ k}_{\rm in}}}$, $\Gamma^0_{{\mathfrak{ k}_{0}}}$ are mutually disjoint plane domains such that:
\begin{enumerate}
\item For $\mathfrak{ k}_{\rm in}\neq \mathfrak{ j}_{\rm in}$, either 
$\overline{\Gamma^{\rm in}_{{\mathfrak{ k}_{\rm in}}}}\cap \overline{\Gamma^{\rm in}_{{\mathfrak{ j}_{\rm in}}}} =\emptyset $ or
$\overline{\Gamma^{\rm in}_{{\mathfrak{ k}_{\rm in}}}}\cap \overline{\Gamma^{\rm in}_{{\mathfrak{ j}_{\rm in}}}} \subset\Gamma^0.$ 
\item For $\mathfrak{ k}_{0}\neq \mathfrak{ j}_{0}$, either
$\overline{\Gamma_{{\mathfrak{ k}_{0}}}}\cap \overline{\Gamma_{{\mathfrak{ j}_{0}}}} =\emptyset $ or
${\rm int}_2(\overline{\Gamma_{{\mathfrak{ k}_{0}}}}\cup \overline{\Gamma_{{\mathfrak{ j}_{0}}}})$ is not a plane domain. 
\item For sets
\bFormula{gbd}
{\mathfrak{g}}^{\rm bd}=(\overline{\Gamma^{\rm in}}\cap\overline{\Gamma^{\rm out}})\cup 
(\overline{\Gamma^{\rm in}}\cap\overline{\Gamma^{0}})\cup (\overline{\Gamma^{\rm out}}\cap\overline{\Gamma^{0}})
\eF
there holds
$$ 
{\mathfrak{g}}^{\rm bd}=\Big(\cup_{k_{\rm bd}=1}^{\overline{k_{\rm bd}}}
\overline{\mathfrak{g}_{k_{\rm bd}}}\Big)\cup 
\Big(\cup_{j_{\rm bd}=1}^{\overline{j_{\rm bd}}} \mathfrak{p}_{j_{\rm bd}}\Big), \; \overline{k_{\rm bd}}, \overline{j_{\rm bd}}\in\mathbb{N},
%\quad {\mathfrak{g}}^{\rm out}=\Big(\cup_{k_{\rm out}=1}^{\overline{k_{\rm out}}}
%\overline{\mathfrak{g}_{k_{\rm out}}}\Big)\cup 
%\Big(\cup_{j_{\rm out}=1}^{\overline{j_{\rm out}}} \mathfrak{p}_{j_{\rm out}}\Big),
$$
where  $\mathfrak{g}_{k_{\rm bd}}$ 
%$\mathfrak{g}_{k_{\rm out}}$  
are open bounded segments in $R^3$, while
$\mathfrak{p}_{j_{\rm bd}}$ 
%$\mathfrak{p}_{j_{\rm  out}}$  
are isolated points. 
\item
Consequently, for any $\ep>0$,
\bFormula{ballA}
|B({\mathfrak{g}}^{\rm bd};\ep)|\aleq \ep^2,\;\mbox{where}\; B(A;\ep)=\{x\in R^3\,|\,d_A(x)={\rm dist}(x,A)<\ep\}.
%\;\mbox{and $A$ stands for ${\mathfrak{g}}^{\rm in}$ or
%  ${\mathfrak{g}}^{\rm out}$ } 
\eF
\end{enumerate}
{\it Step 3: Construction of outer neighborhoods of $\Gamma^{\rm in}$ and ${\rm int}_2\Gamma^0$.} Let $\Gamma$ be  any component in the decomposition (\ref{Gin0}). We denote, for $\ep>0$, $\xi\in\partial\Omega$,
\bFormula{ball+}
T^+(\Gamma;\ep):=\{x_B+s\vc n(\xi)\,|\,x_B\in \Gamma,\;s\in (0,\ep)\},\; B^+(\xi,\ep)= B(\xi,\ep)\cap R^3\setminus\overline\Omega.
\eF 
Clearly,
$$
\forall\xi\in\Gamma,\;\exists\epsilon=\epsilon(\xi)\in (0,\ep),\; B^+(\xi,\epsilon)=T^+(\Gamma;\ep)\cap B(\xi,\epsilon) \quad
\mbox{and}\quad |T^+(\Gamma;\ep)|\aleq \ep.
$$
Further, for any $x\in T^+(\Gamma;\ep)$ there exists a unique $x_B\in\Gamma$ such that $x=x_B +d_\Gamma(x)\vc n(x_B)$. We define
\bFormula{projT}
P: T^+(\Gamma;\ep)\to \Gamma,\;P(x)=x_B.
\eF
In our situation (which is a very particular case of \cite[Theorem 1,2]{Foote} calculable "by hand"),
\bFormula{pd}
P\in C^\infty (T^+(\Gamma;\ep)\cup\Gamma),\; d_\Gamma\in C^\infty (T^+(\Gamma;\ep)\cup\Gamma),\; \Grad d_\Gamma(P(x))=\vc n(P(x)).
\eF

We define the open set
\bFormula{calU+}
{\cal U}={\cal U}_\ep(\Gamma):= \cup_{\xi\in \Gamma} B(\xi;\epsilon(\xi)),\quad
{\cal U}^+={\cal U}_\ep^+(\Gamma):=\cup_{\xi\in \Gamma} B^+(\xi;\epsilon(\xi))\subset T^+(\Gamma;\ep)\cap (R^3\setminus\overline\Omega).
\eF
and realize that this can be done in such a way that
\bFormula{UtildeU}
{\cal U}^+(\Gamma)\cap{\cal U}^+(\tilde\Gamma)=\emptyset,\; \mbox{for any couple $\Gamma\neq\tilde\Gamma$ in the decomposition (\ref{Gin0})}.
\eF

The goal now is to extend the density and velocity fields $(\mathfrak{r},\mathfrak{u})$ from $\Omega$ to
$\overline{{\cal U}^+(\Gamma)}$ for any $\Gamma$ in the decomposition (\ref{Gin0}) in such a way that the extended fields will satisfy the continuity equation 
in the sense of distributions on  ${\rm int}(\Omega\cup \overline{{\cal U}^+(\Gamma)})$. The construction will depend on the fact whether $\Gamma\subset\Gamma^{\rm in}$
or $\Gamma\subset{\rm int}_2\Gamma^0$. 
\\ \\
{\it Step 4:  Properties of the flux of the boundary velocity field.}
\begin{enumerate}
\item We denote by $\mathfrak{X}$ the flux of the vector field $-\mathfrak{u}_B$, i.e. solution of the following family of Cauchy problems  for ODE,
\bFormula{EDO}
\frac{\rm d}{{\rm d}s}\mathfrak{X}(s;x)=-\mathfrak{u}_B(\mathfrak{X}),\;\mathfrak{X}(0;x)=x,\; s\in R,\;x\in R^3.
\eF
It is well known, cf. e.g. \cite[Chapter XI]{Demail}, that,
$$
\mathfrak{X}\in C^1(R\times R^3),\;
\mathfrak{X}(t,\cdot)\;\mbox{is $C^1$ diffeomorphism of $R^3$ onto $R^3$},
$$
in particular,
$$
\forall (s,x)\in R\times R^3,\;\mathfrak{X}(-s,\mathfrak{X}(s;x))=x,\;{\rm det}\Big[\Grad\mathfrak{X}(s
;x)\Big]=
{\rm exp}\Big(-\int_0^s\Div\mathfrak{u}_B(z,\mathfrak{X}(z,x)){\rm d}z\Big)>0.
$$
\end{enumerate}
{\it Step 5: Extension of the density field outside $\Gamma^{\rm in}$}.
\begin{enumerate}
\item 
Let $K\subset \Gamma$ be a compact set (with respect to the trace topology of $R^3$ on $\partial\Omega$) where
$\Gamma\subset\Gamma^{\rm in}$ is any component in decomposition (\ref{Gin0}).
We want to prove that 
\bFormula{IN}
\forall\ep>0,\;\exists\delta_K>0,\;\forall (s,\xi)\in (0,\delta)\times K,
\;
\mathfrak{X}(s,\xi)\subset {\cal U}_\ep^+(\Gamma).
\eF

Indeed:
\begin{enumerate}
\item By the uniform continuity of $\mathfrak{X}$ on compacts of $R^3\times R$ we easily get
$$
\exists \delta>0,\;\mathfrak{X}((0,\delta);K)\subset {\cal U}_\ep(\Gamma).
$$
\item Moreover,
 due to (\ref{EDO}), 
$$
\forall x_B\in \Gamma^{\rm in},\;\exists\delta>0,\;\forall s\in (0,\delta),\;
(\mathfrak{X}(s,x_B)-P(\mathfrak{X}(0,x_B)))\cdot\vc n(x_B)>0.
$$  
By the uniform  continuity of $\vc n$ on 
compacts of $\Gamma$, $P$ on compacts of $T^+(\Gamma)\cup\Gamma$ and
$\mathfrak{X}$ on compacts of $R\times R^3$, we deduce, in particular, that 
$$
\exists \delta>0,\;\forall (s,\xi)\in (0,\delta)\times K,\; \Big(\mathfrak{X}(s,\xi)-P(\mathfrak{X}(s,\xi))\Big)\cdot\vc n(P(
\mathfrak{X}(s,\xi))>0.
$$
\end{enumerate}
Consequently $\mathfrak{X}((0,\delta);K)\subset {\cal U}_\ep(\Gamma)\cap (R^3\setminus\overline\Omega)$ which 
finishes the proof of (\ref{IN})
\item{\it Construction of a diffeomorphism}. We know that $\overline\Gamma= {\cal G}(\overline{\cal O})$ where
${\cal G}:R^2\mapsto R^3$ is affine and ${\cal O}$ is a domain in $R^2$.  
Let $L_n$ be an exhaustive sequence of compacts of
${\cal O}$,
\bFormula{Ln}
L_n\subset {\rm int}_{2}L_{n+1},\;\cup_{n=1}^\infty L_n={\cal O}
\eF
so that $K_n:={\cal G}(L_n)$ is an exhaustive sequence of compacts in $\Gamma$ 
(one can take $L_{n}=\{x\in {\cal O}\,|\,{\rm dist}(x,R^2\setminus\overline{{\cal O}})\ge 1/n\}$).

We define  
\bFormula{calV+}
{\cal V}^+={\cal V}^+_\ep={\cal V}_\ep^+(\Gamma)=\cup_{n\in N}{\cal V}_n ,\;{\cal V}_n:=\mathfrak{X}((0,\delta_{K_n});{\cal G}({\rm int}_2L_n))\subset 
{\cal U}_\ep^+(\Gamma).
\eF
Finally, we define a map,
$$
\Phi: [0,\infty)\times{\cal O}\ni (s,\zeta)\mapsto \mathfrak{X}(s;{\cal G}(\zeta))\in R^3.
$$
We shall prove that $\Phi|_{[0,\delta_{K_n})\times{\rm int}_2L_n}$ is a
bijection of $[0,\delta_{{K_n}})\times {\rm int}_2 L_n$ onto
${\cal V}_n\cup {\rm int}_2L_n$ and a $C^1-$ diffeomorphism of $(0,\delta_{K_n})\times  {\rm int}_2{L_n}$ onto ${\cal V}_n$. In particular,
${\cal V}$ is open.

Since for all
$\zeta\in L_n$, $\Phi(0,\zeta)=\zeta$, in view of the theorem of local inversion, it is enough to show that
$$
\forall \zeta\in {\cal O},\; s\in R,\;
{\rm det}\Big[\partial_s\Phi,\nabla_\zeta\Phi\Big](s,\zeta)\neq 0.
$$ 
Seeing that, $\mathfrak{X}(s;\mathfrak{X}(-s;\xi))=\xi$, we infer
$$
\forall\xi\in R^3,\; s\in R, \quad
\partial_s\mathfrak{X}(s; \mathfrak{X}(-s;\xi))+\mathfrak{u}_B(\mathfrak{X}(-s;\xi))
\cdot\nabla_x\mathfrak{X}(s; \mathfrak{X}(-s;\xi))=0,
$$
i.e., equivalently,
$$
\partial_s\Phi(s,\zeta)= -\mathfrak{u}_B({\cal G}(\zeta))
\cdot\nabla_\xi\mathfrak{X}(s;{\cal G}(\zeta)),\;\mbox{in particular, for all $s>0$, $\zeta\in{\cal O}$,}
$$
we easily find that
$$
\Big[\partial_s\Phi,\partial_{\zeta_1}\Phi, \partial_{\zeta_2}\Phi\Big](s,\zeta)= -
\begin{bmatrix}
[\mathfrak{u}_B({\cal G}(\zeta))]^T\\
[\partial_{\zeta_1}{\cal G}(\zeta)]^T\\
[\partial_{\zeta_1}{\cal G}(\zeta)]^T\\
\end{bmatrix} 
\;
\begin{bmatrix}
\Grad\mathfrak{X}_1(s;{\cal G}(\zeta)), \Grad\mathfrak{X}_2(s;{\cal G}(\zeta)), \Grad\mathfrak{X}_3(s;{\cal G}(\zeta))
\end{bmatrix}, 
$$
where vectors and $\Grad$ are columns.
Whence,
$$
{\rm det}\Big[\partial_s\Phi,\nabla_\zeta\Phi\Big](s,\zeta)=
-\mathfrak{u}_B\cdot\vc n({\cal G}(\zeta))
{\rm exp}\Big(-\int_0^s\Div\mathfrak{u}_B(z,\mathfrak{X}(z;{\cal G}(\zeta))){\rm d}z\Big) >0
$$
for all $\zeta\in {\cal O}$ and $s\in R$.
\item {\it Extension of the density beyond the inflow boundary.}
We may therefore extend the boundary data $\mathfrak{r}_B\in C(\overline\Omega)$ to 
${\cal V}^+={\cal V}_\ep^+(\Gamma)$ by setting 
\bFormula{ecext}
{\mathfrak{r}}_B (\mathfrak{X}(s, {x}_B)) = \mathfrak{r}_B({x}_B){\rm exp}\Big(\int_0^s{\rm div}\mathfrak{u}_B(\mathfrak{ X}(z;{x}_B)){\rm d}z\Big). 
\eF
Clearly,  $\mathfrak{r}_B|_{\cal V^+} \in C^{1}({\cal V^+})$ while $\mathfrak{r}_B\in C(\Omega\cup\Gamma\cup {\cal V}^+)$, 
and
\begin{equation} \label{P2}
\Div (\mathfrak{r}_B \mathfrak{u}_B) = 0 \ \mbox{in}\ {\cal V}^+.
\end{equation}
\end{enumerate}
{\it Step 6: Extension of the density field beyond the slip boundary.} Let now $\Gamma\subset{\rm int}_2\Gamma^0$ be any component of the slip boundary in the decomposition (\ref{Gin0}). We set ${\cal V}^+={\cal V}_\ep^+(\Gamma)={\cal U}^+(\Gamma)$, cf. (\ref{calU+}), and we set in this case simply
\bFormula{0ext}
r_B(x)=0,\;x\in {\cal V}^+.
\eF
We take function,
\bFormula{tf}
\phi(x)=\varphi(x)\chi\Big(\frac{d_{\Gamma}(x)}{\mathfrak{e}}\Big),\;\varphi\in C^1_c(\Omega\cup\Gamma\cup{\cal V}^+),\;0<\mathfrak{e}< \frac 12{\rm dist}\Big({\rm supp\varphi},
\partial(\Omega\cup\Gamma\cup{\cal V}^+)\Big),\;
\eF
$$
\chi\in C^1[0,\infty),\;|\chi'(x)|\le 3,\;
\chi(x)
\left\{\begin{array}{c}
\in [0,1]\\
=1\;\mbox{if $x\in [0,1/2]$}\\
=0\;\mbox{if $x>1$}\\
\end{array}
\right\}
$$
and calculate
\bFormula{P2+}
\Big|\int_{\Omega}\mathfrak{r}\mathfrak{u}\cdot\Grad\phi{\rm d}x\Big|=\Big|\int_{D_\ep}\chi' \mathfrak{r}\frac{\varphi\mathfrak{u}\cdot\Grad d_\Gamma(x)}
{\mathfrak{e}}{\rm d}x\Big|=\Big|\int_{D_\ep}\chi' \mathfrak{r}\frac{\varphi\mathfrak{u}\cdot\Grad d_\Gamma(x)}
{d_{\partial\Omega}(x)}\frac{d_{\partial\Omega}(x)}{d_{\Gamma}(x)}{\rm d}x\Big|
\eF
$$
\aleq \|\mathfrak{ r}\|_{L^2(D_{\mathfrak{e}})}\Big\| \frac{\varphi\mathfrak{u}\cdot\Grad d_\Gamma}{d_{\partial\Omega}}
\Big\|_{L^2(\Omega)},\;D_{\mathfrak{e}}=
\Omega\cap\{\frac{\mathfrak{e}} 2\le d_\Gamma(x)\le 2{\mathfrak{e}}\}.
$$
By virtue of the Hardy inequality (indeed, due to (\ref{pd}), $\varphi{\mathfrak{u}\cdot\Grad d_\Gamma}\in W_0^{1,2}(\Omega)$ - cf. also Ne\v casov\'a et al.
\cite{Necasova}), 
the latter expression in (\ref{P2+}) 
verifies
$$
\int_0^T\Big|\int_{\Omega}\mathfrak{r}\mathfrak{u}\cdot\Grad\phi{\rm d}x\Big|{\rm d}t\aleq \|\mathfrak r\|_{L^2(0,T;L^2(D_{\mathfrak{e}}))}\,\|\mathfrak{u}\|_{L^2(0,T; W^{1,2}(\Omega))}
\to 0\;\mbox{as ${\mathfrak{e}}\to 0$.}
$$
\\ \\
{\it Step 7: Continuity equation extended.} Referring to the decomposition (\ref{Gin0}), we construct ${\cal V}^+(\Gamma_{{\mathfrak{k}}_{\rm in}})$ according
to (\ref{calV+}) and  ${\cal V}^+(\Gamma_{{\mathfrak{k}}_{0}})$ according to (\ref{calU+}), cf. (\ref{0ext}). These open sets are mutually disjoint by virtue of 
(\ref{UtildeU}). Finally, we set,
\bFormula{Vglobal}
\tilde{\cal V}^+:=\Big[\cup_{{\mathfrak{k}}_{\rm in}}^{\overline{{\mathfrak{k}}_{\rm in}}}
\Big({\cal V}^+(\Gamma_{{\mathfrak{k}}_{\rm in}})\cup \Gamma_{{\mathfrak{k}}_{\rm in}}\Big)\Big]\cup \Big[\cup_{{\mathfrak{k}}_{0}}^{\overline{{\mathfrak{k}}_{0}}}
\Big({\cal V}^+(\Gamma_{{\mathfrak{k}}_{0}})\cup \Gamma_{{\mathfrak{k}}_{0}}\Big)\Big],\quad \tilde\Omega=\tilde{\cal V}^+\cup\Omega
\eF
and extend $[\mathfrak{r},\mathfrak{u}]$ from $(0,T)\times\Omega$ to $(0,T)\times\tilde\Omega$ as follows
\bFormula{newru}
(\mathfrak{r},\mathfrak{u})(t,x)=
\left\{
\begin{array}{c}
(\mathfrak{ r},\mathfrak{u})(t,x)\;\mbox{if $(t,x)\in (0,T)\times\Omega$},\\
(\mathfrak{r}_B,\mathfrak{u}_B)(x)\;\mbox{if $(t,x)\in (0,T)\times\tilde{\cal V}^+$}.
\end{array}
\right\},
\eF 
where $r_B$ in $\tilde{\cal V}^+$ is defined through (\ref{ecext}) or (\ref{0ext}), according to the case.
We easily verify by using (\ref{P2}--\ref{P2+}), that the new couple $[\mathfrak{r}, \mathfrak{u}]$ satisfies continuity equation (\ref{NS1})$_1$ in the sense of distributions
on ${\cal D}((0,T)\times\tilde\Omega)$.
\\ \\
{\it Step 8: Application of the DiPerna-Lions regularization, proof of equation (\ref{P3}).}
Next, we use the regularization procedure due to DiPerna and Lions \cite{DL} applying convolution with a family of regularizing kernels
obtaining for the regularized function $[\mathfrak{r}]_{\mathfrak{e}}$,
\begin{equation} \label{P4}
\partial_t[\mathfrak{r}]_{\mathfrak{e}}+\Div ([\mathfrak{r}]_{\mathfrak{e}} \mathfrak{u} ) = R_{\mathfrak{e}} \ \mbox{a.e. in} \ (0,T)\times \tilde\Omega_{{\mathfrak{e}}},
\end{equation}
where
$$
\tilde\Omega_{{\mathfrak{e}}} = \left\{ x \in \tilde\Omega\ \Big| \ {\rm dist}(x, \partial \tilde\Omega ) > {\mathfrak{e}} \right\},
R_{\mathfrak{e}}:=\Div ([\mathfrak{r}]_{\mathfrak{e}} \mathfrak{u} )-\Div ([\mathfrak{r} \mathfrak{u}]_{\mathfrak{e}} ) \to 0 \ \mbox{in} \ L_{\rm loc}^1((0,T)\times\tilde\Omega) \ \mbox{as}\ {\mathfrak{e}} \to 0.
$$
The convergence of $R_{\mathfrak{e}}$ evoked above results from the application of the refined version of the Friedrichs lemma on commutators, see e.g. \cite{DL}
or \cite[Lemma 10.12 and Corollary 10.3]{FeNoB}.

Multiplying equation (\ref{P4}) on $B'([\mathfrak{r}]_{\mathfrak{e}})$, we get
\[
\partial_tB([\mathfrak{r}]_{\mathfrak{e}})+ \Div (B([\mathfrak{r}]_{\mathfrak{e}}) \mathfrak{u} ) + \left( B'([\mathfrak{r}]_{\mathfrak{e}}) [\mathfrak{r}]_{\mathfrak{e}} - B([\mathfrak{r}]_{\mathfrak{e}} ) \right) \Div \mathfrak{u} = B'([\mathfrak{r}]_{\mathfrak{e}}) R_{\mathfrak{e}}
\]
or equivalently,
$$
\int_{\tilde\Omega}B([\mathfrak{r}]_{\mathfrak{e}}(\tau))\varphi(\tau){\rm d}x-\int_{\tilde\Omega}B([\mathfrak{r}(0)]_{\mathfrak{e}})\varphi(0){\rm d}x=
$$
\[
\int_0^T\int_{\tilde\Omega} \Big(B([\mathfrak{r}]_{\mathfrak{e}})\partial_t\varphi+B([\mathfrak{r}]_{\mathfrak{e}}) \mathfrak{u} \cdot \Grad \varphi-\varphi \left( B'([\mathfrak{r}]_{\mathfrak{e}} ) [\mathfrak{r}]_{\mathfrak{e}} - B([\mathfrak{r}]_{\mathfrak{e}}) \right) \Div \mathfrak{u} \Big)\ \dx{\rm d}t
 - \int_0^T
\int_{\tilde\Omega}{ \varphi B'([\mathfrak{r}]_{\mathfrak{e}}) R_{\mathfrak{e}} }{\rm d} x{\ d }t
\]
for all $\tau\in [0,T]$, for any $\varphi \in C^1_c ([0,T]\times\tilde\Omega)$, $0<{\mathfrak{e}}< {\rm dist}({\rm supp}(\varphi),\partial\tilde\Omega)$. 
Thus, letting ${\mathfrak{e}} \to 0$ we get
\bFormula{dod1}
\int_{\tilde\Omega}B(\mathfrak{r}(\tau))\varphi(\tau){\rm d}x-\int_{\tilde\Omega}[B(\mathfrak{r}(0))\varphi(0)]{\rm d}x
\eF
\[
= \int_0^\tau\int_{\tilde\Omega} \Big(B(\mathfrak{r})\partial_t\varphi+B(\mathfrak{r}) \mathfrak{u} \cdot \Grad \varphi
-\varphi \left( B'(\mathfrak{r} ) \mathfrak{r} - B(\mathfrak{r}) \right) \Div \mathfrak{u} \Big)\ \dx{\rm d}t
\]
for all $\tau\in [0,T]$, for any $\varphi \in C^1_c ([0,T]\times\tilde\Omega)$. Now we write
\bFormula{dod2-}
\int_{\tilde\Omega}B(\mathfrak{r}) \mathfrak{u} \cdot \Grad \varphi{\rm d}x=
\intO{B(\mathfrak{r}) \mathfrak{u} \cdot \Grad \varphi}+\int_{\tilde{\cal V}^+} B(\mathfrak{r}) \mathfrak{u} \cdot \Grad \varphi{\rm d}x,
\eF
where, due to (\ref{P2}--\ref{0ext}), the second integral is equal to
\bFormula{dod2}
\int_{\Gamma_{\rm in}}\varphi B(\mathfrak{r}_B)\mathfrak{u}_B\cdot{\vc n}{\rm d}S_x +\int_{\tilde{\cal V}^+}\varphi (\mathfrak{r}_B B'(\mathfrak{r}_B)- B(\mathfrak{r}_B)){\rm div}\mathfrak{u}_B{\rm d}x.
\eF
%We notice that $\varphi$ is vanishing on a neighborhood of $\partial\Omega\setminus \Gamma_{\rm in}$ and that $\Gamma_{\rm in}$ 
%is Lipschitz. This justifies the latter integration by parts although ${\cal U}$ was not shown to be Lipschitz.

Now, we insert the  identities (\ref{dod2-}--\ref{dod2}) into (\ref{dod1}) and let $\ep\to 0$. 
Recall that 
$$
\tilde{\cal V}^+=\tilde{\cal V}^+_\ep\subset 
\Big[\cup_{{\mathfrak{k}}_{\rm in}}^{\overline{{\mathfrak{k}}_{\rm in}}}
\Big({ T}^+(\Gamma_{{\mathfrak{k}}_{\rm in}};\ep)\cup \Gamma_{{\mathfrak{k}}_{\rm in}}\Big)\Big]\cup \Big[\cup_{{\mathfrak{k}}_{0}}^{\overline{{\mathfrak{k}}_{0}}}
\Big({T}^+(\Gamma_{{\mathfrak{k}}_{0}};\ep)\cup\Gamma_{{\mathfrak{k}}_{0}}\Big)\Big],
$$
cf. (\ref{Vglobal}) and (\ref{ball+});
whence
$$
|\tilde{\cal V}^+_\ep|\to 0\;\mbox{as $\ep\to 0$}.
$$
Recalling regularity of $\mathfrak{r}_B$, $\mathfrak{u}_B$, cf. (\ref{P2}--\ref{0ext})) and assumption (\ref{ruB}), and summability of $(\mathfrak{r},\mathfrak{u})$, we deduce finally (\ref{P3}). This finishes proof of of the first item in Lemma \ref{LP2}.
\\ \\
{\it Step 9: Proof of inequality (\ref{P3*}).}
Let us take in equation (\ref{P2}) test function
\bFormula{L9-}
\varphi(t,x)=\varphi_{\ep}(t,x)=\eta_\ep(x)\Big( 1-\chi\Big(\frac {d_{{\mathfrak{g}}^{\rm bd}}(x)}\ep\Big)\Big),\quad
\eta_\ep(x) = \left\{ \begin{array}{l} 1 \ \mbox{if} \ d_{\Gamma^{\rm out}}(x) > \ep\\
\frac{1}{\ep} d_{\Gamma^{\rm out}}(x) \ \mbox{if}\ d_{\Gamma^{\rm out}}(x) \le \ep   \end{array} \right\},
\eF
where $\chi$ is defined in (\ref{tf}) and $\ep$ is a positive sufficiently small number. Since $\eta_\ep$ is a Lipschitz
function, { $\varphi$} is an admissible test function in (\ref{P2}). We calculate,
$$
\Grad\varphi(x)=\left\{ \begin{array}{l}
\frac 1\ep \Big( 1-\chi\Big(\frac {d_{{\mathfrak{g}}^{\rm bd}}(x)}\ep\Big)\Big)\,\Grad d_{\Gamma^{\rm out}}(x) -
\frac 1\ep\,\frac{d_{\Gamma^{\rm out}}(x)}{\ep} \chi'\Big(\frac {d_{{\mathfrak{g}}^{\rm bd}}(x)}\ep\Big)\,\Grad d_{{\mathfrak{g}}^{\rm bd}}(x)\;
\mbox{if $ d_{\Gamma^{\rm out}}(x) < \ep$}\\
-\frac 1\ep \chi'\Big(\frac {d_{{\mathfrak{g}}^{\rm bd}}(x)}\ep\Big)\,\Grad d_{{\mathfrak{g}}^{\rm bd}}(x)\;
\mbox{if $ d_{\Gamma^{\rm out}}(x) > \ep$}
\end{array} \right\}.
$$
By virtue of (\ref{pd}), we can choose $\ep>0$ sufficiently small in such a way that
\bFormula{L9-1}
\mathfrak{u}_B\cdot\Grad d_{\Gamma^{\rm out}}(x)>0\; \mbox{for all}\; x\in \Big(B(\Gamma^{\rm out} ;\ep)\setminus B({\mathfrak{g}}^{\rm bd};\ep)\Big)\cap\Omega.
\eF
and write
$$
\int_0^\tau\intO{\mathfrak{r}\mathfrak{u}\cdot\Grad\varphi}{\rm d}t=
\int_0^\tau\intO{\mathfrak{r}(\mathfrak{u}-\mathfrak{u}_B)\cdot\Grad\varphi}{\rm d}t
+ \int_0^\tau\int_{\Omega\setminus B(\Gamma^{\rm out} ;\ep)}\mathfrak{r}\mathfrak{u}_B\cdot\Grad\varphi{\rm d}x{\rm d}t
$$
\bFormula{L9-2}
+\int_0^\tau\int_{(B(\Gamma^{\rm out} ;\ep)\setminus B({\mathfrak{g}}^{\rm bd};\ep))\cap\Omega}\mathfrak{r}\mathfrak{u}_B\cdot\Grad\varphi{\rm d}x{\rm d}t
+\int_0^\tau\int_{B({\mathfrak{g}}^{\rm bd};\ep)\cap\Omega}\mathfrak{r}\mathfrak{u}_B\cdot\Grad\varphi{\rm d}x{\rm d}t
\eF
Since $\Omega$ is a bounded Lipschitz domain and $\mathfrak{u}-\mathfrak{u}_B\in L^2(0,T;W_0^{1,2}(\Omega))$, we have by the same reasonning as in (\ref{P2+}),
by using essentially the Hardy inequality, Lipschitz continuity of the distance function and estimates (\ref{ballA}), (\ref{ball+}),
$$
\Big|\int_0^\tau\intO{\mathfrak{r}(\mathfrak{u}-\mathfrak{u}_B)\cdot\Grad\varphi}{\rm d}t\Big|\to 0\;\mbox{as $\ep\to 0+$}.
$$
Likewise, but more simply, using (\ref{ballA}), (\ref{ball+}) and uniform estimates,
\bFormula{L9-3}
\Big|\int_0^\tau\int_{\Omega\setminus B(\Gamma^{\rm out} ;\ep)}\mathfrak{r}\mathfrak{u}_B\cdot\Grad\varphi{\rm d}x{\rm d}t\Big|\to 0\;\mbox{as $\ep\to 0+$}
\eF
and
\bFormula{L9-4}
\Big|\int_0^\tau\int_{B({\mathfrak{g}}^{\rm bd};\ep)\cap\Omega}\mathfrak{r}\mathfrak{u}_B\cdot\Grad\varphi{\rm d}x{\rm d}t\Big| \to 0\;\mbox{as $\ep\to 0+$}.
\eF
Finally, due to the choice (\ref{L9-1}),
\bFormula{L9-5}
\int_0^\tau\int_{B(\Gamma^{\rm out} ;\ep)\setminus B({\mathfrak{g}}^{\rm bd};\ep))\cap\Omega}\mathfrak{r}\mathfrak{u}_B\cdot\Grad\varphi{\rm d}x{\rm d}t\le 0.
\eF
Whence letting in (\ref{P3}), $\ep\to 0$ while employing the Lebesgue dominated convergence theorem and (\ref{ballA}), (\ref{ball+}) and uniform estimates,
we get inequality (\ref{P3*}). The proof of Lemma \ref{LP2} is thus complete.

\subsubsection{Discrete renormalized continuity equation}\label{742}
Using in the discrete continuity equation test function $\phi\approx \phi B'(\vr^k)$,
$\phi\in Q(\Omega)$ we obtain after an elementary but laborious calculation 
the discrete renormalized continuity equation, see \cite[Lemma 5.1]{KwNoNUM}. 
\bFormula{r1}
\intO{ D_t B(\vr^k) \phi } 
+  \sum_{K\in {\cal T}}\sum_{\sigma \in {\cal E}(K)\cap{\cal E}_{\rm int}}\int_\sigma F_{\sigma,K}(B(\vr^k),\vu^k)\phi{\rm d}S_x
\eF
$$
{ +  h^\omega\sum_{\sigma\in {\cal E}_{\rm int}}\int_\sigma[[\vr^k]]_\sigma[[\phi B'(\vr^k)]]_{\sigma}{\rm d}S_x}
-\intO{\phi\Big(B(\vr^k)-\vr^k B'(\vr^k)\Big)\Divh\vu^k}
$$
$$
{ + \frac 1{\Delta t}\intO{E_B(\vr^{k-1}|\vr^{k})\phi}}+
\sum_{\sigma\in {\cal E}_{\rm int}}\int_\sigma E_B(\vr_\sigma^{k,+}|\vr^{k,-}_\sigma)
|\vu_{B,\sigma}\cdot\vc n|\phi^+{\rm d S}_x
$$
$$
+
\sum_{\sigma\in{\cal E}^{\rm out}}\int_\sigma B(\vr^{k})\vu_{B,\sigma}\cdot\vc n_\sigma\phi{\rm d}S_x
+\sum_{\sigma\in {\cal E}^{\rm in}}\int_{\sigma}
B(\vr_{B})\vu_{B,\sigma}\cdot\vc n_{\sigma}\phi{\rm d}S_x
$$
$$
+ \sum_{\sigma\in {\cal E}^{\rm in}}\int_{\sigma}
E_B(\vr_{B}|\vr^k)|\vu_{B,\sigma}\cdot\vc n_{\sigma}|\phi{\rm d}S_x
=0\;\mbox{$k=1,\ldots,N$},
$$ 
where $B\in C([0,\infty))\cap C^1(0,\infty)$, $sB'\in C[0,\infty)$.
Using (\ref{pwt+}) and formula (\ref{Up3}), we deduce from (\ref{r1})
\bFormula{r1+}
\int_0^T\intO{ \partial_tB(\utilde{\vr}) \phi }
-\int_{0}^T\intO{\Big(B(\vr)\vu\cdot\Grad\phi 
+\phi\Big(B(\vr)-\vr B'(\vr)\Big)\Divh\vu\Big)}{\rm d}t
\eF
$$
+
\int_{0}^T\int_{\Gamma^{\rm out}} B(\vr)\mathfrak{u}_{B}\cdot\vc n\cdot\phi{\rm d}S_x{\rm d}t
+\int_{0}^T\int_{\Gamma^{\rm in}}
B(\vr_{B})\mathfrak{u}_{B}\cdot\vc n\phi{\rm d}S_x{\rm d}t +  h^\omega\sum_{\sigma\in {\cal E}_{\rm int}}\int_{\Delta t}^\tau\int_\sigma
[[\vr]]_\sigma[[B'(\vr)]]_{\sigma}\widehat\phi^+{\rm d}S_x{\rm d}t
$$
%%%
$$
{ + \frac 1{\Delta t}\int_{0}^T\intO{E_B(\vr(\cdot-\Delta t)|\vr(\cdot))\widehat\phi}}{\rm d}t+
\sum_{\sigma\in {\cal E}_{\rm int}}\int_{0}^T\int_\sigma E_B(\vr_\sigma^{+}|\vr^{-}_\sigma)
|\vu_{B,\sigma}\cdot\vc n|\widehat\phi{\rm d S}_x{\rm d}t
$$
$$
+ \sum_{\sigma\in {\cal E}^{\rm in}}\int_{0}^T\int_{\sigma}
E_B(\vr_{B}|\vr)|\vu_{B,\sigma}\cdot\vc n_{\sigma}|\widehat\phi{\rm d}S_x{\rm d}t=
\int_{0}^\tau<R_{h,\Delta t}^{B},\phi>{\rm d}t
$$
for all $\phi \in C^1([0,T]\times\overline\Omega)$, where
%%%%%%%%%%%%%%%%%%%%%%%%%%%%%
$$
<R_{h,\Delta t}^{B},\phi>=
%1_{[\min\{\Delta t,\tau\},\tau)}
\Big[
\sum_{K\in {\cal T}}\sum_{\sigma\in {\cal E}(K)\cap{\cal E}_{\rm int} }\int_\sigma ( \widehat\phi - \phi ) \ju{B(\vr)}_{{\vc n}_{\sigma,K}}  \ [\vu_\sigma \cdot \vc{n}_{\sigma,K}]^- {\rm d}S_x  
$$
$$
+\sum_{K\in {\cal T}} 
\sum_{\sigma\in {\cal E}(K)}\int_{\sigma} (\phi-\widehat\phi) B(\vr)(\vu - \vu_\sigma ) \cdot \vc{n}_{\sigma,K} \ {\rm dS}_x +\intO{B(\vr)(\phi-\widehat\phi){\rm div}_h\vu} 
$$
$$
{ { -  h^\omega\sum_{\sigma\in {\cal E}_{\rm int}}\int_\sigma[[\vr]]_\sigma B'(\vr^-)[[\widehat\phi ]]_{\sigma}{\rm d}S_x}}
+  \sum_{\sigma\in {\cal E}^{\rm in}}\int_\sigma
(B(\vr)-B(\vr_B))\vu_{B,\sigma}\cdot \vc n_{\sigma,K}(\widehat\phi-\phi){\rm d}S_x
$$
$$
{ +  \sum_{\sigma\in {\cal E}^{\rm out}}\int_\sigma
B(\vr)(\mathfrak{u}_B-\mathfrak{u}_{B,\sigma})\cdot \vc n_{\sigma}\phi{\rm d}S_x
+  \sum_{\sigma\in {\cal E}^{\rm in}}\int_\sigma
B(\vr_B)(\mathfrak{u}_B-\mathfrak{u}_{B,\sigma})\cdot \vc n_{\sigma}\phi{\rm d}S_x
}\Big].
$$
%$$
%+1_{[0,\min\{\Delta t,\tau\})}\Big[-\intO{\Big(B(\vr)\widehat\vu\cdot\Grad\phi  
%+\phi\Big(B(\vr)-\vr B'(\vr)\Big)\Divh\vu\Big)}
%$$
%$$
%+
%\int_{\Gamma^{\rm out}} B(\vr)\mathfrak{u}_{B}\cdot\vc n\cdot\phi{\rm %d}S_x+\int_{\Gamma^{\rm in}}
%B(\vr_{B})\mathfrak{u}_{B}\cdot\vc n\phi{\rm d}S_x +  h^\omega\sum_{\sigma\in {\cal E}_{\rm %int}}\int_\sigma
%[[\vr]]_\sigma[[B'(\vr)]]_{\sigma}\widehat\phi^+{\rm d}S_x
%$$
%%%
%$$
%{ + \frac 1{\Delta t}\intO{E_B(\vr(\cdot-\Delta t)|\vr(\cdot))\widehat\phi}}+
%\sum_{\sigma\in {\cal E}_{\rm int}}\int_\sigma E_B(\vr_\sigma^{+}|\vr^{-}_\sigma)
%|\vu_{B,\sigma}\cdot\vc n|\widehat\phi{\rm d S}_x
%$$
%$$
%+ \sum_{\sigma\in {\cal E}^{\rm in}}\int_{\sigma}
%E_B(\vr_{B}|\vr)|\vu_{B,\sigma}\cdot\vc n_{\sigma}|\widehat\phi{\rm d}S_x\Big].
%$$

\subsubsection{Combining continuous and discrete continuity equations}\label{743}

We now take in (\ref{r1+}) $B=L_\zeta$ with $L_\zeta$ a convex function on $[0,\infty)$,
\bFormula{L}
L_\zeta(\vr)=\vr\log(\vr+\zeta)+a\vr +{b}\;\mbox{where $a, b\ge 0$, $\zeta>0$ are such that $L_\zeta\ge 0$}.
\eF
Clearly, likewise as in Lemma \ref{Lc-consistency}, with this choice of $B$,
 \bFormula{NR1}
\Big|\int_{0}^T<R_{h,h}^{L_\zeta},\phi>{\rm d}t\Big|\aleq c(\zeta) h^{\alpha_r}
\|\phi\|_{L^\infty(0,T; W^{1,\infty}(\Omega))},\;0<\alpha_r<\alpha_c.
\eF

We introduce,
\bFormula{L9}
\phi(t,x)=\phi_{\ep,\delta}(t,x)=\varphi_{\ep}(x) \psi_{\tau,\delta}(t)\;\quad \varphi_{\ep}(x) =\eta_\ep(x)\Big(1-\chi\Big(\frac{d_{{\mathfrak{g}^{\rm bd}}}(x)}\ep\Big)\Big)
\eF
where
$$
\psi_{\tau,\delta}(t) = 
\left\{\begin{array}{l}
1\; \mbox{if $t\in (-\infty,\tau-\delta] $}\\
    1-\frac{t-\tau+\delta }{\delta}\; \mbox{if $t \in [\tau-\delta,\tau] $}\\
    0 \;\mbox{if $t \in (\tau, +\infty) $}
		\end{array}\right\},
		$$
		where $\tau\in (\Delta t,T],$ $0<\ep<\ep_0$, $0<\delta<\tau-\Delta t$, $\ep_0$ sufficiently small and $\eta_\ep$ and $\chi$ are the same as in (\ref{L9-}).

We shall use in (\ref{r1+}) $\phi(t,x)$ as
test functions. Neglecting several non-negative terms at the left hand side, we deduce from 
(\ref{r1+}) with $\Delta t=h$, 
\bFormula{*1*}
\frac 1\delta \int_{\tau-\delta}^\tau\intO{L_\zeta(\vr)\varphi_\ep(x)}{\rm d}t -\intO{L_\zeta(\vr^0)\varphi_\ep(x)}
-\int_{0}^T\intO{L_\zeta(\vr)(\vu-\vu_B)\cdot\Grad\phi}{\rm d}t
\eF
$$ 
-\int_{0}^T\intO{L_\zeta(\vr)\vu_B\cdot\Grad\phi}{\rm d}t 
+\int_{0}^T\intO{\phi \Big(\vr L'_\zeta(\vr)-L_\zeta(\vr)\Big)\Divh\vu}{\rm d}t
$$
$$
+\int_{0}^T\int_{\Gamma^{\rm in}}
L_\zeta(\vr_{B})\mathfrak{u}_{B}\cdot\vc n\phi{\rm d}S_x\aleq \int_{0}^\tau<R_{h,h}^{L_\zeta},\phi>{\rm d}t.
$$

By virtue of (\ref{e0}), (\ref{e3}),
$$
\begin{array}{c}
L_\zeta(\vr_h)\rightharpoonup_*\overline{L_\zeta(\mathfrak{r})}\;\mbox{in $L^\infty(0,T; L^p(\Omega))$, $1<p<\gamma$},
\\ \\
(\vr_hL_\zeta(\vr_h)-L_\zeta(\vr_h))\Divh\vu_h \rightharpoonup
\overline{(\mathfrak{r} L_\zeta(\mathfrak{r})-L_\zeta(\mathfrak{r}))\Div\mathfrak{u}}\;
\mbox{in $L^\infty(0,T; L^p(\Omega))$ for some $p>1$}.
\end{array}
$$
Consequently, letting $h\to 0$ in (\ref{*1*}), and using (\ref{NR1}), we get
\bFormula{*1+*}
\frac 1\delta\int_{\tau-\delta}^\tau\intO{\overline{L_\zeta(\mathfrak{r})}\varphi_\ep(x)}{\rm d}t
-\intO{L_\zeta(\mathfrak{r}_0)\varphi_\ep(x)}-
\int_{0}^T\intO{\overline{L_\zeta(\mathfrak{r})}\mathfrak{u}_B\cdot\Grad\phi}{\rm d}t 
\eF
$$
+\int_{0}^T\intO{\phi \,\overline{(\mathfrak{r} L'_\zeta(\mathfrak{r})-L_\zeta(\mathfrak{r}))\Div\mathfrak{u}}}{\rm d}t
$$
$$
+\int_{0}^T\int_{\Gamma^{\rm in}}
L_\zeta(\mathfrak{r}_{B})\mathfrak{u}_{B}\cdot\vc n\phi{\rm d}S_x
\aleq \limsup_{h\to 0+}\int_0^T\intO{L_\zeta(\vr)(\vu-\vu_B)\cdot\Grad\phi}{\rm d}t.
$$

Combining (\ref{P3}) $B=L_\zeta$ and test function $\varphi=\phi$, and (\ref{*1+*})
one gets

\bFormula{*1++*}
\frac 1\delta\int_{\tau-\delta}^\tau\intO{\overline{L_\zeta(\mathfrak{r})}\varphi_\ep(x)}{\rm d}t
-\frac 1\delta\int_{\tau-\delta}^\tau\intO{L_\zeta(\mathfrak{r})\varphi_\ep(x)}{\rm d}t
-\int_{0}^T\intO{\Big(\overline{L_\zeta(\mathfrak{r})}-L_\zeta(\mathfrak{r})\Big)\mathfrak{u}_B\cdot\Grad\phi}{\rm d}t 
\eF
$$
+\int_{0}^T\intO{\phi \,\Big(\overline{(\mathfrak{r} L'_\zeta(\mathfrak{r})-L_\zeta(\mathfrak{r}))\Div\mathfrak{u}}-(\mathfrak{r} L'_\zeta(\mathfrak{r})-L_\zeta(\mathfrak{r}))\Div\mathfrak{u}\Big)}{\rm d}t
\aleq \Pi_{\tau,\delta,\ep},
$$
where
$$
\Pi_{\tau,\delta,\ep}:=
\limsup_{h\to 0+}\int_0^T\intO{L_\zeta(\vr)(\vu-\vu_B)\cdot\Grad\phi}{\rm d}t - \int_0^T\intO{L_\zeta(\mathfrak{r})(\vu-\vu_B)\cdot\Grad\phi}{\rm d}t.
$$
%+\int_{\tau-\delta}^\tau\intO{L_\zeta(\mathfrak{r})\psi_{\tau,\delta}\mathfrak{u}\cdot\Grad\varphi_\ep}{\rm d}t
%$$
%$$
%-\int_{\tau-\delta}^\tau\intO{\phi \,(\mathfrak{r} L'_\zeta(\mathfrak{r})-L_\zeta(\mathfrak{r}))\Div\mathfrak{u}}{\rm d}t
%+ \int_{\tau-\delta}^\tau\int_{\Gamma^{\rm in}}
%L_\zeta(\mathfrak{r})\mathfrak{u}_{B}\cdot\vc n\phi{\rm d}S_x
%$$

Since $\mathfrak{u}-{\mathfrak{u}_B}\in L^2(I;W^{1,2}_0(\Omega))$, we have, by the same token as in (\ref{P2+}),
$$
|\Pi_{\tau,\delta,\ep}|\aleq\tau o(\ep),\;\mbox{where $\lim_{\ep\to 0+} o(\ep)=0$}
$$
In view of (\ref{L9-1}) and  Lemma \ref{Lemma3}, we reserve to the integral
$$
\int_{0}^T\intO{\Big(\overline{L_\zeta(\mathfrak{r})}-L_\zeta(\mathfrak{r})\Big)\mathfrak{u}_B\cdot\Grad\phi}{\rm d}t
$$ 
the
same treatment as in (\ref{L9-3}--\ref{L9-5}), in order to obtain, 
$$
-\liminf_{\ep\to 0+}\int_{0}^T\intO{\Big(\overline{L_\zeta(\mathfrak{r})}-L_\zeta(\mathfrak{r})\Big)\mathfrak{u}_B\cdot\Grad\phi}{\rm d}t\ge 0.
$$
%$$
%\Big|\intO{L_\zeta(\vr)(\vu-\vu_B)\cdot\Grad\phi}\Big|\aleq
%\int_{\{x\in \Omega|{\rm dist}(x, \Gamma_{\rm out})\le \ep\}}L_\zeta(\vr)\frac{|\mathfrak{v}|}
%{{\rm dist}(x, \Gamma_{\rm out})}{\rm d}x 
%$$
%$$
%+\intO{L_\zeta(\vr)|\vv-\mathfrak{v}|\cdot\Grad\phi}
%\aleq h\|\nabla\mathfrak{v}\|_{L^2(\Omega)}
%\|L_\zeta(\vr)\|_{L^2(\Omega)} +\|\Grad\mathfrak{v}\|_{L^2(\Omega)} |L_\zeta(\vr)\|_{L^2(\{x\in \Omega|{\rm dist}(x, \Gamma_{\rm out})\le \ep\})}
%$$
%where we have used (\ref{errorV}) and the Hardy inequality in $W^{1,2}_0(\Omega)$. We employ estimates (\ref{e0}),
%(\ref{e1}), (\ref{e3}) and (\ref{co1+}), (\ref{co2}), (\ref{cru+}) to get
%$$
%|\Pi_{\tau,\ep,\delta}|\aleq o_1(\ep) + o_2(\delta),\;\lim_{\ep\to 0+} o_1(\ep)=0,\;  \lim_{\delta\to 0+} o_2(\delta)=0.
%$$
%Moreover,
%since for all $x_B\in{\rm int}\sigma$, $\sigma\in {\cal E}_{\rm ext}$,
%$$
%\Grad{\rm dist}(x_B-\ep\vc n(x_B))\to -\vc n(x_B),
%$$
%and since $L_\zeta$ is convex
%we infer,
%$$
%-\liminf_{\ep\to 0+}\int_{0}^T\intO{\Big(\overline{L_\zeta(\mathfrak{r})}-L_\zeta(\mathfrak{r})\Big)\mathfrak{u}_B\cdot\Grad\phi}{\rm d}t\ge 0,
%$$
%where we have used Lemma \ref{Lemma3}.
Thus, inequality (\ref{*1++*}) becomes
\bFormula{*2*}
\frac 1\delta\int_{\tau-\delta}^\tau\intO{\overline{L_\zeta(\mathfrak{r})}}{\rm d}t
-\frac 1\delta\int_{\tau-\delta}^\tau\intO{L_\zeta(\mathfrak{r})(\tau-\delta)}{\rm d}t
\eF
$$
+\int_{0}^T\intO{\phi \,\Big(\overline{(\mathfrak{r} L'_\zeta(\mathfrak{r})-L_\zeta(\mathfrak{r}))\Div\mathfrak{u}}-(\mathfrak{r} L'_\zeta(\mathfrak{r})-L_\zeta(\mathfrak{r}))\Div\mathfrak{u}\Big)}{\rm d}t
\aleq 0.
$$

Now, we let $\zeta\to 0+$ in (\ref{*2*}) by using the lower weak semi-continuity of $L^1$norms. For any $\eta>0$ there is $h_\eta>0$ (decreasing to $0$ with $\eta\to 0$) such that
$$
\|\overline{L(\mathfrak{r})}-\overline{L_\zeta(\mathfrak{r})}\|_{L^1((0,T)\times\Omega)}\le
\liminf_{h\to+}\|{L(\vr_{h})}-
{L_\zeta(\vr_{h})}\|_{L^1((0,T)\times\Omega)}
\aleq \|{L(\vr_{h_\eta})}-
{L_\zeta(\vr_{h_\eta})}\|_{L^1((0,T)\times\Omega)} +\eta
$$
$$
\aleq \Big\|\vr_{h_\eta}\log\Big(\frac{\vr_{h_\eta}}{\vr_{h_\eta}+\zeta}\Big)\Big\|_{L^1((0,T)\times\Omega)}+\eta, 
\quad L(\mathfrak{r})=\mathfrak{r}\log\mathfrak{r}+
a\mathfrak{r}+b.
$$
where
$$
\Big\|\vr_{h_\eta}\log\Big(\frac{\vr_{h_\eta}}{\vr_{h_\eta}+\zeta}\Big)\Big\|_{L^1((0,T)\times\Omega)}\to 0 
$$
by the Lebesgue dominated convergence theorem. Consequently,
$$
\|\overline{L(\mathfrak{r})}-\overline{L_\zeta(\mathfrak{r})}\|_{L^1((0,T)\times\Omega)}\to 0.
$$
Similar reasoning in the other terms in (\ref{*2*}) turns inequality (\ref{*2*}) to
\bFormula{*3*}
\frac 1\delta\int_{\tau-\delta}^\tau\intO{\overline{L(\mathfrak{r})}}{\rm d}t
-\intO{L(\mathfrak{r})(\tau-\delta)}
\eF
$$
+\int_{0}^T\intO{\phi \,\Big(\overline{\mathfrak{r} \Div\mathfrak{u}}-\mathfrak{r}\Div\mathfrak{u}\Big)}{\rm d}t
\aleq 0.
$$

Finally we let $\delta\to 0$ in (\ref{*3*}). 
%while in all other terms,
%we shall employ the Lebesgue dominated convergence theorem. 
The final inequality reads,
\bFormula{*4*}
\intO{\Big(\overline{\mathfrak{r}\log\mathfrak{r}}-\mathfrak{r}\log\mathfrak{r}\Big)(\tau)}
\aleq\int_{0}^\tau\intO{ \Big( \mathfrak{r}\Div\mathfrak{u}-\overline{\mathfrak{r}\Div\mathfrak{u}}
\Big)}{\rm d}t
\eF
for almost all $\tau\in (0,T)$, where we have used the Theorem on Lebesgue points.

\subsection{Strong convergence of the density sequence. Limit in the momentum equation ends.}

The right hand side of the expression (\ref{*4*}) is less or equal than zero by virtue of the effective viscous flux identity (\ref{evf}) and Lemma \ref{Lemma4}. 
On the other hand, by Lemma \ref{Lemma3}, $\overline{\mathfrak{r}\log\mathfrak{r}}-\mathfrak{r}\log\mathfrak{r}\ge 0$ whence
the left hand side is greater or equal to zero. Consequently,
$$
\overline{\mathfrak{r}\log\mathfrak{r}}-\mathfrak{r}\log\mathfrak{r}\ge 0\;\mbox{a.e. in $(0,T)\times\Omega$.}
$$
Lemma \ref{Lemma3}, again, implies
\bFormula{vrstrong}
\vr_h\to\mathfrak{r}\;\mbox{a.e. in $(0,T)\times\Omega$.}
\eF
The latter, in conjunction with (\ref{co5}) yields $\overline{p(\mathfrak{r})}= p(\mathfrak{r})$ in (\ref{D3*}).
We have thus proved that the weak limit $(\mathfrak{r},\mathfrak{u})$ satisfies the momentum equation (\ref{D3}).

\subsection{Limit in the energy identity}

Neglecting several positive terms at the left hand side and using the definition of { $\Pi^V$ and (\ref{fit}) in order to replace $\vu_{B,\sigma}=\tilde{\mathfrak{u}}_{B,\sigma}$ by $\mathfrak{u}_B$ in the
remaining boundary integrals,} we can rewrite the energy identity (\ref{ebalance}) 
in the following form -- see \cite[Section 7.3]{KwNoNUM},
\bFormula{ee-v}
\Big[\intO{\Big(\frac 12\vr |\widehat\vv|^2 + H_h(\vr)\Big) }\Big]_0^\tau +\int_0^\tau\intO{\mathbb{S}(\nabla_h\vu):\nabla_h\vv}{\rm d}t
\eF
$$
+\int_0^\tau\int_{\Gamma^{\rm out}} H_h(\vr){\mathfrak{u}}_{B}\cdot\vc n
{\rm d}S_x {\rm d}t
 \aleq
  -\int_0^\tau\int_{\Gamma^{\rm in}} { H_h(\vr_B)}{\mathfrak{u}}_{B}\cdot\vc n
{\rm d}S_x {\rm d}t
$$
$$
- \int_0^\tau\intO{\Big(\vr\widehat\vu\otimes\widehat\vu+p_h(\vr)\mathbb{I}\Big):\nabla_h\vu_B}{\rm d}t
+\int_0^\tau\intO{\vr\widehat\vu\cdot\nabla_h\vu_B\cdot\widehat\vu_B}{\rm d}t  
 + { R}_h^E[\vr,\vu],
 $$
where
$$
\Big[{ R}_h^E[\vr,\vu]\Big](\tau)=
\int_\tau^{\tau_m}\Big[\int_{\Gamma^{\rm in}} { H_h(\vr_B)}|\vu_{B}\cdot\vc n|
{\rm d}S_x 
$$
$$
- \intO{\Big(\vr\widehat\vu\otimes\widehat\vu+p_h(\vr)\mathbb{I}\Big):\nabla_h\vu_B}
+\intO{\vr\widehat\vu\cdot\nabla_h\vu_B\cdot\widehat\vu_B}\Big]{\rm d}t,\;
%$$
%$$
%{ +\int_0^\tau\int_{\Gamma^{\rm out}} H_h(\vr)(\mathfrak{u}_B-\tilde{\mathfrak{u}}_B){\rm d}S_x{\rm d}t
%+\int_0^\tau\int_{\Gamma^{\rm in}} H_h(\vr)(\mathfrak{u}_B-\tilde{\mathfrak{u}}_B){\rm d}S_x{\rm d}t,
%}\;
\tau\in(\tau_{m-1},\tau_m],\; m=1,\ldots,N ; 
$$
whence
\bFormula{enb3}
\Big|{ R}_h^E[\vr,\vu]\Big|\aleq h,\;{ \mbox{ where $\Delta t=h$} }
%- \sum_{\sigma\in {\cal E}^{\rm in}}\int_0^\tau\int_\sigma \vr_B\vu_{B,\sigma}\cdot\vc n_\sigma |\widehat\vv|^2{\rm d}S_x{\rm d}t\aleq h
\eF
by virtue of (\ref{e0}--\ref{e0+}), { (\ref{ru}), (\ref{ruu})}.

In view of the convergence established in (\ref{rg+1}) and (\ref{vrstrong}), (\ref{co2}), (\ref{co4}), (\ref{cru}), (\ref{cruu}), and
(\ref{globalV}--\ref{errorV}) it is rudimentary to pass to the limit in the inequality (\ref{ee-v}) in order to get the energy inequality
(\ref{D4}). It is to be noticed that one must use the lower weak continuity of the convex functionals with respect to the weak convergence
when passing to the limit at the left hand side (and eventually integrate the left hand side over short time intervals to be able to use
the theorem on Lebesgue points to pass to the limit in the term $\intO{\frac 12\vr |\widehat\vv|^2}$, cf. e.g. \cite{KwNo} for more details).
Theorem \ref{TN2} is proved.

\section{Appendix: Numerical and Theoretical background}\label{Prel}
\subsection{Preliminaries from numerical analysis} \label{Preln}

 In this part, we recall  several classical inequalities related to the discrete functional spaces which will be used throughout the paper.

\subsubsection{Some useful elementary inequalities}

We recall Jensen's inequalities
\bFormula{JensenV}
\begin{array}{c}
\|\widehat v\|_{L^q(K)}\aleq \|v\|_{L^q(K)}\;\mbox{for all $v\in L^q(K)$, $1\le q\le\infty$},\\
%\|\widehat v\|_{L^q(\Omega)}\aleq \|v\|_{L^q(\Omega)}\;\mbox{for all $v\in L^q(K)$, $1\le q\le\infty$},\\
\| v_\sigma\|_{L^q(\sigma)}\aleq \|v\|_{L^q(\sigma)}\;\mbox{for all $v\in L^q(\sigma)$, $1\le q\le\infty$}.
\end{array}
\eF
together with the error estimate
\bFormula{errorQ}
\forall v\in W^{1,q}(\Omega),\; 
\left\{\begin{array}{c}
        \left\| v - \widehat v \right\|_{L^q(K)} \aleq h \| \Grad v \|_{L^q(K; R^3)}\;\forall v\in W^{1,q}(K)\\
       \left\| v - \widehat v \right\|_{L^q(\Omega)} \aleq h \| \Grad v \|_{L^q(\Omega; R^3)}\; \forall v\in W^{1,q}(\Omega)
       \end{array}
       \right\},\; 1\le q\le \infty.
\eF

We also recall the Poincar\'e type inequalities of  mesh elements,
\bFormula{PoincareV-} 
\begin{array}{c}
\|v-v_\sigma\|_{L^q(K)}\aleq h \| \Grad v \|_{L^q(K; R^{3})}, \ \forall\sigma\in {\cal E}(K),\,v\in W^{1,q}(K),\,1\le q\le\infty,\\
\|v-v_K\|_{L^q(K)}\aleq h \| \Grad v \|_{L^q(K; R^{3})}, \ \forall K\in {\cal T},\,v\in W^{1,q}(K),\,1\le q\le\infty.
\end{array}
\eF

\subsubsection{Properties of piecewise constant and  Crouzeix-Raviart finite elements}\label{Prop}

We report the estimates of jumps on mesh elements,
\bFormula{PoincareV} 
\begin{array}{c}
\|[[\widehat v]]_{\sigma=K|L}\|_{L^q(\sigma)}\aleq h \| \nabla_h v \|_{L^q(K\cup L; R^{3})}, \ \forall v\in V(\Omega),\,1\le q\le\infty,\\
\|v|_K-v|_L \|_{L^q(K\cup L)}\aleq h \| \nabla_h v \|_{L^q(K\cup L; R^{3})}, \ \forall \sigma=K|L,\, v\in V(\Omega),\,1\le q\le\infty,
\end{array}
\eF
%and its global version,
%\bFormula{jumps}
%\left\{
%\begin{array}{c}
%\Big(\sum_{\sigma=K|L\in {\cal E}_{\rm int}}\frac{\|v|_K- v|_L\|^q_{L^q(\sigma)}}{h^{q-1}}\Big)^{\frac 1q}\aleq \|\nabla_h v\|_{L^q(\Omega)},\\
%\|\widehat v\|_{Q^{1,q}(\Omega)}\aleq  \|\nabla_h v\|_{L^q(\Omega)}
%\end{array}\right\},
%\;\mbox{for all $v\in V(\Omega)$, $1\le q<\infty$},
%\eF
see  Gallouet et al. \cite[Lemma 2.2]{GaHeLa}.

Further we recall a global version of the Poincar\'e inequality on $V(\Omega)$
\bFormula{globalV}
\|v-\widehat v\|_{L^q(\Omega)}\aleq h\|\nabla_h v\|_{L^q(\Omega)}
\;\mbox{for all $v\in V(\Omega)$, $1\le q<\infty$},
\eF  
along with the global error estimate 
\bFormula{errorV}
\left\| v - \tilde v  \right\|_{L^q(\Omega)} +
h \left\| \Gradh \left( v - \tilde v \right) \right\|_{L^q(\Omega;R^3)}    \aleq  h^m \left\| v \right\|_{W^{m,q}
(\Omega)}
\eF 
for any $v \in W^{m,q} (\Omega)$, $m = 1,2,\  1 \le q \le \infty$,
see Karper \cite[Lemma 2.7]{Ka} or Crouzeix and Raviart \cite{CrRa}.

Next, we shall deal with the Sobolev properties of piecewise constant and Crouzeix-Raviart finite elements.
To this end we introduce a discrete (so called broken) Sobolev $H^{1,p}$-(semi)norm on $Q(\Omega)$,
\bFormula{?}
\| g \|_{Q^{1,p}(\Omega)}^p = \sum_{\sigma \in {\cal E}_{\rm int}} \int_{\sigma} \frac{ \ju{g}_\sigma^p }{h^{p-1}} \ {\rm dS}_x,\; 1\le p<\infty.
\eF
and a discrete (so called broken) Sobolev $H^{1,p}$-(semi)norm on $V(\Omega)$,
\bFormula{EM2}
\| {v} \|_{V^{1,p}(\Omega)}^p = \intO{  |\Gradh {v} |^p }.
\eF

 Related to the $Q^{1,p}$-norm, we report the following discrete Sobolev and Poincar\'e type inequalities: We have
\bFormula{SobolevQ}
\forall g\in Q(\Omega),\;\| g \|_{L^6(\Omega)} \aleq \| {g} \|_{Q^{1,2}(\Omega)} 
+ \| g \|_{L^2(\Omega)},
\eF
%and
%\bFormula{PoincareQ}
%\forall M,\Gamma>0,\, \gamma>1, \,\exists C>0,\,\forall g\in Q(\Omega),
%\eF
%$$
%\| g \|_{L^2(\Omega)} \leq C(M,\Gamma) \left( \| {g} \|_{Q^{1,2}(\Omega)} + \intO{ r |g| } \right)
%$$
%with any $r$,
%$$
%r\ge 0,\;
%0 < \intO{ r } = M ,\ \intO{ r^\gamma } \leq \Gamma.
%$$
see Bessemoulin-Chatard et al.  \cite[Theorem 6]{ChHiDr}. 
%for the former and
%\cite[Lemma 2.2]{FeKaNo} for the latter inequality.
Likewise, we have the Discrete Sobolev inequality,
\bFormula{SobolevV}
\forall v\in V(\Omega),\,
\left\{\begin{array}{c}
\| {v} \|_{L^{p^*}(K)} \aleq \| \Grad{v} \|_{L^p(K)}+
\|v\|_{L^p(K)},\\
\| {v} \|_{L^{p^*}(\Omega)} \aleq \| {v} \|_{V^{1,p} (\Omega)}+
\|v\|_{L^p(\Omega)}
\end{array}\right\},\; 1\le p\le\infty,
\eF
\begin{equation}\label{SobolevV0}
\forall v\in V_0(\Omega),\,\| {v} \|_{L^{p^*}(\Omega)} \aleq \| {v} \|_{V^{1,p}(\Omega)},
\end{equation}
(where $p^*$ is the Sobolev conjugate exponent to $p$), see e.g. \cite[Lemma 9.3]{GaHeMaNo}.

Last but not least, it is well known -- see
Christiansen et al. \cite[Theorem 5.67]{ChMKOw} for the very general formulation -- that the functions in $Q(\Omega)$ -- can be approximated by smooth functions. We report this result in the form formulated in Proposition 2 and formulas (8.6), (8.8) in \cite{FeKaPo}. To this end we extend
 $g\in Q_h(\Omega)$  by $0$ outside $\Omega$ calling the new functions again $g$, cf. (\ref{Q0}), and we denote
by $[g]_\ep$, $\ep>0$ the mollified $g$  via the standard mollifier (a convolution of $g$ with the
standard mollifying kernel). We have the following: For any ${\cal K}$ a compact subset
of $\Omega$, there exists $C=C(K)\aleq 1$ such that for any $0<h\aleq {\rm dist}({\cal K};R^3\setminus\Omega)$ and  $g\in Q_h(\Omega)$
\bFormula{regQ}
\|[g]_h-g\|_{L^2({\cal K})}\aleq h\|g\|_{Q^{1,2}(\Omega)},\;\|\nabla_x [g]_h\|_{L^2({\cal K})}
\aleq \|g\|_{Q^{1,2}(\Omega)}
\eF
{
Similar statement holds for functions in $v\in V_{h,0}(\Omega)$ (extended by $0$ outside $\Omega$),
cf. \cite[Formula (8.7) and (8.9)]{FeKaNo},
\bFormula{regV}
\|[v]_h-v\|_{L^2({\Omega})}\aleq h\|v\|_{V^{1,2}(\Omega)},\;\|\nabla_x [v]_h\|_{L^2(\Omega)}
\aleq \|v\|_{V^{1,2}(\Omega)}.
\eF
}

Finally, we report the well-known identities for the Crouzeix-Raviart finite elements,
\bFormula{vv1}
\intO{ \Divh \tilde{\vc{u} }\ w } = \intO{ \Div \vc{u} \ {w} } \ \mbox{for any} \ {w} \in Q(\Omega)\ \mbox{and}\,\vc{u}\in V(\Omega;R^3),
\eF
{
\bFormula{proj}
\intO{ \Gradh v \otimes \Gradh \tilde\varphi } = \intO{\Gradh v \otimes \Grad \varphi } \ \mbox{for all}\ v \in V_{h}(\Omega),\
\varphi \in W^{1,1} (\Omega),
\eF
see \cite[Lemma 2.11]{Ka}.
}

%%%%
\subsubsection{Trace and ``negative'' $L^p-L^q$ estimates for finite elements}\label{TNE}

We start by the classical trace estimate,
\bFormula{trace}
\| v \|_{L^q(\partial K)} \aleq  \frac{1}{h^{1/q}} \left( \| v \|_{L^q(K)} + h \| \Grad v \|_{L^q(K; R^3)} \right),\
1 { \leq q} \le \infty
\ \mbox{for any}\ v \in C^1({ \overline K}),
\eF
The following can be easily  obtained from the previous one by means of the scaling arguments.
\bFormula{traceFE}
\| w \|_{L^q(\partial K)} \aleq \frac{1}{h^{1/q}} \| w \|_{L^q(K)} \ \mbox{for any}\ 1 \leq q { \le} \infty,
\  w \in P_m,
\eF
where $P_m$ denotes the space of polynomials of order $m$.

In a similar way, from the local estimate
\bFormula{inter}
\| w \|_{L^p(K)} \aleq h^{3 \left( \frac{1}{p} - \frac{1}{q} \right) } \| w \|_{L^q(K)} \ 1 \leq q < p \leq  \infty, \ w \in P_m,
\eF
we deduce the global version
\bFormula{interG}
\| w \|_{L^p(\Omega)} \leq c h^{3 \left( \frac{1}{p} - \frac{1}{q} \right) } \| w \|_{L^q(\Omega)} \ 1 \leq q < p \leq  \infty,
\ \mbox{for any} \ w|_K \in P_m (K), \ K \in {\cal T}.
\eF
In particular, for a piecewise constant function in $[0,T=N\Delta t)$, $a(t)=\sum_{n=0}^{N-1} a_n 1_{I_n}(t)$, where $I_n=[n\Delta t,(n+1)\Delta t)$ we have
\bFormula{time}
\| a \|_{L^p(I_n)} \aleq (\Delta t)^{ \left( \frac{1}{p} - \frac{1}{q} \right) } \| a \|_{L^q(I_n)} \ 1 \leq q < p \leq  \infty,
\eF
and
\bFormula{timeG}
\| a \|_{L^p(0,T)} \aleq  (\Delta t)^{ \left( \frac{1}{p} - \frac{1}{q} \right) } \| a\|_{L^q(0,T)} \ 1 \leq q < p \leq  \infty.
\eF
\subsubsection{Some formulas related to upwinding}
The following formulas can be easily verified by a direct calculation, see e.g. \cite[Section 2.4]{FeKaNo}
\begin{enumerate}
 \item Local conservation of fluxes:

\begin{equation}\label{Up1}
 \begin{array}{c}
  \forall \sigma=K|L\in {\cal E}_{\rm int},\;F_{\sigma,K}[g,\vu]=-F_{\sigma,L}[g,\vu],\\
  \forall \sigma \in {\cal E}_{\rm int},\ {\rm Up}_{\sigma,\vc n}[{ g}, \vu]=-
  {\rm Up}_{\sigma,-\vc n}[{ g}, \vu]\ \mbox{and}\ [[g]]_{\vc n_\sigma}=-[[g]]_{-\vc n_\sigma}.
 \end{array}
\end{equation}
\item For all $r,g\in Q(\Omega)$, $\vu\in V(\Omega,R^3)$,
\begin{equation}\label{Up2}
 \sum_{K\in {\cal T}}r_K\sum_{\sigma\in {\cal E}_{\rm int}}|\sigma| F_{\sigma,K}[g,\vu]=
  -\sum_{\sigma\in{\cal E}_{\rm int}} {\rm Up}_{\sigma}[g,\vu][[r]]_\sigma.
  \end{equation}
\item For all $r,g  \in Q(\Omega)$, $\vu,\vu_B\in V(\Omega;R^3)$,  $\vu-\vu_B\in V_0(\Omega;R^3)$,
 $\phi \in C^1(\Ov{\Omega})$,
\bFormula{Up3}
\intO{ g \vu \cdot \Grad \phi } = -\sum_{K\in {\cal T}}\sum_{\sigma\in {\cal E}(K)\cap{\cal E}_{\rm int}} F_{\sigma,K}[g,\vu] r
%+
%\sum_{\sigma\in {\cal E}_{\rm int}} \int_\sigma {\rm Up}_\sigma[g, \vu] \ju{ \widehat\phi }_\sigma {\rm d}S_x
\eF
$$
+ \sum_{K\in{\cal T}} \sum_{\sigma\in {\cal E}(K)\cap{\cal E}_{\rm int}}\int_\sigma (r -  \phi ) \ \ju{ g }_{\sigma,\vc n_{\sigma,K}} \ [\vu_\sigma \cdot \vc{n}_{\sigma,K}]^- {\rm d}S_x 
$$
$$
+ \sum_{K \in {\cal T}} \sum_{\sigma\in {\cal E}(K)}\int_\sigma  g (\vu -
\vu_\sigma ) \cdot \vc{n}_{\sigma,K}(\phi-r){\rm d}S_x
+ \intO{ (r - \phi) g \Divh \vu }
$$
$$
+  \sum_{\sigma\in {\cal E}^{\rm out}}\int_\sigma
g\vu_{B,\sigma}\cdot \vc n_{\sigma}(\phi-r){\rm d}S_x
+  \sum_{\sigma\in {\cal E}^{\rm in}}\int_\sigma
g\vu_{B,\sigma}\cdot \vc n_{\sigma}(\phi-r){\rm d}S_x.
$$
\end{enumerate}

\subsection{Preliminaries from  analysis}

%\subsection{Some elements of functional, convex and harmonic analysis}

In the existence proof, we need several properties related to the space $[W^{1,q'}(\Omega)]^*$ We start by
the Riesz representation formula, which can be readily obtained from the Hahn-Banch theorem and the definition of dual norms:

\begin{Lemma}\label{LRiesz}
Let $1<q<\infty$. Then for any ${\cal G}\in [W^{1,q'}(\Omega)]^*$ there exists a unique $[g_0,\vc g]\in L^q(\Omega;R^{4})$ such that
\bFormula{RieszW*}
\forall u\in W^{1,q'}(\Omega),\;<{\cal G};u>=\intO{\Big(g_0 +\vc g\cdot\Grad u\Big)};
\eF
$$
\mbox{if $\|u\|_{W^{1,q'}(\Omega)}:=\|u\|_{L^{q'}(\Omega)}+\sum_{i=1}^3\|\partial_i u\|_{L^{q'}(\Omega)}$ then
$\|{\cal G}\|_{[W^{1,q'}(\Omega)]^*}=\max_{i=0,\ldots,3}\|g_i\|_{L^q(\Omega)}$}.
$$
\end{Lemma}
This formula implies immediately,
\bFormula{W*dens}
C^\infty_c(\Omega)\;\mbox{is dense in}\;[W^{1,q'}(\Omega)]^*.
\eF
Indeed, if (\ref{RieszW*})  is a representation of ${\cal G}\in [W^{1,q'}(\Omega)]^*$ and $[g^n_0,\vc g^n]\in
C^\infty_c(\Omega;R^{4})$ is a sequence converging to $[g_0,\vc g]$ in $L^q(\Omega;R^{4})$ then
 $C^\infty_c(\Omega;R^{4})\ni {\cal G}^n:= g_0^n-\Div\vc g^n\to {\cal G}$ in $[W^{1,q'}(\Omega)]^*$.

We will need the properties of the pseudodifferential operator
\bFormula{calA}
{\cal A}[u]=\Grad\Delta^{-1}[u] := \mathcal{F}^{-1}_{\xi \to x} \left[
\frac{-{\rm i}\xi}{|\xi|^2 } \mathcal{F}_{x \to \xi} [1_\Omega u] \right]\;\;u\in C^\infty_c(\Omega),\;\mbox{where
${\cal F}_{x\to\xi}$ denotes the Fourier transform}.
\eF
We have the following lemma.

\begin{Lemma}\label{LcalA}
For all $u\in C^{\infty}_c(\Omega)$ we have,
$$
\|{\cal A}[u]\|_{L^2(R^3)}\aleq \|u\|_{L^1(\Omega)}+\|u \|_{L^2(\Omega)}
$$
$$
\|{\cal A}[u]\|_{W^{1,p}(\Omega)}\aleq \|u\|_{L^p(\Omega)},\; 1<p<\infty,
$$
and for all $[g_0,\vc g]\in C^\infty_c(\Omega;R^{4})$
$$
\|{\cal A}[g_0-\Div \vc g]\|_{L^q(\Omega)}\aleq \|g_0-\Div \vc g\|_{[W^{1,q'}(\Omega)]^*},\; 1<q<\infty.
$$
\end{Lemma}

Indeed, the first estimate follows just from the splitting
$\frac 1{|\xi|}{\cal F}_{x\to\xi}[u]=\frac {1_{|\xi|\le 1}}{|\xi|}{\cal F}_{x\to\xi}[u]+
\frac {1_{|\xi|> 1}}{|\xi|}{\cal F}_{x\to\xi}[u]$, where the first term in the sum
is in $L^2(R^3)$ provided ${\cal F}_{x\to\xi}[u]\in L^\infty(R^3)$ (i.e. $u\in L^1(R^3)$), while
the second term is in $L^2(R^3)$ provided ${\cal F}_{x\to\xi}[u]\in L^2(R^3)$ (i.e. $u\in L^2(R^3)$).
 
The second estimate for  $\Grad{\cal A}[u]$ in $L^p(R^3)$ follows from the classical H\"ormander-Mikhlin multiplier theorem. According to the Sobolev embeddings in the homogenous Sobolev spaces setting,  cf. Galdi \cite[Theorem II.6.1]{Ga}, there is $a_\infty\in R$ such that, if $1<s<3$
in particular, $\|{\cal A}[u]-a_\infty\|_{L^{s^*}(R^3)}\aleq \|\Grad{\cal A}[u]\|_{L^s(R^3)}$,
$s^*=\frac{3s}{3-s}$; whence the first inequality in Lemma \ref{LcalA} yields $a_\infty=0$.
Consequently, since $s<s^*$ and $\Omega$ is bounded, the second inequality follows so far for $1<p<3$. If $p\ge 3$, we have according to the preceding ${\cal A}[u]$ bounded in $W^{1,\frac{3p}{3+p}}(\Omega)$ by $\|u\|_{L^{\frac{3p}{3+p}}(\Omega)}\aleq \|u\|_{L^{p}(\Omega)}$ and we can use Sobolev imbedding $W^{1,\frac{3p}{3+p}}(\Omega)\hookrightarrow L^p(\Omega)$ to conclude.

Finally, we deduce the third estimate in Lemma \ref{LcalA} by calculating
$$
\intO{{\cal A}_k[g_0-\Div \vc g]\phi}= \intO{{\cal A}_k[g_0]\phi}+\intO{ {\cal A}_k[\Div\vc g]\phi}
$$
$$
\aleq \|\phi\|_{L^{q'}(\Omega)}(\|{\cal A}_k[g_0]\|_{L^q(\Omega)}+\|\Div{\cal A}_k [\vc g]\|_{L^q(\Omega)})
\aleq \max_{i=0,\ldots,3}\{\|g_i\|_{L^q(\Omega)}\}\|\phi\|_{L^{q'}(\Omega)},
$$
In the above, we have used besides Lemma \ref{LRiesz} also the H\"older inequality and the already proved second estimate of Lemma \ref{LcalA}.

We can therefore extend the operator (\ref{calA}) from
$C^\infty_c(\Omega)$ functions to $L^p(\Omega)\cap
[W^{1,q'}(\Omega)]^*$. We resume this observation in the following Lemma.

\begin{Lemma}\label{LcalA+}
The operator $\cal A$ explicitly defined in (\ref{calA})
is a continuous linear operator from
$$
L^p(\Omega)\;\mbox{to $W^{1,p}(\Omega)$ and from
$[W^{1,q'}(\Omega)]^*$ to $L^q(\Omega)$}, 1<q,p<\infty. 
$$
\end{Lemma}

The next lemma deals with a particular solution of the eqation ${\rm div}\vc w=r$ in  $\Omega$,  $\vc w|_{\partial\Omega}=0$ 
called the Bogovskii solution, \cite{Bog}.

\begin{Lemma}\label{LBog}
Let $\Omega$ be a bounded Lipschitz domain. There exists a linear operator ${\cal B}$ defined on ${C}^\infty_c(\Omega)$
with the following properties
\bFormula{Z0}
{\cal B}(C^\infty_c(\Omega))\subset  C^\infty_c(\Omega)
\eF
\bFormula{Z2}
\forall r\in {C^\infty_c(\Omega)},\quad\left\| \mathcal{B}[r] \right\|_{W^{1,p}(\Omega;R^3)} \aleq \| r\|_{L^p(\Omega)},\ 1 < p < \infty,
\eF
\bFormula{Z3}
\forall {\cal G}:= g_0-\Div\vc g, \;[g_0,\vc g]\in {C^\infty_c(\Omega;R^4)},\quad
\left\| \mathcal{B}[{\cal G}] \right\|_{L^q(\Omega;R^3)} \aleq \| {\cal G} \|_{[W^{1,q'}(\Omega)]^*},\ 1 < q < \infty.
\eF
\bFormula{Z1}
\Div \mathcal {B}[r] = r \ \mbox{for any} \ r \in C^\infty_c(\Omega),\ \intO{ r } = 0.
\eF
\end{Lemma}

On Lipschitz domains, the Bogovskii solution is given by an explicit formula involiving a sigular kernel which is particularly "acessible" if the domain is star-shaped and which allows to provide the proof of Lemma \ref{LBog} via an explicit
(but involved) calculation. We refer to Galdi \cite[Chapter 3]{Ga} for a detailed proof of the properties (\ref{Z0}), (\ref{Z2}), (\ref{Z1}), and to Geissert, Heck and Hieber \cite{GeHeHi} for (\ref{Z3}). 

Using the density of $C^\infty_c(\Omega)$ in $L^q(\Omega)$ and in $[W^{1,q'}(\Omega)]^*$, cf. Lemma \ref{RieszW*} and formula (\ref{W*dens}), we can extend the operator $\cal B$ introduced in  Lemma \ref{LBog} to a continuous linear operator
on $L^q(\Omega)\cap[W^{1,q'}(\Omega)]^*$:

\begin{Lemma}\label{LBog+}
The operator ${\cal B}$ introduced in Lemma \ref{LBog} is a continuous linear operator
$$
\mbox{from ${L^p(\Omega)}$ to $W^{1,p}_0(\Omega)$ and from $[W^{1,q'}(\Omega)]^*$ to $L^q(\Omega)$, $1<q,p<\infty$}
$$
and
$$
\Div \mathcal {B}[r] = r \ \mbox{for any} \ r \in L^p(\Omega),\ \intO{ r } = 0,\; 1<p<\infty.
$$
\end{Lemma}

The next theorem involving commutator of Riesz operators may be seen as a consequence
of the celebrated Div-Curl lemma above, see Murat, Tartar \cite{Mu} and
\cite[Section 6]{Fe} or \cite[Theorem 10.27]{FeNoB} for the below adapted formulation

\begin{Lemma}\label{rieszcom} Let
\[
\vc{ V}_n \rightharpoonup \vc{ V} \ \mbox{ in}\ L^p(\R^3; \R^3),
\]
\[
\vc{ U}_n \rightharpoonup \vc{ U} \ \mbox{ in}\ L^q(\R^N; \R^N),
\]
where $ \frac{1}{p} + \frac{1}{q} = \frac{1}{s} < 1$. Then
\[
\vc{ U}_n \cdot \Grad\Delta^{-1}\Div[\vc{V}_n ] -\vc V_n\cdot \Grad\Delta^{-1}\Div
[\vc{
U}_n]\cdot \vc{ V}_n \rightharpoonup \vc{ U}\cdot \Grad\Delta^{-1}\Div[\vc{ V}] -
\vc V\cdot\Grad\Delta^{-1}\Div[\vc{ U}] \ \mbox{ in}\ L^s(\R^N).
\]
\end{Lemma}

Finally, the last two lemmas are well known results from convex analysis, see e.g. Lemma 2.11 and Corollary 2.2 in
Feireisl \cite{Fe}.

\begin{Lemma}\label{Lemma2}
Let $O \subset R^d$, $d\ge 2$, be a measurable set and $\{ \vc{v}_n
\}_{n=1}^{\infty}$ a sequence of functions in $L^1(O; R^M)$ such
that
\[
\vc{v}_n \rightharpoonup \vc{v} \ \mbox{ in}\ L^1(O; R^M).
\]
Let $\Phi: R^M \to (-\infty, \infty]$ be a lower semi-continuous
convex function.

Then $\Phi(\vc v):O\mapsto R$ is integrable and
$$
\int_{O} \Phi(\vc v){\rm d} x\le \liminf_{n\to\infty} \int_{O} \Phi(\vc v_n){\rm d} x.
$$
\end{Lemma}

\begin{Lemma}\label{Lemma3}

Let $O \subset R^d$, $d\ge 2$ be a measurable set and $\{ \vc{v}_n
\}_{n=1}^{\infty}$ a sequence of functions in $L^1(O; R^M)$ such
that
\[
\vc{v}_n \to \vc{v} \ \mbox{weakly in}\ L^1(O; R^M).
\]
Let $\Phi: \R^M \to (-\infty, \infty]$ be a lower semi-continuous
convex function such that $\Phi (\vc{v}_n ) \in L^1(O)$ for
any $n$, and
\[
\Phi (\vc{v}_n) \to \Ov{\Phi(\vc{v})}
\ \mbox{weakly in}\ L^1(O).
\]

Then

\bFormula{FNeb4}
\Phi (\vc{v})  \leq \Ov{\Phi (\vc{v})} \ \mbox{a.e. in}\ O.
\eF

If, moreover, $\Phi$ is strictly convex on an open convex set
$U \subset R^M $, and
\[
\Phi(\vc{v}) = \Ov{\Phi (\vc{v})} \ \mbox{a.e. on}\ O,
\]
then

\bFormula{FNeb5}
\vc{v}_n  (\vc{y}) \to \vc{v} (\vc{y})
\ \mbox{for a.a.}\ \vc{y} \in \{
\vc{y} \in O \ | \ \vc{v}(\vc{y}) \in U \}
\eF
extracting a subsequence as the case may be.
\end{Lemma}

\begin{Lemma} \label{Lemma4}
Let $O$ be a domain in $R^d$, $P,G:Q\mapsto R$ be a couple of non decreasing and continuous functiona on $[0,\infty)$.  Assume that  $\vr_n\in
L^1(O;[0,\infty))$ is a sequence such that
$$
\left.\begin{array}{c}
P(\vr_n) \rightharpoonup \overline{P(\vr)}, \\
G(\vr_n) \rightharpoonup \overline{G(\vr)}, \\
P(\vr_n)G(\vr_n) \rightharpoonup \overline{P(\vr)G(\vr)}
\end{array} \right\} \mbox{ in } L^1(O).
$$
Then
$$
\overline{P(\vr)}\, \, \overline{G(\vr)} \leq
\overline{P(\cdot,\vr)G(\cdot,\vr)}
$$
a.e. in $O$.
\end{Lemma}

\end{document}